\documentclass[a4paper,onecolumn,superscriptaddress,10pt,accepted=2026-10-17]{compositionalityarticle}
\pdfoutput=1

\usepackage[utf8]{inputenc}
\usepackage[english]{babel}
\usepackage[T1]{fontenc}

\newcommand{\authorsforheader}{Carranza, Kapulkin, Kershaw, Laubenbacher and Wheeler}
\newcommand{\compositionalityIssue}{6}
\newcommand{\compositionalityVolume}{8}
\newcommand{\paperdoinumber}{10.46298/compositionality-\compositionalityVolume-\compositionalityIssue}
\newcommand{\paperdoi}{https://doi.org/\paperdoinumber}
\newcommand{\receiveddate}{2025-08-13}
\newcommand{\accepteddate}{2026-10-17}
\usepackage{amsmath}         % align-like environments
\usepackage{amsthm,thmtools,thm-restate} % theorem-like environments
\usepackage{enumitem}        % enumerate-like environments

\usepackage{graphicx,graphics}
\graphicspath{{figures/}}

\usepackage{tikz}
  \usetikzlibrary{matrix,arrows,positioning,decorations.markings}
  
\usepackage{tikz-cd}

\usepackage{float} %to force figure into right place
\usepackage{caption, subcaption} % for subfigures
\usepackage{wrapfig} % to wrap figures around text

\usepackage{hyperref}                              % hyperlinks
\usepackage[capitalize,nosort]{cleveref}           % named references

\usepackage{amssymb}
\usepackage[mathscr]{euscript}
\usepackage{relsize}
\usepackage{stmaryrd}
\usepackage{stackengine}
\usepackage{mathdots} % for inverse diagonal dots

\usepackage{indentfirst} % indent the first paragraph
\usepackage{multicol}

\usepackage{booktabs}

\tikzstyle{vertex}=[circle, fill, minimum size=4pt, inner sep=0pt, outer sep=4pt]
\tikzstyle{colorvertex}=[circle, draw, minimum size=6pt, inner sep=0pt]
\definecolor{vertex1}{RGB}{255,0,0}
\definecolor{vertex2}{RGB}{0,0,255}
\definecolor{vertex3}{RGB}{0,255,0}
\definecolor{vertex4}{RGB}{255,255,0}
\definecolor{vertex5}{RGB}{255,0,255}

\newcommand{\from}{\colon}
\newcommand{\ito}{\hookrightarrow}

\declaretheorem[style=definition,numberwithin=section]{definition}
\declaretheorem[style=definition,numberlike=definition]{example}
\declaretheorem[style=definition,numberlike=definition]{remark}

\declaretheorem[style=plain,numberlike=definition]{corollary}
\declaretheorem[style=plain,numberlike=definition]{lemma}
\declaretheorem[style=plain,numberlike=definition]{proposition}
\declaretheorem[style=plain,numberlike=definition]{theorem}

\Crefname{corollary}{Corollary}{Corollaries}
\Crefname{definition}{Definition}{Definitions}
\Crefname{lemma}{Lemma}{Lemmas}
\Crefname{proposition}{Proposition}{Propositions}
\Crefname{remark}{Remark}{Remarks}
\Crefname{theorem}{Theorem}{Theorems}
\Crefname{notation}{Notation}{Notations}
\crefname{paragraph}{Paragraph}{Paragraphs}

\declaretheorem[style=definition,numbered=no,name=Definition]{definition*}
\declaretheorem[style=definition,numbered=no,name=Example]{example*}
\declaretheorem[style=definition,numbered=no,name=Remark]{remark*}
\declaretheorem[style=definition,numbered=no,name=Notation]{notation*}

\declaretheorem[style=plain,numbered=no,name=Corollary]{corollary*}
\declaretheorem[style=plain,numbered=no,name=Lemma]{lemma*}
\declaretheorem[style=plain,numbered=no,name=Proposition]{proposition*}
\declaretheorem[style=plain,numbered=no,name=Theorem]{theorem*}

\newcommand{\Ftwo}{\mathbb{F}_2} % the field with 2 elements

\newcommand{\DiGraph}{\mathsf{DiGraph}} % The category of digraphs
\newcommand{\MultiDiGraph}{\mathsf{MultiDiGraph}} % the category of multidigraphs
\newcommand{\Set}{\mathsf{Set}} % The category of sets
\newcommand{\DDS}{\mathsf{DDS}} % The category of discrete dynamical systems
\newcommand{\cySet}{\mathsf{cySet}} % The category of cycle sets

\newcommand{\cat}[1]{\mathsf{#1}} % font for generic categories
\newcommand{\Vect}{\mathsf{Vect}} % The category of vector spaces
\newcommand{\NN}{\mathbb{N}} % the natural numbers
\newcommand{\ZZ}{\mathbb{Z}}
\newcommand{\BN}{B\NN} % The monoid of natural numbers, as a one-object category
\newcommand{\BNobj}{\ast} % The unique object of the category \BN
\newcommand{\Cy}{\mathbf{C}_\bullet} % The category of cycles
\newcommand{\GG}{\mathbf{G}} % The small category with objects {V, E}
\newcommand{\BZ}[1]{B(\ZZ / #1)} % the cyclic group Z/n as a one-object category
\newcommand{\ZnSet}[1]{\Set^{\BZ{#1}}}
\newcommand{\Cyhat}[1]{\widehat{C}_{#1}} % the cycle set represented by C_#1

\renewcommand{\succ}{\mathrm{succ}} % the successor operation
\newcommand{\Gfun}{Q} % the reflection functor from multidigraphs to digraphs
\newcommand{\modmap}[2]{\mu_{#1, #2}} % the map C_m -> C_n given by mod
\newcommand{\autmap}[2]{\rho_{#1, #2}}
\newcommand{\semiprod}[1]{\mathbin{\rtimes_{#1}}} % semi-direct product symbol - argument for subscript
\newcommand{\Orb}{\operatorname{Orb}} % set of orbits of a G-set

\newcommand{\DEdge}{E^\circ} % notation for the directed edge graph

\newcommand{\op}{\mathsf{op}}
\newcommand{\slice}{\downarrow}
\DeclareMathOperator*{\colim}{colim}
\newcommand{\adj}{\dashv} % adjoint symbol
\newcommand{\pbtick}{\lrcorner}

\newcommand{\id}[1][]{\mathrm{id}_{#1}} % identity morphism
\newcommand{\proj}{\operatorname{proj}}
\newcommand{\FF}{\mathbb{F}}
\newcommand{\supp}{\operatorname{supp}} % the support of an object in a regular category
\newcommand{\Aut}{\operatorname{Aut}} % automorphism group
\newcommand{\Orbnd}{\operatorname{Orb}_{\mathrm{nd}}} % non-degenerate orbits
\begin{document}

\title{Categorical foundations of discrete dynamical systems}

% Alphabetical order (by surname)
\author{Daniel Carranza}
\affiliation{Department of Mathematics, Johns Hopkins University, Baltimore, Maryland, United States}
\author{Krzysztof Kapulkin}
\email{k.kapulkin@vanderbilt.edu}
\affiliation{College of Connected Computing \& Department of Mathematics, Vanderbilt University, Nashville, Tennessee, United States}
\author{Nathan Kershaw}
\affiliation{Department of Mathematics, University of Western Ontario, London, Ontario, Canada}
\author{Reinhard Laubenbacher}
\affiliation{Department of Medicine, University of Florida, Gainesville, Florida, United States}
\author{Matthew Wheeler}
\affiliation{Department of Medicine, University of Florida, Gainesville, Florida, United States}

\begin{abstract}
  We develop categorical foundations of discrete dynamical systems, aimed at understanding how the structure of the system affects its dynamics.
  We introduce the notion of cycle sets to analyze attractors of a system, and use this to generalize multiple decomposition theorems of Kadelka, Veliz-Cuba, Murrugarra, and the last two authors from Boolean networks to arbitrary discrete dynamical systems.
\end{abstract}

\maketitle

\section*{Introduction}

\subsection{Motivation}

Discrete dynamical systems are a common generalization of many mathematical modeling frameworks from the natural sciences, including: Boolean networks, used in biology to represent gene regulatory networks; Petri nets, a common modelling framework from biology and biochemistry; and cellular automata, used as models in large-scale brain networks in neuroscience.
A discrete dynamical system consists of a (not necessarily finite) set $X$ along with a function $f \colon X \to X$.
In the case of a Boolean network, the set $X$ is taken to be $\mathbb{F}_2^n$, for a Petri net, it is $P^\mathbb{N}$ for some set $P$, while for cellular automata, one takes $X$ to be $A^G$ for a set $A$ and a group $G$.

Given a dynamical system $(X, f)$, we typically want to analyze all possible trajectories of elements of $X$. 
That is, for a given $x \in X$, we wish to understand the sequence of elements $(x, f(x), f^2(x), \ldots)$.
This is typically done using the following ``pipeline:''
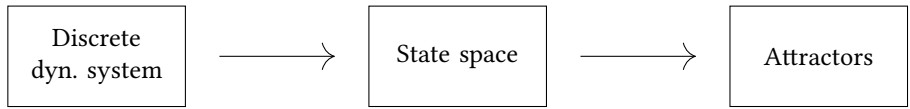
\begin{center}
    \begin{tikzpicture}[box/.style={draw, text width=5em, align=center, outer sep=3ex, inner sep=2ex, minimum height=3.8em}]
        \node[box] (A) {Discrete dyn.\ system};
        \path (A) -- +(4.8, 0) node[box] (B) {State space};
        \path (B) -- +(4.8, 0) node[box] (C) {Attractors};

        \draw[-{>[scale=2.0]}] (A) -- (B);
        \draw[-{>[scale=2.0]}] (B) -- (C);
    \end{tikzpicture}
\end{center}

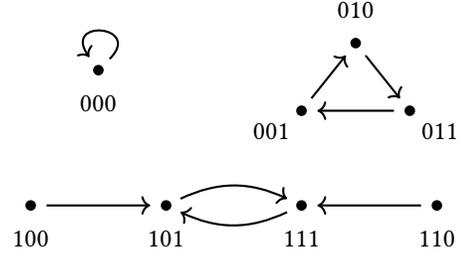
\begin{wrapfigure}{R}{0.42\textwidth}
\centering

\begin{tikzpicture}[
    edge/.style={->,thick}, % shorten keeps arrowheads off the vertices
    every label/.style={black},
    scale=0.9
]

% ------------- vertices -------------
\node[vertex,label={below:$000$}] (000) at (1,0)      {};
\node[vertex,label={[label distance=-1.1ex]below left:$001$}] (001) at (4,-0.6)      {};
\node[vertex,label={above:$010$}] (010) at (4.8, 0.4)      {};
\node[vertex,label={[label distance=-1.1ex]below right:$011$}] (011) at (5.6,-0.6)  {};
\node[vertex,label={below:$100$}] (100) at (0,-2)     {};
\node[vertex,label={below:$101$}] (101) at (2,-2)     {};
\node[vertex,label={below:$111$}] (111) at (4,-2)     {};
\node[vertex,label={below:$110$}] (110) at (6,-2)     {};

% ------------- edges -------------
\path[edge]
  % round, tidy self-loop
  (000) edge[out=45, in=125, looseness=6] (000)
  % 3-cycle
  (001) edge (010)
  (010) edge (011)
  (011) edge (001)
  % tail into 2-cycle
  (100) edge (101)
  % 2-cycle
  (101) edge[bend left=25] (111)
  (111) edge[bend left=25] (101)
  % exit arrow
  (110) edge (111);
\end{tikzpicture}
\caption{The state space graph of a discrete dynamical system (in fact, a Boolean network) $f\colon \{0,1\}^3 \to \{0,1\}^3$ given by $f(x_1,x_2,x_3) = (x_1,x_2 \oplus x_3, x_1 \vee x_2)$.} \label{system-to-state-space}
\end{wrapfigure}

Given a dynamical system, one associates to it a directed graph (from here on, simply a digraph), called its \emph{state space graph}.
Its vertices are the elements of $X$ with edges connecting each $x \in X$ to $f(x)$ (see \cref{system-to-state-space}).
In a way, the state space graph carries the exact same information as the dynamical system itself, but this information is more readily available for analysis, as the attractors of a dynamical system are simply the cycles of its state space graph (in the sense of graph theory).

Instances of such considerations can be found across all examples of discrete dynamical systems: Boolean networks \cite{chaves2018analysis,veliz2014steady,gan-albert}, Petri nets \cite{gamache2025algorithm}, and cellular automata \cite{liu2024application}.
Each of these works develops methods specific to its particular framework.
Indeed, Boolean networks, Petri nets, and cellular automata have been, generally speaking, analyzed one-by-one, rather than en masse, and thus the tools developed for this purpose do not readily generalize to other contexts.

The goal of this paper is to change this paradigm by developing a foundational framework that permits an en-masse analysis.
The cornerstone of this framework is the language of category theory, a branch of mathematics well suited for understanding abstract structures.
As a result, our analysis focuses on the structural aspects of a system and thus establishes a clear connection between the structure of a system and its dynamics.

The key motivation behind it, coming from biology, is the program of understanding the notion of modularity, widely accepted as a key feature of biological systems \cite{huitzil:modularity,lorenz:modularity,melo:modularity,wagner:modularity}.
Although there is no agreed-upon definition of modularity, several were proposed, e.g., \cite{newman:modularity,kadelka2023modularity}.
The latter in particular can be used for decomposition of the structure of a dynamical system.
Building on it, we prove a generalization of the main theorem from \cite{kadelka2023modularity}.

We should say at the outset that the idea of leveraging category theory to study discrete dynamical system at a higher level of generality (for instance, identifying discrete dynamical systems as certain kinds of functors $\BN \to \Set$) is not new.
It is a running example in Lawvere and Schanuel's \emph{Conceptual Mathematics} \cite{lawvere-schanuel}, where systems of this kind appear under the name of \emph{endomaps}, and it appears again in Lawvere and Rosebrugh's \emph{Sets for Mathematics} \cite{lawvere-rosebrugh}; it is also one of the six guiding examples of La Palme Reyes, Reyes, and Zolfaghari's \emph{Generic figures and their glueings} \cite{la-palme-reyes-reyes-zolfaghari}; some of the basic adjunctions we recall in \cref{sec:framework} can already be found in those texts.
Our aim is not to introduce this point of view, but to develop it further, to the point where it becomes a working tool for the analysis of attractors, e.g. we introduce new mathematical objects to study these systems (cycle set) and use these notions to express new decomposition results for understanding attractors of a system.

\subsection{Organization and language of the paper}

The main body of this paper, contained in \cref{sec:framework,Section:Decomposition,sec:wiring,sec:conclusion}, is written in the language of category theory and requires a working knowledge of the subject at the level of a standard text, with \cite{riehl:context} being the canonical choice.
Within each section, we aim to showcase three contributions.

% \begin{enumerate}
    % \item 
    \cref{sec:framework} contributes the notion of a \emph{cycle set} (\cref{cycle-set}) and the \emph{attractor functor} (\cref{attractor-functor}) which assigns to a system the cycle set of its attractors, and establishes the properties of this functor that we need, culminating in a recognition theorem (\cref{characterization-digraph-attractors}) characterizing exactly which cycle sets arise as the cycle sets of attractors of digraphs.
    % \item 
    
    Using these notions, \cref{Section:Decomposition} contributes the notion of semi-direct products of dynamical systems (\cref{semidirect-product}) and gives a decomposition theorem (\cref{modularity-theorem}), which says that the attractor functor carries a semi-direct product decomposition of a system to a decomposition of its attractors. This generalizes the main theorem of \cite{kadelka2023modularity} from Boolean networks to arbitrary discrete dynamical systems: the theorem becomes an instance of the general fact that right adjoints preserve limits.
    % \item 
    
    \cref{sec:wiring} makes the decomposition theorem effective, contributing recognition theorems for when semi-direct product decompositions are known to exist. The key result (\cref{semiprod-iff-wiring-map}) shows that a system $f \from A^n \to A^n$ admits a semi-direct product decomposition if and only if its \emph{wiring diagram} $W(f)$ admits a graph map to a fixed two-vertex digraph. This is a condition on a graph with $n$ vertices, whereas the state space has $\lvert A \rvert^n$ vertices.
% \end{enumerate}

The first and third contributions both exist to support the second: the first supplies the language in which the decomposition theorem can be stated, while the third supplies the means of applying it.
\cref{sec:conclusion} contains the paper's conclusion and a discussion of future work.

We hope that the results of this paper can be of interest also to those without background in category theory.
For their benefit, in the remainder of this introduction, we summarize our results, making minimal use of the language of category theory.
An excellent reference for category theory aimed at researchers from other sciences is \cite{spivak:category-theory}.

\subsection{Summary of results}

\paragraph{The category of discrete dynamical systems.}

The benefit of our general approach is that we can take advantage of the framework of category theory, which provides a systematic language in which one can make sense of mathematical structures.
A category consists of a collection of objects, seen as the basic structures of a mathematical theory, and a collection of morphisms, seen as the transformations between these structures.
Thus, for instance, we have the category of sets and functions, the category of groups and group homomorphisms, the category of topological spaces and continuous functions, and the category of graphs and graph maps.

The advantage of this language is that many constructions across mathematics are instances of some definition stated in an abstract category.
For example, the cartesian product of sets, the Tychonoff topology on a product of topological spaces, and the Kronecker product of graphs are instances of a general definition of a product in an arbitrary category, instantiated to the categories of sets, topological spaces, and graphs, respectively.
By employing the framework of category theory to study discrete dynamical systems, we can canonically understand the underlying structure, instead of relying on ad hoc approaches.

The category of discrete dynamical systems can be defined as the category of functors $\BN \to \Set$, denoted $\DDS := \Set^{\BN}$.
Here $\BN$ is the monoid of natural numbers, meaning it has a unique object $\ast$ and a morphism from $\ast$ to itself for every $n \in \mathbb{N}$.
Unwinding this definition, a discrete dynamical system is a pair $(X, f\colon X \to X)$.
A morphism, or map, between discrete dynamical systems $(X,f)$ and $(Y,g)$ is a function $\alpha \colon X \to Y$ such that $\alpha \circ f = g \circ \alpha$, i.e., making the following diagram commute:
\[ \begin{tikzcd}
    X \ar[d, "f"'] \ar[r, "\alpha"] & Y \ar[d, "g"] \\
    X \ar[r, "\alpha"] & Y\text{.}
\end{tikzcd} \]
One benefit of defining the category this way is we immediately get that it is both complete and cocomplete, meaning it has all limits and colimits. 
Broadly speaking, this allows us to perform certain constructions of dynamical systems, such as products and quotients.

To study the properties of a discrete dynamical system, we can associate to it a directed graph, called its \emph{state space graph}.
Explicitly, this graph has vertex set $X$ and a directed edge $x \to f(x)$ for all $x \in X$.

The directed graphs in question also form a category, which we will denote by $\DiGraph$.
By a \emph{digraph}, we mean a set $V$ of vertices together with a subset $E \subseteq V \times V$ of edges; that is, digraphs for us are simple and directed, with at most one edge from a given vertex to a given vertex, with self-loops allowed.
A \emph{graph map} $G \to H$ is a function $V(G) \to V(H)$ on vertices that carries edges to edges.
Digraphs and graph maps form the category $\DiGraph$.
Like $\DDS$, it can be described as a category of presheaves --- on the category with two objects $V$ and $E$ and two morphisms $V \rightrightarrows E$ --- except that one must cut down to those presheaves for which the map $E \to V \times V$ is injective; we spell this out in \cref{sec:framework}, and explain this design choice in \cref{rmk:simple-digraphs}.

\begin{remark} \label{rmk:state-space-terminology}
    The term ``state space'' often refers to the underlying set (or space) $X$ on which the dynamics takes place; in particular, the same underlying space may carry many different dynamics.
    What we call the \emph{state space graph} of $(X, f)$ --- also known in the literature as the \emph{transition graph} or the \emph{functional graph} of $f$ --- is the set $X$ \emph{together with} the edges recording $f$, so that every discrete dynamical system has a unique state space graph.
    We will consistently say ``state space graph'' when we mean the digraph.
\end{remark}

The passage from a system to its state space graph is functorial: it is given by the state space functor $S \colon \DDS\to \DiGraph$.
We show in \cref{cor:S_preserves_lim_colim} that this functor preserves both limits and colimits, meaning it preserves many of the constructions we can perform in $\DDS$.
Furthermore, this functor is full and faithful (\cref{state-space-functor-full-faithful}), and in fact identifies $\DDS$ with the full subcategory of $\DiGraph$ spanned by those digraphs in which every vertex has exactly one outgoing edge (\cref{state-space-full-subcat}).
In other words, no information is lost or gained in passing from a system to its state space graph, and the digraphs that arise this way can be recognized by a simple local condition.

Below, we discuss four classes of examples of discrete dynamical systems.

\begin{example}[Boolean networks]
    A \textit{Boolean} network is a collection of Boolean variables, along with an update function for each variable that is dependent on the states of the entire collection.
    Equivalently, a Boolean network is a function $f=(f_1,\dots,f_n)\colon \mathbb{F}_2^n \to \mathbb{F}_2^n$, and hence fits into our framework.
    Boolean networks were introduced by Kauffman \cite{kauffman:boolean-networks} as models of gene regulatory networks.
\end{example}

\begin{example}[Petri nets] \label{ex:petri-nets}
    A \emph{Petri net} (see, e.g., \cite{murata:petri-nets}) is a quadruple $(P, T, F, W)$ where $P$ and $T$ are finite sets (of \emph{places} and \emph{transitions}, respectively), $F$ is a subset of $(P \times T) \cup (T \times P)$, and $W \colon F \to \mathbb{N}$ is a \emph{weight} function.
    The subset $F$ can be thought of as the set of arrows indicating the inputs and outputs (from $P$) of each transition $t \in T$, and the weight function records how many tokens each such arrow consumes or produces.
    A \emph{marking} is a function $M \colon P \to \mathbb{N}$, recording how many tokens sit at each place.
    A transition $t$ is \emph{enabled} at $M$ if $W(p, t) \leq M(p)$ for all $p$ with $(p,t) \in F$, and firing it produces the new marking
    \[ M'(p) = M(p) - W(p, t) + W(t, p), \]
    with the convention that $W(p,t) = 0$ when $(p,t) \notin F$, and similarly for $W(t,p)$.

    As stated, this does \emph{not} yet define an object of $\DDS$: at a given marking, several transitions may be enabled at once, so that a marking has no well-defined ``next'' marking, whereas an object of $\DDS$ is by definition deterministic.
    What is needed is a \emph{firing policy} resolving this choice.
    Two standard ones are: fix a total order $t_1 < \dots < t_r$ on $T$ and fire the enabled transitions in that order, giving the endofunction $f = f_{t_r} \circ \dots \circ f_{t_1}$ of the set $\mathbb{N}^P$ of markings, where $f_t$ fires $t$ if it is enabled and is the identity otherwise; or fire, at each marking, the least enabled transition in a fixed order.
    Either policy makes $(\mathbb{N}^P, f)$ an object of $\DDS$
    (a non-deterministic net is naturally modelled by a multidigraph rather than a digraph; see the discussion in \cref{sec:conclusion}).
\end{example}

\begin{example}[Cellular automata]
    A \textit{cellular automaton} (see, e.g., \cite{ceccherini-silberstein-coornaert}) consists of a set $A$, a group $G$, and a map $\tau: A^G \to A^G$ such that there exists a finite set $S$ and a function $\mu\colon A^S \to A$ so that $\tau(x)(g) = \mu((g \cdot x)|_S)$ for all $x \in A^G$ and $g \in G$. 
    A famous example is Conway's game of life \cite{games1970fantastic}, which as a cellular automaton can be interpreted as a function $\tau \colon \{0,1\}^{\mathbb{Z}^2} \to \{0,1\}^{\mathbb{Z}^2}$.
    Here $A = \{0, 1\}$ (dead or alive), $G = \ZZ^2$, and $S = \{-1, 0, 1\}^2$ is the \emph{Moore neighbourhood} of a cell, consisting of the cell itself and its eight neighbours; the local rule $\mu \colon \{0,1\}^S \to \{0,1\}$ sends $y$ to $1$ precisely when either $y(0,0) = 1$ and exactly two or three of the eight neighbours are alive, or $y(0,0) = 0$ and exactly three of them are alive.
    As cellular automata are in particular endofunctions on a set, namely $A^G$, they are also objects in the category of discrete dynamical systems. 
\end{example}

\begin{example}[Continuous dynamical systems]
More generally, a (not necessarily discrete) dynamical system (see, e.g., \cite{katok-hasselblatt}) is given by a map $\varphi \colon T\times X\to X$ for a monoid $T$ and a set $X$ such that $\varphi(s+t,x)=\varphi(s,\varphi(t,x))$, where $T$ usually represents time and $X$ represents the space of all states of the system.
In many cases, $T$ is considered to be a continuous variable with $T=\mathbb{R}$ being a usual choice.
Given a monoid $T$, we can look at all submonoids $S\subseteq T$ that have a single generator $e$.
For each such submonoid, we have a discrete dynamical system $\varphi_e \colon X\to X$ given by $\varphi_e(x)=\varphi(e,x).$ Intuitively, we could consider a system of differential equations whose solution is given by the dynamical system $\Phi \colon \mathbb{R}\times X\to X$.
Given a discretization of time, say by using microsecond intervals, we have a discrete dynamical system $\Phi_\mu \colon X\to X$ where $\Phi_\mu(x)$ outputs the state occurring one microsecond following state $x$.
Many differential equations must be discretized and simulated by computational means, and such discretized systems fit into our $\DDS$ framework. 

\end{example}

\paragraph{Cycle sets and attractors.}

A large part of the present work relies on the tool of \emph{cycle sets}.
A cycle set $K$ is a collection of sets $K_n$ for all $n \geq 1$, which one thinks of as the \emph{$n$-cycles} of some mathematical object of interest.
The primary example would be the cycles of a digraph.
For digraphs coming from the state space of a discrete dynamical systems, these cycles represent (finite) attractors of the system, which can be either fixed points or cycles.

% A morphism of cycle sets is a collection of functions
When developing the theory of cycle sets, our priority is to establish that these objects form a well-behaved category (denoted by $\cySet$), so that results from category theory instantiate to meaningful or practical structure results about discrete dynamical systems.
This emphasis requires us to insert structure that may seem unintuitive at first glance.
For instance, one consequence of our definition of cycle set is that every 2-cycle can also be viewed as a 4-cycle, a 6-cycle, an 8-cycle, etc.
We refer to such cycles as \emph{degenerate} cycles.
Similarly, every set $K_n$ of $n$-cycles comes with an action of the group $\ZZ / n$, which encodes the process of ``rotating'' a cycle.
In the example of taking the cycles of a digraph, this means that each cycle is stored with a chosen order of vertices, so that a $3$-cycle $(v_1, v_2, v_3)$ determines 3 different elements of $K_3$ (including $(v_2, v_3, v_1)$ and $(v_3, v_1, v_2)$).
Moreover, it determines a degenerate 6-cycle $(v_1, v_2, v_3, v_1, v_2, v_3)$, a degenerate $9$-cycle, and in general, a degenerate $3k$-cycle by repeating $k$ times.

The benefit of working with such an object is that the category of cycle sets $\cySet$ can be described succinctly as a presheaf category, so that for instance, it has limits and colimits which are preserved by the functor sending a cycle set $K$ to its set of $n$-cycles $K_n$ (\cref{ev-preserves-limits}).
We refer to the functor assigning a cycle set to a digraph as the \emph{attractor functor}, and denote it by $A \from \DiGraph \to \cySet$.
The definition of cycle sets as functors allows us to easily deduce that the attractor functor preserves both limits (\cref{attractor-functors-right-adj-preserves-lims}) and coproducts (\cref{attractor-functor-preserves-coprods}),
which moreover implies that the passage from discrete dynamical systems to cycle sets (via first taking the state space $S$ then taking the attractor functor $A$) also preserves limits.
Importantly, these results would fail without our conventions that cycles in a digraph are ordered, and that an $n$-cycle produces an infinite number of degenerate cycles.

We conclude \cref{sec:framework} by resolving the question of ``reversing'' the pipeline of digraphs to cycle sets: given a cycle set $K$, does there exist a digraph whose cycles recover the starting cycle set?
In general, this is not true.
For instance, in an arbitrary cycle set, it can happen that two distinct elements of $K_n$ can be ``degenerated'' to the same cycle in $K_{2n}$.
It can also happen that an $n$-cycle in $K_n$ has a non-trivial stabilizer with respect to the $\ZZ / n$-action (meaning it is similar to a degenerate cycle in that certain rotations send the cycle to itself), but the cycle itself is non-degenerate.
% However unlike with cycle sets arising from a digraph (and thus from a network in $\DDS)$, the maps $K(\mu_{m,n})$ do not have to be injective nor is it necessary for $nk$-cycle in $K_{nk}$ with periodicity $n$ required to have a corresponding $n$-cycle in $K_n$ as illustrated in Examples \ref{ex:A-not-B}-\ref{ex:not-AB}. 
These cases are somewhat pathological in that they cannot occur for cycle sets coming from a digraph, and our recognition theorem (\cref{characterization-digraph-attractors}) establishes that these are the \emph{only} possible obstructions.
That is, if neither of the two situations presented above occurs in a cycle set (formalized in \cref{def:property-AB} as Property A and Property B, respectively) then one can always construct a digraph whose cycles match the starting cycle set up to isomorphism (i.e.\ up to a relabeling of the elements of $K_n$).
% equivalently that $K$ can be written as the union of representables indexed by the non-degenerate cycles. 

\paragraph{Decomposition theorem.}

% Discrete dynamical systems are often used to model physical systems \cite{larkin:diagram, biologist:radio}. 
% These systems often come with unique structure that we can take advantage of when trying to analyze the attractors.
% One way to take advantage of this is to decompose the system into 
It is beneficial to study the attractors of a discrete dynamical system by decomposing the system into smaller, more manageable subsystems.
The reason this matters in practice is worth spelling out, since it is what motivated the theorem below in the first place.
A Boolean network on $n$ variables has a state space graph with $2^n$ vertices, so that finding its attractors by exhaustive search is hopeless for the networks arising in systems biology, where $n$ is routinely in the dozens or hundreds.
On the other hand, published gene-regulatory network models may admit some kind of modularity. 
For instance, their wiring diagrams may break up into strongly connected components, with the components arranged in a directed acyclic pattern of influence, reflecting the biological expectation that cellular processes are organized into modules \cite{huitzil:modularity,lorenz:modularity,melo:modularity,wagner:modularity}.
If the attractors of the whole network can be assembled from the attractors of these much smaller modules, then a computation of size $2^n$ is replaced by several computations of size $2^{n_i}$ with $\sum n_i = n$, which is the difference between infeasible and routine.
Reductions of this kind are developed in \cite{kadelka2023modularity} for Boolean networks, which introduces the notion of a semi-direct product of Boolean networks and proves the following theorem:
\begin{theorem}[\cite{kadelka2023modularity}]\label{old-decomposition-theorem}
    If a Boolean network $F$ can be written as a semi-direct product $F=F_1\rtimes_p F_2$, then \[A(F)=\bigsqcup_{c \in A(F_1)} c \oplus A(F_2^{c})\]
\end{theorem}
Here, $A(F)$ denotes the attractors of the network $F$, and $F_2^{c}$ is a type of twisted network which encodes $F_2$ as a new network being driven by the attractor $c$.
Loosely speaking, the theorem can be rephrased as an equality 

\[A(F\rtimes_p G)= A(F)\rtimes_p A(G). \]

In other words, the semi-direct product structure is preserved when sending a discrete dynamical system to its attractors. 
Indeed, making the right-hand side of this equation formal leads to the equality stated previously.

We use our categorical framework to prove a generalization of this result for all discrete dynamical systems. The first step is to give a categorical definition of a semi-direct product.
Our definition works in any category $\cat{C}$ with finite products: given morphisms and objects $f\colon X \to X$, $g \colon E \times Y \to Y$, and $p\colon X \to E$, the \textit{semi-direct product of $f$ and $g$ along $p$} is given by the composite:
\[ \begin{tikzcd}
        X \times Y \ar[rr, dotted, "f \semiprod{p} g"] \ar[rd, "{(\proj_X, p \circ \proj_X, \proj_Y)}"'] & {} & X \times Y \\
        {} & X \times E \times Y \ar[ur, "{(f \circ \proj_X, g \circ \proj_{E \times Y})}"'] & {}
\end{tikzcd} \] 
In particular, if $\cat{C} = \Set$, then the semi-direct product $f \semiprod{p} g \colon X \times Y \to X \times Y$ is a discrete dynamical system. Using basic properties of many of the objects we have defined, we are able to recover a generalization of the aforementioned theorem in \cref{modularity-theorem}.
Namely, if $(X \times Y, f \semiprod{p} g)$ is a semi-direct product, we get an isomorphism of $\mathbb{Z}/n$-sets:
\[ AS(X \times Y, f \semiprod{p} g)_n \cong \coprod\limits_{[c] \in \Orb (AS(X, f)_n)} AS(\ZZ / k \times Y, \succ \semiprod{pc} g)_n. \]
% In many real-world cases, discrete dynamical systems can be decomposed into a series of semi-direct products. 
This result allows one to more effectively analyze the attractors of a system by decomposing them into a disjoint union of much smaller sets.

\paragraph{Wiring diagram.}

In \cref{sec:wiring}, we formalize the notion of a wiring diagram for an object in $\DDS$.
In the context of biological networks, the wiring diagram is a digraph that records when the expression of one gene is affected by another.

\begin{remark}[On the term ``wiring diagram''] \label{rmk:wiring-terminology}
    The term \emph{wiring diagram} is unfortunately overloaded.
    In the applied category theory literature \cite{vagner-spivak-lerman,schultz-et-al}, a wiring diagram is a morphism in an operad: it is an \emph{interaction pattern} specifying how several open subsystems, each with its own inputs and outputs, are to be plugged into one another to form a composite system.
    % That is emphatically not what the term means here.
    In contrast, for us --- following the usage standard in the Boolean network and systems biology literature that our results generalize, e.g.\ \cite{kadelka2023modularity,chaves2018analysis,shilts:wiring_diagram} --- the wiring diagram of a \emph{single, closed} system $f \from A^n \to A^n$ is a digraph on the $n$ coordinates recording which coordinates influence which others.
    % It is an invariant extracted from one system, not a recipe for combining several.
    % Readers who prefer an unambiguous name may substitute \emph{dependency graph}, \emph{interaction graph}, or \emph{causal network} throughout; we retain ``wiring diagram'' only for continuity with \cite{kadelka2023modularity}, whose results we are generalizing.
\end{remark}
When modelling such a network using a Boolean function $f \from \FF_2^n \to \FF_2^n$, the wiring diagram equivalently records when the output of a coordinate function $f_j \from \FF_2^n \to \FF_2$ depends on the value of the $i$-th input coordinate (meaning $f_j(x_1, \dots, x_{i-1}, 0, x_{i+1} \dots, x_n) \neq f_i(x_1, \dots, x_{i-1}, 1, x_{i+1} \dots, x_n)$).
More precisely, the wiring diagram has $n$ vertices $\{ 1, \dots, n \}$ with an edge $i \to j$ if $f_j$ depends on the $i$-th coordinate.
A key result of \cite{kadelka2023modularity} is that a semi-direct product structure in a Boolean function can be detected by a subgraph in the wiring diagram with no incoming edges.
The categorical perspective enables us to generalize the wiring diagram and the notion of semi-direct product to the case of morphisms $f \from A^n \to A^n$ well beyond the case of sets and set functions.
Moreover, our perspective gives a new result: such an $f$ admits the structure of a semi-direct product if and only if its wiring diagram $W(f)$ admits a digraph map to the digraph $\DEdge$ with two looped vertices $0, 1$ and an edge $0 \to 1$ (\cref{semiprod-iff-wiring-map}).
Combining this with the decomposition theorem yields \cref{wiring-map-implies-attractors-decompose}: the mere existence of a graph map $W(f) \to \DEdge$ --- a finite, and typically small, check --- already guarantees that the attractors of $f$ decompose.

We state and prove these results for an arbitrary morphism $f \from A^n \to A^n$ in an arbitrary \emph{regular} category.
This, in particular, requires a generalization of what it means for the morphism $f_j \from A^n \to A$ to ``depend'' on the $i$-th coordinate (\cref{def:independent-input}).
Regularity is a mild technical hypothesis which implies that product projections are regular epimorphisms; every topos, every abelian category, and every variety of algebras is regular.
In \cref{subsec:wiring-examples}, we explain how these results specialize to several categories of interest: in $\Set$, it recovers the classical wiring diagram of a Boolean or multi-state network; in vector spaces, it says that a linear system decomposes precisely when its matrix can be put in block triangular form; in $G$-sets, it produces an equivariant refinement in the spirit of \cite{golubitsky-stewart}; and in $\DDS$ itself, it applies to networks whose nodes are themselves dynamical systems.

\subsection{Related work}
The program of using category theory to study dynamical systems has seen much interest in recent years and we conclude the introduction by reviewing some of the relevant work in the field.
Tao \cite{tao:three-categories} discusses three categories of dynamical systems: discrete, topological, and measurable, differentiating between the morphisms between these systems.
Abstractly, the categories of topological and measurable dynamical systems (with respect to an arbitrary topological/measurable monoid, not necessarily $\NN$) were shown to arise as certain categories of (co)algebras \cite{behrisch-kerkhof-poschel-schneider-siegmund}.
However, the notion of isomorphism provided by the category of discrete dynamical systems is too strong and one is instead interested in `shift equivalences' as studied by Bush \cite{bush:thesis-rutgers}, which in turn are an instance of generalized interleaving distances from persistence homology, as shown by Bubenik, de Silva, and Scott \cite{bubenik-de-silva-scott}.

% Several researchers have used category theory to analyze discrete dynamical systems, albeit their questions, and consequently methods, are quite different from ours.
% In particular, 
The works of both Schultz, Spivak, and Vasilakopoulou \cite{schultz-et-al}, and Fong, Soboci{\'n}ski, and Rapisarda \cite{fong-et-al}, focus on the use of monoidal categories with a view towards control theory.
The operadic notion of a wiring diagram studied by Vagner, Spivak, and Lerman \cite{vagner-spivak-lerman} belongs to the same circle of ideas; as noted in \cref{rmk:wiring-terminology}, it is a different notion from the one bearing that name here.

The category $\DDS = \Set^{\BN}$ is a Grothendieck topos, and it has been studied as such.
Hora and Kamio \cite{hora-kamio} classify its quotient toposes, using it to answer a question of Lawvere by exhibiting a Grothendieck topos with a proper class of quotient toposes.
While our concerns are different from theirs --- we are interested in the attractors of individual systems rather than in the global structure of the topos --- their work is a good illustration of how much structure this category carries.

More broadly, nothing in \cref{sec:framework} is special to the monoid $(\NN, +)$.
For any monoid $M$ one has the topos $\Set^{BM}$ of $M$-sets, and the adjoint triples we use are instances of the adjoint triple $u_! \adj u^* \adj u_*$ induced by a functor $u$ between small categories (\cref{adjoint-triple}). 
In the case of a monoid homomorphism, these are precisely the adjunctions studied systematically by Hemelaer and Rogers \cite{hemelaer-rogers}, whose paper we refer the reader to for a far more thorough treatment than we give here.
We have chosen to keep $M = (\NN, +)$ fixed throughout because it is the case relevant to the applications that motivate us, and because the notion of a cycle set does depend on more than the monoid alone; but the reader will have no trouble identifying which of our results are instances of the general theory.

Finally, we mention the work of Golubitsky and Stewart \cite{golubitsky-stewart} on coupled cell systems.
They study systems of ordinary differential equations organized over a network, with the aim of understanding how the symmetries of the network --- encoded in a groupoid rather than a group, so as to capture partial symmetries --- constrain the dynamics, in particular producing synchrony.
Our strategy is similar in spirit, in that a combinatorial invariant of the network is used to predict features of the dynamics, but the invariant (the wiring diagram, rather than the symmetry groupoid), the property predicted (decomposition of attractors, rather than synchrony), and the class of systems (discrete, rather than ODEs) are all different.
This list is, of course, far from exhaustive.

\subsection*{Acknowledgements}

We thank the editor and the two anonymous referees for their exceptionally careful and generous reports, which corrected several errors and led to substantial improvements in both the results and their presentation.
In particular, the treatment of independent inputs via regular categories and support objects in \cref{sec:wiring}, the separation of the wiring diagram material into its own section, and a number of the examples were suggested by the referees.

\section{Framework} \label{sec:framework}

In this section, we introduce the categorical framework for studying discrete dynamical systems.
We begin by constructing the category $\DDS$ of discrete dynamical systems and showing that it is a reflective and co-reflective subcategory of the category of digraphs via the state space functor (\cref{state-space-admits-both-adjoints,state-space-functor-full-faithful,state-space-full-subcat}).
We then introduce the notion of a cycle set (\cref{cycle-set}) and construct the functor taking a digraph to its set of attractors (\cref{attractor-functor}).
We provide an in-depth categorical analysis of the category of cycle sets, which makes our framework robust.

\subsection{The category of discrete dynamical systems}

Let $\BN$ denote the monoid of natural numbers under addition, viewed as a category with one object (which we denote by $\BNobj$).

\begin{definition} \label{def:DDS}
    The \emph{category of discrete dynamical systems}, denoted $\DDS$, is the category of presheaves on $\BN$, that is, the category of functors $\BN^\op \to \Set$.
\end{definition}

Since $(\NN, +)$ is a commutative monoid, we have an equality of categories $\BN^\op = \BN$, and hence an equality of functor categories $\Set^{\BN^\op} = \Set^{\BN} = \DDS$.
We will therefore not distinguish between covariant and contravariant functors out of $\BN$, and will freely write $\BN \to \Set$ for an object of $\DDS$; the reader should keep in mind that ``presheaf'' means contravariant functor, and that it is only this coincidence that lets us be cavalier about it.

As a functor category into $\Set$, $\DDS$ has limits and colimits; in fact it is a Grothendieck topos.
We will make use of this in \cref{sec:wiring}, where we need only the weaker fact that $\DDS$, like every topos, is a \emph{regular} category.
The topos $\Set^{\BN}$ is itself an object of study: Hora and Kamio \cite{hora-kamio} have classified its quotient toposes, and the adjunctions we describe below are instances of a general theory of geometric morphisms between toposes of monoid actions developed by Hemelaer and Rogers \cite{hemelaer-rogers}.

Since $\BN$ is freely generated by the object $\ast$ and the endomorphism $1 \from \ast \to \ast$, a functor $F \from \BN \to \Set$ is uniquely determined by a set $X$, together with a single endofunction $f \from X \to X$ (which represents the update function of the system).
From this data, the corresponding functor $F \from \BN \to \Set$ is defined by sending the object $\ast$ to $X$ and the morphism $n \from \ast \to \ast$ to the $n$-fold composition $f^{n}$; in particular, the function $f$ is recovered as the morphism $F(1) \from F(\ast) \to F(\ast)$.

% \begin{example}
%     One example of a discrete dynamical system, and one that we will refer to frequently, is that of a Boolean network. In general a functor $F\in\DDS$ is a Boolean network if the state space $X$ is in bijection with a Boolean algebra. For instance, the function $F(X,Y,Z,W) = (Y,X,Y\, \text{XOR}\, W, Z)$ is a Boolean network. One can equivalently think of $F$ as a polynomial function $F: \FF_2^4 \to \FF_2^4$ where now $F(x,y,z,w)=(y,x,x+w,z).$. In Figure~\ref{fig:tt}, we've represented this network using its truth table which completely determines the network.
%     \fxnote{Working example of a Boolean network.}
%     \fxnote{Example of multistate network.}
%     \fxnote{Example of petri nets}
%     \fxnote{Example of discrete differential network}
% \end{example}
% \begin{figure}[h]
%     \centering
%     \includegraphics[width=0.8\textwidth]{Rplot}\\
%     \label{fig:tt}
% \end{figure}

% \begin{example}
%     \fxnote{Example of discrete differential network: Predator-Prey} A classical dynamical system is the predator-prey model, and a version of this model, the discrete predator-prey model \cite{Satoh-Numerical,Al Kaff-Exploring Chaos}, is a discrete dynamical system. The discrete predator-prey model is given by
%     \[\begin{cases}
%         x_{n+1}\,=\,ax_n(1-x_n)-cx_ny_n,\\
%         y_{n+1}\,=\, by_n(1-y_n)-dx_ny_n.
%     \end{cases}\]
%     To see the discrete predator prey model as an object in $\DDS$, we can represent it as the function 
%     \[F(x,y)=(ax(1-x)-cxy,by(1-y)-dxy)\]
%     which acts on the state space $X=[0,1]^2.$
% \end{example}
A functor $\BN \to \Set$ is also the same data as a set $X$ equipped with a monoid action by the natural numbers $(\NN, +)$.
We occasionally borrow this notation, writing $n \cdot x$ for the value $f^n(x)$.

A natural transformation between two functors $(X, f)$ and $(Y, g)$ is a function $\alpha \from X \to Y$ making the square
\[ \begin{tikzcd}
    X \ar[d, "f"'] \ar[r, "\alpha"] & Y \ar[d, "g"] \\
    X \ar[r, "\alpha"] & Y
\end{tikzcd} \]
commute.
A priori, this condition merely asserts naturality at the morphism $1 \from \ast \to \ast$, however naturality at an arbitrary morphism $n \from \ast \to \ast$ follows from this by concatenating squares.
Recall that a natural transformation is an isomorphism precisely when its components are isomorphisms, from which we deduce the following useful fact:
\begin{proposition}
    A morphism $\alpha \from (X, f) \to (Y, g)$ in $\DDS$ is an isomorphism precisely when the underlying function $\alpha \from X \to Y$ is a bijection.
\end{proposition}
Again, using the language of monoid actions, we occasionally refer to such a natural transformation as an \emph{$\NN$-equivariant map}, and write the naturality condition as $\alpha(n \cdot x) = n \cdot \alpha(x)$ for all $n \in \NN$ and $x \in X$.

% For the remainder of our development, we will not distinguish between the category $\BN$ and $(\BN)^\op$, as all our constructions will be invariant under this isomorphism.
% Let $(\NN, \succ)$ denote the discrete dynamical system given by the natural numbers $\NN$ with the endofunction $\succ \from \NN \to \NN$ defined by $n \mapsto n+1$.
% Unfolding to verify that the functor $(\NN, \succ)$ is \emph{representable}.
% The dynamical system $(\NN, \succ)$ is a \emph{representable functor}, represented by the object $\ast \in \BN$.
\begin{example} 
    An important example of a discrete dynamical system is the \emph{representable} functor $\BN(\ast, -) \from \BN \to \Set$.
    Unfolding the definitions, this functor corresponds to the pair $(\NN, \succ)$, i.e.\ the dynamical system given by the natural numbers equipped with the endofunction $n \mapsto n+1$.
\end{example}
\begin{example}
    For $n \in \NN$, the pair $(\ZZ / n, \succ)$ is a dynamical system where $\ZZ / n$ denotes the set $\{ 1, \dots, n \}$ and $\succ$ denotes the endofunction $i \mapsto i + 1 \ \operatorname{mod} n$.
    This dynamical system is the quotient of $(\NN, \succ)$ by the equivariant map $(\NN, \succ) \to (\NN, \succ)$ defined by $i \mapsto i+n$.
    More precisely,
    \[ (\ZZ / n, \succ) \cong \colim \left( \begin{tikzcd}[column sep = 2.9em]
        (\NN, \succ ) \ar[r, yshift=0.8ex, "{i \mapsto i + n}"] \ar[r, yshift=-0.8ex, "\id"'] & (\NN, \succ )
    \end{tikzcd} \right). \]
\end{example}

% We now fix a small dynamical system that we will return to at each stage of the development, computing in turn its state space graph, its cycle set of attractors, the $(\ZZ/n)$-sets obtained by evaluating that cycle set, its degenerate cycles, and its digraph realization; in \cref{sec:wiring} we will use it again to illustrate the decomposition theorem.

% The following will be a running example 
% We will use the following example throughout the paper as a running 
\begin{example}[Running example]\label{ex:running}
    Let $X := \{ a, b, c \}$ and let $f \from X \to X$ be given by
    \[ f(a) = a, \qquad f(b) = c, \qquad f(c) = b. \]
    Thus $a$ is a fixed point and $b, c$ form a $2$-cycle.
    We write $(X, f)$ for the resulting object of $\DDS$.
\end{example}
We will use \cref{ex:running} as a running example throughout the paper.

Under the identification $\Set^{\BN^\op} = \Set^{\BN}$ noted after \cref{def:DDS}, the covariant and contravariant representables $\BN(\ast, -), \BN(-, \ast)$ are in fact identical.
This identification gives us access to the \emph{Yoneda embedding}, the functor $\BN \to \DDS$ which sends the object $\ast$ to $(\NN, \succ)$ and the morphism $n$ to $\succ^n$, which we view as a \emph{natural transformation} $(\NN, \succ) \to (\NN, \succ)$; this is justified since $\succ$ is an $\NN$-equivariant map, i.e.\ it commutes with itself.

There is a forgetful functor $U \from \DDS \to \Set$ which sends $(X, f)$ to $X$.
It follows from the (co)Yoneda lemma that this functor is representable.
\begin{proposition} \label{forgetful-representable}
    The forgetful functor $U \from \DDS \to \Set$ is representable, represented by $(\NN, \succ)$. \qed
\end{proposition}
This means that, for a discrete dynamical system $(X, f)$, a morphism $(\NN, \succ) \to (X, f)$ is the same data as an element of the set $X$.
Concretely, given an element $x \in X$, the corresponding function $(\NN, \succ) \to (X, f)$ is defined by $0 \mapsto x$ and $n \mapsto f^n(x)$. 
One can think of the choice $0 \mapsto x$ as a choice of initial condition; the morphism $(\NN,\succ) \to (X,f)$ then maps out the entire trajectory of $x$ under the system. 

% All three of the adjoint triples appearing in this paper --- the one below, the one for the state space functor (\cref{state-space-admits-both-adjoints}), and the one for the evaluation functor on cycle sets (\cref{ev-admits-adjoints}) --- arise by one and the same mechanism.
% We introduce the following as our primary way of constructing adjunctions 
To construct adjoint functors to and from $\DDS$, we use the following proposition.

\begin{proposition} \label{adjoint-triple}
    Let $u \from \cat{A} \to \cat{B}$ be a functor between small categories.
    Then the restriction functor
    \[ u^* \from \Set^{\cat{B}^\op} \to \Set^{\cat{A}^\op}, \qquad u^*(F) := F \circ u^\op, \]
    admits both a left adjoint $u_!$ and a right adjoint $u_*$, given by left and right Kan extension along $u^\op$ respectively:
    \[ u_! \adj u^* \adj u_* . \]
    Explicitly, these are computed by the coend and end formulas
    \[ u_!(F)(b) \cong \int^{a \in \cat{A}} \cat{B}(b, ua) \times F(a), \qquad u_*(F)(b) \cong \int_{a \in \cat{A}} F(a)^{\cat{B}(ua, b)}. \]
    In particular, $u^*$ preserves all limits and colimits.
\end{proposition}
\begin{proof}
    This is standard; see, e.g., \cite[\S 6.2]{riehl:context}.
    The Kan extensions exist because $\Set$ is cocomplete and complete and $\cat{A}$ is small, and the displayed formulas are the usual pointwise ones.
    % The final claim holds because a functor admitting both adjoints preserves whatever limits and colimits exist.
\end{proof}

\begin{remark} \label{rmk:adjoint-triple-monoids}
    When $\cat{A}$ and $\cat{B}$ are one-object categories, i.e.\ monoids, and $u$ is a monoid homomorphism, \cref{adjoint-triple} specializes to an adjoint triple between toposes of monoid actions.
    These are studied systematically in \cite{hemelaer-rogers}, and several of the adjunctions below are instances of that theory.
\end{remark}

Recall that the category of sets coincides with the category of presheaves on the terminal category $\mathbf{1}$, i.e.\ the category with one object and only its identity morphism.
The forgetful functor $U \from \DDS \to \Set$ is exactly the functor $u^*$ of \cref{adjoint-triple} for the inclusion $u \from \mathbf{1} \ito \BN$, from which it follows that $U$ admits both adjoints.
We briefly describe these adjoints in the following proposition.
\begin{proposition} \leavevmode
    \begin{enumerate}
        \item The forgetful functor admits a left adjoint $\Set \to \DDS$ which sends a set $X$ to the pair $(X \times \NN, \id[X] \times \succ)$.
        \item The forgetful functor admits a right adjoint $\Set \to \DDS$ which sends a set $X$ to the pair $(X^\NN, \succ^*)$ \qed
    \end{enumerate}
\end{proposition} 
Here, $X^\NN$ denotes the set of functions $\NN \to X$, and $\succ^*$ is the endofunction $X^\NN \to X^\NN$ which sends $\varphi \in X^\NN$ to $\varphi \circ \succ$.

\subsection{Digraphs and the state space functor}
% There are two important directed graphs related to a Boolean network, the {\em wiring diagram} and the {\em state space graph}, and they play an important role in the decomposition theory developed in \cite{kadelka2023modularity}. The constructions of these graphs generalize naturally to the setting of discrete dynamical systems. Here, we lay down the categorical foundation for these constructions.

Define a category $\GG$ with two objects $\{ V, E \}$ and two non-identity morphisms $s, t \from V \to E$.
\[ \GG := \begin{tikzcd}
    V \ar[r, yshift=0.8ex, "s"] \ar[r, yshift=-0.8ex, "t"'] & E
\end{tikzcd} \]
\begin{definition}
    % A \emph{multidigraph} is a functor $\GG^\op \to \Set$.
    The \emph{category of multidigraphs} $\MultiDiGraph$ is the category of functors $\Set^{\GG^\op}$.
\end{definition}
A \emph{multidigraph} is an object of $\MultiDiGraph$ (i.e.\ a functor $\GG^\op \to \Set$) and a \emph{multigraph map} is a morphism (i.e.\ a natural transformation).
Unfolding this definition recovers the ``usual'' definition of a multidigraph, since a functor $G \from \GG^\op \to \Set$ consists of a set of \emph{vertices} (the value of $G$ at $V$) and a set of \emph{edges} (the value of $G$ at $E$) equipped with \emph{source} and \emph{target} maps (the values of $G$ on $s$ and $t$, respectively).
We write $V(G)$ for the set of vertices and $E(G)$ for the set of edges.
\[ \begin{tikzcd}[sep = large]
    V(G) & E(G) \ar[l, yshift=0.8ex, "\textrm{source}"'] \ar[l, yshift=-0.8ex, "\textrm{target}"]
\end{tikzcd} \]
When speaking of edges in a multidigraph, we write $e \from v \to w$ or $v \xrightarrow{e} w$ to mean $e$ is an element of $E(G)$ with source $v \in V(G)$ and target $w \in V(G)$.

A \emph{digraph} is a multidigraph such that no two (distinct) edges have the same source and target.
We rephrase this condition in the following definition.
\begin{definition}
    A \emph{digraph} is a multidigraph with the property that the map
    \[ (\textrm{source}, \textrm{target}) \from E(G) \to V(G) \times V(G) \]
    is a monomorphism.
    The \emph{category of digraphs} $\DiGraph$ is the full subcategory $\MultiDiGraph$ spanned by digraphs.
\end{definition}
We refer to a morphism in $\DiGraph$ as a \emph{graph map}.

Every multidigraph can be ``collapsed'' to a digraph by identifying edges in $E(G)$ with the same source and target.
This process gives a functor $\Gfun \from \MultiDiGraph \to \DiGraph$ which is left adjoint to the inclusion $i \from \DiGraph \ito \MultiDiGraph$.

Define a functor $\Sigma \from \GG \to \BN$ by sending the morphisms $s, t$ in $\GG$ to $0, 1$ in $\BN$, respectively.
\[ \begin{tikzcd}
    V \ar[r, yshift=0.8ex, "s"] \ar[r, yshift=-0.8ex, "t"'] & E \ar[rr, mapsto, shorten=2em, "\Sigma"] & {} & \ast \ar[r, yshift=0.8ex, "0"] \ar[r, yshift=-0.8ex, "1"'] & \ast
\end{tikzcd} \]
Applying \cref{adjoint-triple}, this induces a functor $\Sigma^* \from \DDS \to \MultiDiGraph$ by pre-composition, which admits both a left adjoint $\Sigma_! \from \MultiDiGraph \to \DDS$ and a right adjoint $\Sigma_* \from \MultiDiGraph \to \DDS$.
We give explicit descriptions of each of these functors.

% \begin{proposition} \label{sigma-adj-explicit-descriptions}
    \begin{itemize}
        \item Given a dynamical system $(X, f)$, the multidigraph $\Sigma^*(X, f)$ has
        \begin{itemize}
            \item as vertices, elements of $X$;
            \item an edge $x \to f(x)$ for all $x \in X$.
        \end{itemize}
        \item Given a multidigraph $G$, the dynamical system $\Sigma_!(G)$ has
        \begin{itemize}
            \item as its underlying set, the Cartesian product $V(G) \times \NN$, quotiented by the equivalence relation generated by $(v, n+1) \sim (w, n)$ for every edge $e \from v \to w$ in $G$ and $n \in \NN$;
            \item the endofunction on the quotient $\big( V(G) \times \NN \big) / \sim$ sends an equivalence class $[v, n]$ to $[v, n+1]$.
        \end{itemize}
        \item Given a multidigraph $G$, the dynamical system $\Sigma_*(G)$ has
        \begin{itemize}
            \item as its underlying set, the set of infinite paths 
            $(v_0 \xrightarrow{e_0} v_1 \xrightarrow{e_1} v_2 \xrightarrow{e_2} \cdots)$
            in $G$ with a distinguished starting point.
            \item the endofunction takes a path
            $(v_0 \xrightarrow{e_0} v_1 \xrightarrow{e_1} v_2 \xrightarrow{e_2} \cdots)$
            and discards the starting edge, outputting the path
            $(v_1 \xrightarrow{e_1} v_2 \xrightarrow{e_2} v_3 \xrightarrow{e_3} \cdots)$. \qed
        \end{itemize}
    \end{itemize}
% \end{proposition}

From the explicit description, we see that $\Sigma^*(X, f)$ is always a digraph, i.e.\ there are no two distinct edges with the same source and target.
Thus, we may view $\Sigma^*$ as taking values in $\DiGraph$ rather than $\MultiDiGraph$.
In light of this, we make the following definition.
\begin{definition}
    The \emph{state space} functor, denoted $S \from \DDS \to \DiGraph$, is the unique functor fitting into the commutative triangle
    \[ \begin{tikzcd}
        \DDS \ar[rr, "\Sigma^*"] \ar[rd, dotted, "S"'] & {} & \MultiDiGraph \\
        {} & \DiGraph \ar[ur, hook] & {}
    \end{tikzcd} \]
\end{definition}
One verifies that the adjoint triple $\Sigma_! \adj \Sigma^* \adj \Sigma_*$ descends to a related adjoint triple $(\Sigma_! \circ i) \adj S \adj (\Sigma_* \circ i)$ for the state space functor.
We denote the leftmost functor by $F := \Sigma_! \circ i$ and the rightmost functor by $P := \Sigma_* \circ i$.
\begin{proposition} \label{state-space-admits-both-adjoints}
    The state space functor $S \from \DDS \to \DiGraph$ fits into an adjoint triple $F \adj S \adj P$. \qed
\end{proposition}
In particular, the descriptions given before remain valid for the state space functor and its adjoints.
\begin{proposition}
    Let $(X, f)$ be a dynamical system and $G$ be a digraph.
    \begin{enumerate}
        \item The state space digraph $S(X, f)$ has
        \begin{itemize}
            \item as vertices, elements of $X$;
            \item an edge $x \to f(x)$ for all $x \in X$.
        \end{itemize}
        \item The dynamical system $F(G)$ has
        \begin{itemize}
            \item as its underlying set, the Cartesian product $V(G) \times \NN$, quotiented by the equivalence relation generated by $(v, n+1) \sim (w, n)$ for every edge $e \from v \to w$ in $G$ and $n \in \NN$;
            \item the endofunction on the quotient $\big( V(G) \times \NN \big) / \sim$ sends an equivalence class $[v, n]$ to $[v, n+1]$.
        \end{itemize}
        \item The dynamical system $P(G)$ has
        \begin{itemize}
            \item as its underlying set, the set of infinite paths 
            $(v_0 \to v_1 \to v_2 \to \cdots)$
            in $G$ with a distinguished starting point.
            \item the endofunction takes a path
            $(v_0 \to v_1 \to v_2 \to \cdots)$
            and discards the starting edge, outputting the path
            $(v_1 \to v_2 \to v_3 \to \cdots)$. \qed
        \end{itemize}
    \end{enumerate}
\end{proposition}
% \begin{proof}
%     Follows from \cref{state-space-admits-both-adjoints}.
% \end{proof}
\begin{example}
    The state space of the dynamical system $(\NN, \succ)$ is the digraph with infinitely many vertices connected in a (directed) line.
    \[ \begin{tikzpicture}
        \node[vertex, label={0}] (0) {};
        \path (0) -- +(1, 0) node[vertex, label={1}] (1) {};
        \path (1) -- +(1, 0) node[vertex, label={2}] (2) {};
        \path (2) -- +(1, 0) node[vertex, label={3}] (3) {};
        \path (3) -- +(1, 0) node (dots) {$\ldots$};

        \path (0) -- +(-0.4, 0.2) node[anchor=east] {\large $S(\NN, \succ) = $};

        \draw[thick, ->]
            (0) edge (1)
            (1) edge (2)
            (2) edge (3)
            (3) edge (dots);
    \end{tikzpicture} \]
\end{example}
\begin{example}
    The state space of the dynamical system $(\ZZ / n, \succ)$ is the cycle graph with $n$ vertices connected in a cycle.
    \[ \begin{tikzpicture}
        \coordinate (center) {};
        \path (center) -- +(234:0.8) node[vertex, label={south west:0}] (0) {};
        \path (center) -- +(306:0.8) node[vertex, label={south east:1}] (1) {};
        \path (center) -- +(18:0.8) node[vertex, label={east:2}] (2) {};
        \path (center) -- +(90:0.8) node (dots) {$\ldots$};
        \path (center) -- +(162:0.8) node[vertex, label={west:$n-1$}] (last) {};

        \path (last) -- +(-1.35, -0.3) node[anchor=east] {\large $S(\ZZ / n, \succ) = $};

        \draw[thick, ->]
            (0) edge (1)
            (1) edge (2)
            (2) edge (dots)
            (dots) edge (last)
            (last) edge (0);
    \end{tikzpicture} \]
\end{example}
\begin{example}[Running example, continued] \label{ex:running-state-space}
    The state space graph of the system $(X, f)$ of \cref{ex:running} is
    \[ \begin{tikzpicture}
        \node[vertex, label={[label distance=-1ex]below:$a$}] (A) {};
        \path (A) -- +(1.6, 0) node[vertex, label={[label distance=-1ex]below:$b$}] (B) {};
        \path (B) -- +(1.2, 0) node[vertex, label={[label distance=-1ex]below:$c$}] (C) {};

        \path (A) -- +(-0.5, 0.15) node[anchor=east] {\large $S(X, f) = $};

        \draw[thick, ->] (A) to [out=45, in=125, looseness=6] (A);
        \draw[thick, ->] (B) to [bend right] (C);
        \draw[thick, ->] (C) to [bend right] (B);
    \end{tikzpicture} \]
    We observe that every vertex has exactly one outgoing edge (see \cref{state-space-full-subcat}).
\end{example}

A functor which admits both adjoints always preserves both limits and colimits.
\begin{corollary} \label{cor:S_preserves_lim_colim}
    Both the functor $\Sigma^* \from \DDS \to \MultiDiGraph$ and the state space functor $S \from \DDS \to \DiGraph$ preserve limits and colimits. \qed
\end{corollary}

An important property of the state space functor is that it is full and faithful.
\begin{theorem} \label{state-space-functor-full-faithful}
    Both the functor $\Sigma^* \from \DDS \to \MultiDiGraph$ and the state space functor $S \from \DDS \to \DiGraph$ are full and faithful.
\end{theorem}
\begin{proof}
    Since the inclusion $\DiGraph \ito \MultiDiGraph$ is full and faithful, it suffices to show $S$ is full and faithful.

    For any discrete dynamical system $(X, f)$, the state space digraph $S(X, f)$ has the property that every vertex has a unique edge pointing out.
    Thus, an infinite path $v_0 \to v_1 \to v_2 \to \dots$ in $S(X, f)$ is uniquely determined by where it sends the vertex $0 \in S(X, f)$.
    It follows from this that the unit $\eta_{(X, f)} \from (X, f) \to P(S(X, f))$ of the $(S \adj P)$-adjunction, whose explicit formula is given by 
    \[ \eta_{(X, f)}(x) = (x \to f(x) \to f^2(x) \to \dots), \]  
    is a bijection.
\end{proof}
\begin{corollary} \label{state-space-full-subcat}
    The state space functor identifies $\DDS$ as the full subcategory of $\DiGraph$ spanned by digraphs such that every vertex has a unique edge pointing out.
\end{corollary}
\begin{proof}
    By \cref{state-space-functor-full-faithful}, the state space functor identifies $\DDS$ as the full subcategory of $\DiGraph$ spanned by its essential image.
    An endofunction on a set $X$ is exactly a binary relation such that every element has a unique edge pointing out, from which it follows that this is the essential image of the state space functor.
\end{proof}
\Cref{state-space-full-subcat} says that passing from a system to its state space graph loses nothing and adds nothing: the systems $(X,f)$ correspond exactly to those digraphs in which every vertex has exactly one outgoing edge, and morphisms of discrete dynamical systems correspond to graph maps.

\begin{remark}[On the choice of simple digraphs] \label{rmk:simple-digraphs}
    In category theory, the category $\MultiDiGraph$ is generally preferred to $\DiGraph$: it is a presheaf category, and it is in this sense that every category has an underlying digraph.
    For instance, the underlying digraph of $\BN$ has a single vertex and one loop for every $n \in \NN$, and so is not simple.
    We have nonetheless chosen to work with simple digraphs, for two reasons.
    First, the digraphs of primary interest in this paper are simple: state space graphs (\cref{state-space-full-subcat}), directed cycles, wiring diagrams (\cref{sec:wiring}), and the digraph $\DEdge$.
    Second, simple digraphs are the standard notion in the literature on Boolean networks that this paper generalizes, and we hope to keep the paper accessible to readers with that background.

    Little depends on this choice, however.
    The functor $\Sigma^*$ is already full and faithful as a functor into $\MultiDiGraph$ and preserves all limits and colimits (\cref{cor:S_preserves_lim_colim,state-space-functor-full-faithful}).
    The attractor functor of \cref{subsec:attractors} may be defined on multidigraphs by the same formula, $A(G)_n := \MultiDiGraph(C_n, G)$, and the results of \cref{subsec:attractors,subsec:degenerate,subsec:characterization} hold for it with the same proofs.
    Finally, since $\Gfun$ is left adjoint to the inclusion, a multidigraph $G$ admits a map to a simple digraph $H$, such as $\DEdge$, if and only if $\Gfun(G)$ does.
    The reader who prefers to do so may therefore read $\DiGraph$ as $\MultiDiGraph$ throughout.
\end{remark}
% We stress that this is a statement about the state space \emph{graph}, which encodes $f$ in its edges, and not about the underlying set $X$ of states, which of course carries no dynamics at all and typically supports many different systems (cf.\ \cref{rmk:state-space-terminology}).

\subsection{Cycle sets} \label{subsec:cycle-sets}

In this subsection, we introduce a central concept of the paper, namely that of a cycle set.
Cycle sets form a presheaf category with the indexing category given by the full subcategory of digraphs spanned, perhaps unsurprisingly, by the cycles.
This subsection is concerned with cycle sets as such.
% ; the functor relating them to digraphs is constructed in \cref{subsec:attractors}.

For $n \geq 1$, the \emph{directed $n$-cycle} $C_n$ is the digraph with vertices $\{ 0, \dots, n-1 \}$ connected in a cycle.
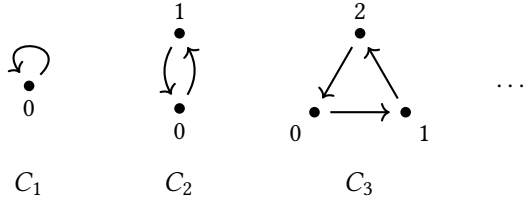
\[ \begin{tikzpicture}
    \node[coordinate] (C1center) {};
    \node[vertex, label={[label distance=-1ex]below:0}] (C1A) at (C1center) {};

    \path (C1center) -- +(2, 0.2) node[coordinate] (C2center) {};
    \path (C2center) -- +(90:0.5) node[vertex, label={[label distance=-1ex]above:1}] (C2T) {};
    \path (C2center) -- +(-90:0.5) node[vertex, label={[label distance=-1ex]below:0}] (C2B) {};

    \path (C2center) -- +(2.4, -0.2) node[coordinate] (C3center) {};
    \path (C3center) -- +(90:0.7) node[vertex, label={[label distance=-1ex]above:2}] (C3T) {};
    \path (C3center) -- +(210:0.7) node[vertex, label={[label distance=-1ex]below left:0}] (C3BL) {};
    \path (C3center) -- +(330:0.7) node[vertex, label={[label distance=-1ex]below right:1}] (C3BR) {};

    \path (C3center) -- +(2, 0) node (dotsA) {$\ldots$};

    \draw[thick, ->] (C1A) to [out=45, in=125, looseness=6] (C1A);

    \draw[thick, ->] (C2B) to [bend right] (C2T);
    \draw[thick, ->] (C2T) to [bend right] (C2B);
    
    \draw[thick, ->] (C3BL) to [bend right=0] (C3BR);
    \draw[thick, ->] (C3BR) to [bend right=0] (C3T);
    \draw[thick, ->] (C3T) to [bend right=0] (C3BL);

    % Labels
    \path (C1center) -- +(0, -1.3) node (C1name) {\large $C_1$};
    \path (C1name) -- +(2, 0) node (C2name) {\large $C_2$};
    \path (C2name) -- +(2.4, 0) node {\large $C_3$};
    % \path (dotsA) -- +(0, -1.1) node {$\ldots$};
\end{tikzpicture} \]
Let $\Cy$ denote the full subcategory of $\DiGraph$ spanned by the directed $n$-cycles.
When $n$ divides $m$, there is a graph map $\modmap{m}{n} \from C_m \to C_n$ defined by
\[ \modmap{m}{n}(k) := k \ \operatorname{mod} n. \]
Additionally, each cycle $C_n$ comes with $n$ automorphisms $\autmap{n}{0}, \dots, \autmap{n}{n-1}$ defined by
\[ \autmap{n}{i}(k) := k + i \ \operatorname{mod} n. \]
It turns out that every map between directed cycles is a composite of these two in a unique way.
\begin{proposition} \label{graph-maps-between-cycles}
    \begin{enumerate}
        \item A graph map $C_m \to C_n$ exists if and only if $n$ divides $m$.
        \item Given a graph map $\varphi \from C_m \to C_n$, there exists a unique $i \in \{ 0, \dots, n-1 \}$ such that
        \[ \varphi = \autmap{n}{i} \circ \modmap{m}{n}. \]
        That is,
        \[ \varphi(k) = k + i \ \operatorname{mod} n. \]
    \end{enumerate}
\end{proposition}
\begin{proof}
    We prove (2) first, then prove (1).

    The graphs $C_m$ and $C_n$ have the property that every vertex admits a unique edge pointing out.
    % Thus, a graph map $\varphi \from C_m \to C_n$ is uniquely determined by where it sends the vertex $0 \in C_m$.
    Since $\varphi$ is a graph map (i.e.\ it preserves edges), we deduce an equality
    \[ \varphi(k) = \varphi(0) + k \ \operatorname{mod} n \]
    for all $k \in C_m$ by an induction argument on $k$.
    This proves the desired formula by setting $i := \varphi(0)$.
    Uniqueness of $i$ follows by definition.
    
    For (1), we take the above formula as the definition of a set function $\varphi \from V(C_m) \to V(C_n)$.
    This assignment is a valid graph map if and only if $\varphi$ preserves the edge $(m - 1) \to 0$ in $C_m$, as all other edges are preserved by construction. 
    This edge is preserved if and only if $\varphi(m-1) = \varphi(0) - 1 \operatorname{mod} n$.
    By construction, we have that $\varphi(m-1) = \varphi(0) + m - 1 \operatorname{mod} n$, hence the desired equality holds if and only if $m = 0 \operatorname{mod} n$; that is, if and only if $n$ divides $m$.
\end{proof}
From \cref{graph-maps-between-cycles}, we may give an explicit description of the automorphism groups of the cycle graphs.
\begin{corollary}
    The group of automorphisms of $C_n$ is isomorphic to the cyclic group $(\ZZ / n, +)$ of order $n$.
    The isomorphism $\ZZ / n \xrightarrow{\cong} \operatorname{Aut}(C_n)$ is given by the assignment $i \mapsto \autmap{n}{i}$. \qed
\end{corollary}
We also have a presentation of the category $\Cy$ of cycles in terms of generators and relations.
\begin{corollary} \label{cycle-cat-presentation}
    The category $\Cy$ is isomorphic to the category whose
    \begin{itemize}
        \item objects are the symbols $C_n$ indexed by positive integers $n \in \ZZ_{\geq 1}$;
        \item morphisms are generated by the arrows
        \[ \modmap{m}{n} \from C_m \to C_n \]
        whenever $n$ divides $m$, along with an arrow
        \[ \autmap{n}{1} \from C_n \to C_n \]
        for all $n$, subject to the relations
        \[ \begin{array}{c@{\qquad}c@{\qquad}c}
            \modmap{m}{n} \modmap{\ell}{m} = \modmap{\ell}{n} & (\autmap{n}{1})^{n} = \id[C_n] = \modmap{n}{n} & \modmap{m}{n} \autmap{m}{1} = \autmap{n}{1} \modmap{m}{n}
        \end{array} \]
        for all $\ell, m, n \geq 1$ with $n$ dividing $m$ dividing $\ell$.
    \end{itemize}
\end{corollary}
\begin{proof}
    Let $\cat{A}$ denote the small category defined by the given presentation.
    % Observe that item (1) of \cref{graph-maps-between-cycles} holds for $\cat{A}$ by definition of the generating morphisms, and statement (2) can be proven using the given relations (where $\autmap{n}{i}$ is interpreted as a composite $(\autmap{n}{1})^i$).
    Since the given relations hold in $\Cy$, we have a canonical functor $\cat{A} \to \Cy$ which is a bijection on objects.
    It is surjective on morphisms by \cref{graph-maps-between-cycles}.

    Regarding injectivity on morphisms, observe that item (1) of \cref{graph-maps-between-cycles} holds for $\cat{A}$ by definition of the generating morphisms, and item (2) can be proven using the given relations (where $\autmap{n}{i}$ is interpreted as a composite $(\autmap{n}{1})^i$).
    The uniqueness property implies the desired injectivity.
    % This fact, combined with \cref{graph-maps-between-cycles}, defines the desired isomorphism.
\end{proof}
% Note that a graph map $C_m \to C_n$ exists if and only if $n$ divides $m$, and is always determined by where it sends the vertex 1, since every vertex in $C_n$ has a unique edge pointing out.

\begin{definition} \label{cycle-set}
    A \emph{cycle set} is a functor $\Cy^\op \to \Set$.
    The \emph{category of cycle sets} $\cySet$ is the functor category $\Set^{\Cy^\op}$.
\end{definition}
% Given an $n$-cycle $x$ and a morphism $f\varphi \from C_m \to C_n$, we write $x\varphi$ for the $m$-cycle given by evaluating $K(\varphi) \from C_n \to C_m$ at $x$.
As a category of functors into $\Set$, the category of cycle sets is a presheaf category, and thus admits all limits and colimits.

Given a cycle set $K \in \cySet$, we write $K_n$ for the set $K(C_n)$.
We refer to an element $x \in K_n$ as an \emph{$n$-cycle} of $K$. 
Given a morphism $\varphi \from C_m \to C_n$ in $\Cy$ and $x \in K_n$, we write $x\varphi$ for the $m$-cycle $K(\varphi)(x) \in K_m$.
Via \cref{cycle-cat-presentation}, a cycle set is uniquely defined by the data of:
\begin{itemize}
    \item a collection of sets $K_n$ for $n \geq 1$;
    \item a function $K(\modmap{m}{n}) \from K_n \to K_m$ for all $n$ dividing $m$; and
    \item a bijection $K(\autmap{n}{1}) \from K_n \to K_n$ for all $n$
\end{itemize}
satisfying the equalities
\[ \begin{array}{c@{\qquad}c@{\qquad}c}
    x\modmap{m}{n} \modmap{\ell}{m} = x\modmap{\ell}{n} & x(\autmap{n}{1})^{n} = x = x\modmap{n}{n} & x\modmap{m}{n} \autmap{m}{1} = x\autmap{n}{1} \modmap{m}{n}
\end{array} \]
for all $x \in K_n$ and $m, n \geq 1$ with $n$ dividing $m$.

The operators $\modmap{m}{n}$ are what allow a cycle to be regarded as a longer cycle by repetition.
We introduce the following vocabulary to aid in understanding the examples below.
\begin{definition} \label{def:degeneracy}
    Let $K$ be a cycle set.
    \begin{enumerate}
        \item A cycle $x \in K_m$ is \emph{a degeneracy of} a cycle $y \in K_n$ if $x = y\modmap{m}{n}$.
        \item A cycle $x \in K_m$ is \emph{non-degenerate} if the only solution to the equation $x = y\modmap{m}{n}$ is $n = m$ and $y = x$.
        \item The \emph{minimal length} of a cycle $x \in K_n$ is
        \[ \min \{ m \geq 1 \mid \text{$x = y\modmap{m}{n}$ for some $y \in K_m$} \}. \]
    \end{enumerate}
\end{definition}

\begin{example} \label{ex:cycle-set-basic}
    % We describe two cycle sets explicitly,
    Consider the following two cycle sets:
    \begin{enumerate}
        \item Let $K$ be given by $K_n := \{ \ast_n \}$ for all $n$, with all operators the unique maps.
        Every cycle is a degeneracy of $\ast_1$, and $\ast_1$ is the only non-degenerate cycle; every cycle has minimal length $1$.
        \item Let $K$ be given by
        \[ K_n := \begin{cases} \{ u_n, v_n \} & \text{if $n$ is even,} \\ \varnothing & \text{if $n$ is odd,} \end{cases} \]
        with $K(\modmap{m}{n})$ being the identity function for $n$ even dividing $m$, and with $K(\autmap{n}{1})$ being the transposition swapping $u$ and $v$.
        One checks the relations of \cref{cycle-cat-presentation} directly: $(\autmap{n}{1})^n$ is the identity because $n$ is even.
        Here, $u_2$ and $v_2$ are non-degenerate $2$-cycles, whereas all other cycles are degenerate.
        % , and they are exchanged by the $(\ZZ/2)$-action, so they lie in a single orbit.
    \end{enumerate}
    % Item (2) already illustrates that a cycle set need not have any $1$-cycles at all, even though every set $K_n$ with $n$ even is non-empty; this is the sort of behaviour that \cref{characterization-digraph-attractors} will show cannot occur for the attractors of a digraph.
\end{example}

For a fixed $n \geq 1$, the set of $n$-cycles $K_n$ admits an action of the cyclic group $(\ZZ / n, +)$ by the assignment
\[ i \cdot x := x\autmap{n}{i}. \]
This extends to a functor taking values in the category of $(\ZZ / n)$-sets, which we denote by $\operatorname{ev}_n \from \cySet \to \ZnSet{n}$.

This functor is an instance of \cref{adjoint-triple}: by \cref{graph-maps-between-cycles}, the full subcategory of $\Cy$ spanned by the single object $C_n$ has endomorphism monoid $\Aut(C_n) \cong \ZZ / n$; that is, it is the one-object category $\BZ{n}$.
Writing $\iota_n \from \BZ{n} \ito \Cy$ for this inclusion, we have $\operatorname{ev}_n = \iota_n^*$, and so \cref{adjoint-triple} supplies an adjoint triple
\[ (\iota_n)_! \adj \operatorname{ev}_n \adj (\iota_n)_* . \]
Unwinding the coend and end formulas gives the following explicit descriptions.

\begin{proposition}\label{ev-admits-adjoints} \leavevmode
    \begin{enumerate}
        \item The functor $\operatorname{ev}_n$ admits a left adjoint $L \from \ZnSet{n} \to \cySet$ whose set of $k$-cycles is given by:
        \[ L(X)_k := \big( \{ (x, \varphi) \mid x \in X, \ \varphi \from C_k \to C_n \} \big) \big/ \sim, \]
        where $\sim$ is generated by $(i \cdot x, \varphi) \sim (x, \autmap{n}{i} \varphi)$ for all $i \in \ZZ/n$.
        For $\psi \from C_m \to C_k$, we define $[x, \varphi]\psi := [x, \varphi \psi]$, which is well defined on equivalence classes.
        \item The functor $\operatorname{ev}_n$ admits a right adjoint $R \from \ZnSet{n} \to \cySet$ whose $k$-cycles are given by
        \[ R(X)_k := \begin{cases}
            \{ x \in X \mid k \cdot x = x \} & \text{ if $k$ divides $n$} \\
            \{ \ast \} & \text{otherwise}.
        \end{cases} \]
        For $m, k \geq 1$ with $k$ dividing $m$ dividing $n$, we define $x \modmap{m}{k}$ by $x\modmap{m}{k} := x$ and $x\autmap{k}{1}$ by $x\autmap{k}{1} := 1 \cdot x$.
    \end{enumerate}
\end{proposition}
\begin{proof}
    Both statements follow from \cref{adjoint-triple} applied to $\iota_n$.
    For (1), the coend formula gives
    \[ L(X)_k \cong \int^{\ast \in \BZ{n}} \Cy(C_k, C_n) \times X, \]
    and a coend over a one-object groupoid is precisely the quotient of $\Cy(C_k, C_n) \times X$ by the diagonal action, which is the stated equivalence relation.
    For (2), the end formula gives $R(X)_k \cong \operatorname{Hom}_{\ZnSet{n}}\big( \Cy(C_n, C_k), X \big)$, the set of $(\ZZ/n)$-equivariant maps.
    By \cref{graph-maps-between-cycles}, the $(\ZZ/n)$-set $\Cy(C_n, C_k)$ is empty unless $k$ divides $n$, in which case it is the transitive $(\ZZ/n)$-set $\ZZ/k$; equivariant maps out of it are exactly the elements of $X$ fixed by $k$.
\end{proof}
% \begin{remark} \label{rmk:left-adjoint-quotient}
%     The quotient in item (1) is essential and easy to overlook: without it, there may fail to be a unit morphism $X \mapsto L(X)_n$, since $\operatorname{ev}_n L(X)$ would then be $\ZZ/n \times X$ rather than $X$, and the unit of the adjunction would fail to be an isomorphism.
%     Only the right adjoint $R$ is used in what follows.
% \end{remark}
\begin{corollary}\label{ev-preserves-limits}
    The functor $\operatorname{ev}_n \from \cySet \to \ZnSet{n}$ preserves limits and colimits. \qed
\end{corollary}

\subsection{Attractors of a digraph} \label{subsec:attractors}

Having developed cycle sets on their own terms, we now connect them to digraphs.
From the inclusion $\Cy \ito \DiGraph$, we obtain the \emph{attractor functor} $A \from \DiGraph \to \cySet$ by mapping out.
\begin{definition} \label{attractor-functor}
    The \emph{attractor functor} $A \from \DiGraph \to \cySet$ is the functor which sends a digraph $G$ to the cycle set $A(G)$ whose
    \begin{itemize}
        \item $n$-cycles are graph maps $C_n \to G$, i.e.
        \[ A(G)_n := \{ f \from C_n \to G \}; \]
        and whose
        \item action on a morphism $\varphi \from C_m \to C_n$ is given by pre-composition, i.e. given $f \in A(G)_n$ and $\varphi \from C_m \to C_n$, we define $f\varphi \in A(G)_m$ by
        \[ f\varphi := f \circ \varphi. \]
    \end{itemize}
    On a graph map $\gamma \from G \to H$, the induced morphism $A(\gamma) \from A(G) \to A(H)$ is given by post-composition, $f \mapsto \gamma \circ f$.
\end{definition}

\begin{remark} \label{rmk:attractor-restricted-representable}
    It is worth stating that $A(G)$ can be identified with the representable presheaf $\DiGraph(\mathord{-}, G)$ restricted along the inclusion $\Cy \ito \DiGraph$, i.e.
    \[ A(G) = \DiGraph(\mathord{-}, G) \big|_{\Cy} . \]
    This observation can be used to immediately recognize why $A(G)$ is functorial in the variable $G$.
    Readers familiar with simplicial sets (or other well-studied presheaf categories) may recognize the attractor functor as a ``nerve''-type construction (see \cref{rmk:nerve-realization}).
    % In the language of \cref{adjoint-triple}, $A$ is the \emph{nerve} associated to the functor $\Cy \ito \DiGraph$.
\end{remark}

The attractor functor admits a left adjoint defined via left Kan extension, which we refer to as the \emph{digraph realization} functor.
\begin{definition} \label{digraph-realization}
    The \emph{digraph realization} functor is the functor $|\mathord{-}| \from \cySet \to \DiGraph$ defined by left Kan extension of the inclusion $\Cy \ito \DiGraph$ along the Yoneda embedding $\Cy \ito \cySet$.
    \[ \begin{tikzcd}
        \Cy \ar[r, hook] \ar[d, hook] & \DiGraph \\
        \cySet \ar[ur, dotted, "|\mathord{-}|"']
    \end{tikzcd} \]
\end{definition}

\begin{remark} \label{rmk:nerve-realization}
    The pair $\big( |\mathord{-}| \adj A \big)$ is an instance of the \emph{nerve--realization} paradigm: any functor $u \from \cat{A} \to \cat{B}$ from a small category into a cocomplete category induces a functor $\cat{B}(u(\mathord{-}), \mathord{-}) \from \cat{B} \to \Set^{\cat{A}^\op}$ with left adjoint given by left Kan extension along the Yoneda embedding.
    Taking $u$ to be $\Cy \ito \DiGraph$ recovers exactly \cref{attractor-functor,digraph-realization}.
    The familiar example is $u \from \Delta \to \mathsf{Top}$, which yields the singular simplicial set and geometric realization adjunction; the notation $|\mathord{-}|$ is chosen to evoke this.
    % In particular, the adjunction below requires no separate proof.
\end{remark}

The colimit defining $|K|$ may be expressed as a quotient of a disjoint union of cycle graphs.
Explicitly, the digraph $|K|$ is isomorphic to the quotient
\[ |K| \cong \Big{(} \coprod\limits_{n \geq 1} \coprod\limits_{x \in K_n} C_n \Big{)} \Big{/} \sim ,  \]
where $\sim$ is the equivalence relation generated by identities $\big( m, x\varphi, i \big) \sim \big( n, x, \varphi(i) \big)$ for all $\varphi \from C_m \to C_n$ in $\Cy$ and $x \in K_m$.

\begin{remark} \label{rmk:realization-explicit}
    Let us unwind this equivalence relation, since $|K|$ is otherwise rather opaque.
    One starts with one copy of the cycle graph $C_n$ for each $n$-cycle $x \in K_n$ (typically, a cycle set contains infinitely many cycles) and then glues.
    By \cref{graph-maps-between-cycles}, every $\varphi \from C_m \to C_n$ factors as $\autmap{n}{i}\modmap{m}{n}$, so the relation is generated by two kinds of identification, corresponding to the two families of generators:
    \begin{itemize}
        \item \emph{degeneracies collapse}: taking $\varphi = \modmap{m}{n}$, the copy of $C_m$ indexed by the degenerate cycle $x\modmap{m}{n}$ is glued onto the copy of $C_n$ indexed by $x$, wrapping around it $m/n$ times.
        Thus the infinitely many degeneracies of $x$ contribute no new vertices or edges.
        \item \emph{rotations of a cycle collapse}: taking $\varphi = \autmap{n}{i}$, the copy of $C_n$ indexed by $x\autmap{n}{i}$ is identified with the copy indexed by $x$, via a rotation by $i$.
        Thus an entire orbit of the $(\ZZ/n)$-action contributes only one cycle graph.
    \end{itemize}
    The upshot is that, after gluing, $|K|$ has exactly one copy of $C_k$ for each non-degenerate orbit of $K_k$ --- at least when $K$ is well behaved, which is precisely the content of \cref{characterization-digraph-attractors}(3).
\end{remark}

\begin{proposition} \label{attractor-functors-right-adj-preserves-lims}
    The attractor functor $A \from \DiGraph \to \cySet$ is right adjoint to the digraph realization functor $|\mathord{-}| \from \cySet \to \DiGraph$.
    In particular, the attractor functor preserves limits. \qed
\end{proposition}
Moreover, the attractor functor preserves coproducts since the cycle graphs $C_n$ are connected and coproducts in $\cySet$ are computed level-wise.
\begin{proposition} \label{attractor-functor-preserves-coprods}
    The attractor functor $A \from \DiGraph \to \cySet$ preserves coproducts. \qed
\end{proposition}

\Cref{attractor-functor-preserves-coprods} naturally raises the question of whether $A$ preserves all colimits, and hence (because $\cySet$ and $\DiGraph$ are locally presentable) whether it admits a right adjoint.
It does not.

\begin{proposition} \label{attractor-no-right-adjoint}
    The attractor functor $A \from \DiGraph \to \cySet$ does not preserve pushouts, and consequently does not admit a right adjoint.
\end{proposition}
\begin{proof}
    Let $P$ denote the digraph with two vertices $u, v$ and a single edge $u \to v$, and let $D$ denote the digraph with one vertex and no edges.
    Consider the pushout of digraphs
    \[ \begin{tikzcd}
        D \sqcup D \ar[r, "{(u, v)}"] \ar[d] & P \ar[d] \\
        D \ar[r] & C_1
    \end{tikzcd} \]
    in which the left map collapses the two points and the right-hand vertical map identifies $u$ with $v$; the pushout is the digraph $C_1$, a single looped vertex.
    Applying $A$, the digraphs $D \sqcup D$, $D$, and $P$ contain no cycles at all, so $A$ sends each of them to the initial cycle set, which is empty at every level.
    The pushout of the resulting square is therefore again the initial cycle set, whereas $A(C_1)$ is non-empty --- it has a $1$-cycle.
    Hence the square is not preserved.
    % For the second claim, a right adjoint would force $A$ to preserve all colimits, in particular pushouts.
\end{proof}

\begin{remark} \label{rmk:right-adjoint-on-DDS}
    The counterexample above leaves open the more interesting question of whether the composite $AS \from \DDS \to \cySet$ admits a right adjoint, since the digraphs $D$ and $P$ used above are not state space graphs (their vertices have no outgoing edge), and $\DDS$ is closed in $\DiGraph$ under neither of the colimits used.
    We do not know the answer, and record it as an open question in \cref{sec:conclusion}.
\end{remark}
% Writing this colimt as
% Explicitly, this functor sends a cycle set $K$ to the digraph which 
% \fxnote{Give an explicit description of this functor.}

\begin{example} \label{ex:attractor-of-a-digraph}
    Let $G$ be the digraph
    
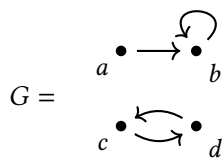
\[ \begin{tikzpicture}
        \node[coordinate] (C1center) {};
        \path (C1center) -- +(0.5, 0) node[vertex, label={[label distance=-1ex]below right:$b$}] (C1R) {};
        \path (C1center) -- +(-0.5, 0) node[vertex, label={[label distance=-1ex]below left:$a$}] (C1L) {};

        \path (C1center) -- +(0, -1) node[coordinate] (C2center) {};
        \path (C2center) -- +(-0.5, 0) node[vertex, label={[label distance=-1ex]below left:$c$}] (C2L) {};
        \path (C2center) -- +(0.5, 0) node[vertex, label={[label distance=-1ex]below right:$d$}] (C2R) {};

        \draw[thick, ->] (C1L) to (C1R);
        \draw[thick, ->] (C1R) to [out=45, in=125, looseness=6] (C1R);

        \draw[thick, ->] (C2R) to [bend right] (C2L);
        \draw[thick, ->] (C2L) to [bend right] (C2R);

        \path (C1center) -- node[coordinate, pos=0.5] (middle) {} (C2center);
        \path (middle) -- +(-1.25, -0.1) node[anchor=east] {\large $G = $};
    \end{tikzpicture} \]
    The $n$-cycles of $A(G)$ are graph maps $C_n \to G$, which is exactly the data of a cycle $(v_1 \to v_2 \to \dots \to v_n \to v_1)$ of length $n$ in $G$ with possible repeats.
    We will denote such a cycle by $(v_1 v_2 \dots v_n)$ for brevity, and compute the sets $A(G)_n$ for small $n$.
    \begin{itemize}
        \item For $n = 1$, there is one such cycle $(b)$, so 
        \[ A(G)_1 = \{ (b) \}. \]
        \item For $n = 2$, the cycles of length 2 are $(bb)$, $(cd)$, and $(dc)$.
        Thus, 
        \[ A(G)_2 = \{ (bb), (cd), (dc) \}. \]
        The action of $\ZZ / 2$ on $A(G)_2$ swaps $(cd)$ with $(dc)$ and leaves $(bb)$ as a fixed point. 
        \item For $n = 3$, the only cycle of length 3 is 
        \[ A(G)_3 = \{ (bbb) \} \]
        and the action by $\ZZ / 3$ is trivial.
        \item For $n = 4$, the cycles of length 4 are 
        \[ A(G)_4 = \{ (bbbb), (cdcd), (dcdc) \}. \]
        The action of $\ZZ/4$ swaps $(cdcd)$ with $(dcdc)$ and leaves $(bbbb)$ as a fixed point.
    \end{itemize}
    Note that $(b)\modmap{2}{1} = (bb)$ and $(cd)\modmap{4}{2} = (cdcd)$.
    % Note that $(bb)$ is the 2-cycle obtained by applying $A(G)(\modmap{2}{1}) \from A(G)_1 \to A(G)_2$ to the 1-cycle $(b)$.
    In general, the function $A(G)(\modmap{m}{n}) \from A(G)_n \to A(G)_m$ takes a cycle of length $n$ and repeats it until it becomes a cycle of length $m$.
    In this example, the function $A(G)(\modmap{4}{2}) \from A(G)_2 \to A(G)_4$ is a bijection, and is moreover equivariant when viewing the $(\ZZ / 2)$-set $A(G)_2$ as a $(\ZZ / 4)$-set via the surjective group homomorphism $\ZZ / 4 \to \ZZ / 2$.

    Finally, we compute the digraph realization $|A(G)|$.
    Using \cref{rmk:realization-explicit}, we see that $|A(G)|$ has one copy of $C_k$ for each non-degenerate orbit of $A(G)_k$.
    There are exactly two: the orbit of the $1$-cycle $(b)$, and the orbit $\{ (cd), (dc) \}$ of $2$-cycles.
    Hence
    \[ \begin{tikzpicture}
        \node[coordinate] (C1center) {};
        \node[vertex, label={[label distance=-1ex]below:$b$}] (C1A) at (C1center) {};

        \path (C1center) -- +(2, 0) node[coordinate] (C2center) {};
        \path (C2center) -- +(-0.5, 0) node[vertex, label={[label distance=-1ex]below:$c$}] (C2L) {};
        \path (C2center) -- +(0.5, 0) node[vertex, label={[label distance=-1ex]below:$d$}] (C2R) {};

        \path (C1center) -- +(-0.55, 0.15) node[anchor=east] {\large $|A(G)| \cong $};
        % \path (C1center) -- +(1.15, 0.05) node {\large $\sqcup$};

        \draw[thick, ->] (C1A) to [out=45, in=125, looseness=6] (C1A);
        \draw[thick, ->] (C2L) to [bend right] (C2R);
        \draw[thick, ->] (C2R) to [bend right] (C2L);
    \end{tikzpicture} \]
    that is, $|A(G)| \cong C_1 \sqcup C_2$.
    Comparing with $G$, we see that the vertex $a$ and the edge $a \to b$ have disappeared: they lie on no cycle, and the attractor functor sees only the cyclic part of $G$.
    The unit map $\eta_{A(G)} \from A(G) \to A|A(G)|$ is nevertheless an isomorphism, in accordance with \cref{characterization-digraph-attractors}, since $G$ and $C_1 \sqcup C_2$ have the same cycles.
\end{example}

\begin{example}[Running example, continued] \label{ex:running-attractors}
    We compute the attractors of the system $(X, f)$ of \cref{ex:running}, whose state space graph was drawn in \cref{ex:running-state-space}.
    An $n$-cycle of $AS(X,f)$ is a closed walk of length $n$ in $S(X,f)$, so
    \[ AS(X,f)_1 = \{ (a) \}, \qquad AS(X,f)_2 = \{ (aa), (bc), (cb) \}, \qquad AS(X,f)_3 = \{ (aaa) \}, \]
    and in general $AS(X,f)_n$ consists of the degeneracy $(a \dots a)$ together with, when $n$ is even, the two degeneracies of $(bc)$ and $(cb)$.
    Evaluating at $n = 2$ gives the $(\ZZ/2)$-set
    \[ \operatorname{ev}_2 AS(X,f) = \{ (aa), (bc), (cb) \}, \]
    on which the generator of $\ZZ/2$ fixes $(aa)$ and swaps $(bc)$ with $(cb)$; thus there are two orbits, of sizes $1$ and $2$.
    The non-degenerate cycles are $(a)$, $(bc)$, and $(cb)$, with minimal lengths $1$, $2$, and $2$ respectively; every other cycle is a degeneracy of one of these.
    The digraph realization is $|AS(X,f)| \cong C_1 \sqcup C_2$, which in this case is isomorphic to $S(X,f)$ itself, since every state of $(X,f)$ already lies on a cycle.
\end{example}

If a digraph is the state space of a discrete dynamical system $(X, f)$, an $n$-cycle in $AS(X, f)$ can be identified with a morphism of dynamical systems $(\ZZ / n, \succ) \to (X, f)$.
\begin{proposition} \label{cycles-of-a-dynamial-system-correspondence}
    Let $(X, f)$ be a dynamical system.
    For any $n \geq 1$, there is a one-to-one correspondence between $n$-cycles of $AS(X, f)$ and morphisms of dynamical systems $(\ZZ / n, \succ) \to (X, f)$.
\end{proposition}
\begin{proof}
    By definition, an $n$-cycle of $AS(X, f)$ is a graph map $C_n \to S(X, f)$.
    The cycle graph $C_n$ is the state space of the discrete dynamical system $(\ZZ / n, \succ)$.
    The desired correspondence thus follows since the state space functor is full and faithful (\cref{state-space-functor-full-faithful}).
\end{proof}

\subsection{Degenerate cycles in the attractors of a digraph} \label{subsec:degenerate}

Cycle sets which arise as the cycle set of attractors of a digraph satisfy certain convenient properties, which we now establish, using the vocabulary of \cref{def:degeneracy}.
% \begin{example} \label{ex:generators-of-digraph-cycles}
%     Returning to the digraph in \cref{ex:attractor-of-a-digraph}, the cycles $(b)$, $(bb)$, $(bbb)$, and $(bbbb)$ are all degeneracies of $(b)$, which is non-degenerate.
%     The cycles $(cd)$ and $(dc)$ are non-degenerate, and $(cdcd), (dcdc)$ are degeneracies on them, respectively.
%     % the cycle $(b)$ is a generator of $(b)$, $(bb)$, $(bbb)$, and $(bbbb)$; this shows the generating length of each of these cycles is 1.
%     % The cycle $(cd)$ is a generator of $(cd)$ and $(cdcd)$, hence these cycles have generating length $2$.
% \end{example}
% When studying attractors of a digraph, we are often interested in the injective $n$-cycles, i.e.\ sequences $(v_1 \dots v_n)$ that are not repeats of a smaller-length cycle.
% Note that, in \cref{ex:attractor-of-a-digraph}, the minimal length of each cycle coincides with the size of its orbit.
% This fact holds for the attractors of any digraph.
\begin{lemma} \label{attractor-cycle-autmap-degeneracy}
    Let $G$ be a digraph and $f \in A(G)_n$ be an $n$-cycle.
    If $f\autmap{n}{k} = f$ for some $k$ dividing $n$ then there exists $g \in A(G)_{k}$ such that $f = g\modmap{n}{k}$.
\end{lemma}
\begin{proof}
    The equality $f\autmap{n}{k} = f$ unfolds to an equality
    \[ f(i) = \textstyle f(i +_{\ZZ / n} k) \]
    for all $i \in C_n$.
    With this, we may define a graph map $g \from C_{k} \to G$ by $g(i) := f(i)$.
    By construction, we have that $g \modmap{n}{k} = f$, as desired.
\end{proof}
The next theorem is the technical heart of this subsection.
Recall from \cref{rmk:realization-explicit} that a cycle set carries a great deal of redundant data: each attractor of a system appears not once but infinitely often, as its own degeneracies, and each of those appears once for every rotation.
To count attractors, it would be beneficial to know that this redundancy is bookkeeping and nothing more, i.e.\ that each cycle can be traced back to a well-defined ``genuine'' attractor.
That is exactly what \cref{attractor-digraph-injetive-unique-degeneracies} provides: item (2) says every cycle has a unique non-degenerate ancestor, and item (3) identifies the length of that ancestor with the size of the orbit, so that the genuine attractors consisting of $k$ distinct states (i.e.\ those represented by a non-degenerate $k$-cycle) are counted by the non-degenerate orbits of $A(G)_k$.
It is what makes the counting formula of \cref{non-deg-cycle-from-AB} possible, and it is used again in the proof of the decomposition theorem (\cref{modularity-theorem}) to make sense of the ``distinguished element'' of an orbit.
Crucially, the theorem is false for general cycle sets: \cref{ex:A-not-B,ex:B-not-A,ex:A-without-unique-degens} exhibit cycle sets in which it fails, and \cref{characterization-digraph-attractors} shows that this failure is precisely what distinguishes a general cycle set from the cycle set of attractors of a digraph.

\begin{theorem} \label{attractor-digraph-injetive-unique-degeneracies}
    Let $G$ be a digraph.
    \begin{enumerate}
        \item For any morphism $\varphi \from C_m \to C_n$, the function $A(G)(\varphi) \from A(G)_n \to A(G)_m$ is injective.
        \item Every cycle in $A(G)$ is a degeneracy of a unique non-degenerate cycle.
        \item The minimal length of an $n$-cycle $f \in A(G)_n$ is equal to the size of the orbit containing $f$.
    \end{enumerate}
\end{theorem}
\begin{proof}
    Item (1) follows since every map between cycles $\varphi \from C_m \to C_n$ is surjective, hence an epimorphism.

    For (2), fix non-degenerate cycles $g \in A(G)_m$ and $h \in A(G)_n$, and suppose $f \in A(G)_k$ is a $k$-cycle such that $f = g\modmap{k}{m} = h\modmap{k}{n}$.
    If $n = m$ then this follows from (1), since $A(G)(\modmap{k}{m})$ is injective.
    Thus, it suffices to show that $n = m$.

    Using the identities in \cref{cycle-cat-presentation}, we calculate
    \[ f\autmap{k}{m} = g\modmap{k}{m}\autmap{k}{m} = g\autmap{m}{m}\modmap{k}{m} = g\modmap{k}{m} = f, \]
    and an analogous calculation gives $f\autmap{k}{n} = f$.
    This shows that $m$ and $n$ are elements of the stabilizer subgroup $\operatorname{Stab}(f) \subseteq \ZZ / k$ of $f$.
    As $\operatorname{Stab}(f)$ is a subgroup of a cyclic group, it follows that the greatest common divisor $\gcd(m, n)$ is an element of $\operatorname{Stab}(f)$.
    That is, writing $d := \gcd(m, n)$, we have an equality $f\autmap{k}{d} = f$.

    We now compute
    \[ g\autmap{m}{d}\modmap{k}{m} = g\modmap{k}{m}\autmap{k}{d} = f\autmap{k}{d} = f = g\modmap{k}{m}, \]
    By (1), the function $A(G)(\modmap{k}{m})$ is injective, which implies $g\autmap{m}{d} = g$.
    An analogous computation gives that $h\autmap{n}{d} = h$.
    If either $d \neq m$ or $d \neq n$ then one of $g$ or $h$ is not non-degenerate by \cref{attractor-cycle-autmap-degeneracy}.
    This would be a contradiction, therefore $n = m$.

    For (3), the size of the orbit containing $f$ is $n / m$, where $m$ is the order of the stabilizer subgroup of $f$.
    As a subgroup of a cyclic group, the stabilizer of $f$ is generated by the element $(n/m) \in \ZZ / n$, and this element is the minimal positive integer such that $(n/m) \cdot f = f$.
    By \cref{attractor-cycle-autmap-degeneracy}, there exists $g \in A(G)_{n/m}$ such that $f$ is a degeneracy of $g$.
    By minimality of $(n/m)$, it follows that $g$ is non-degenerate.
    By (2), we conclude that the minimal length of $f$ is $n/m$.
\end{proof}

\subsection{Characterizing cycle sets coming from digraphs} \label{subsec:characterization}

This subsection is devoted to a single result, \cref{characterization-digraph-attractors}, and since its proof occupies several pages we begin by saying what it asserts and why one should want it.

The attractor functor $A \from \DiGraph \to \cySet$ is not surjective on objects: it is easy to exhibit cycle sets that are not the cycle set of attractors of any digraph (see \cref{ex:A-not-B,ex:B-not-A,ex:not-AB} below).
The question this subsection answers is: \emph{which ones are?}
The answer is that there are exactly two obstructions, both concerning degeneracies, and both already met in \cref{subsec:degenerate}:
\begin{itemize}
    \item two distinct cycles must not become equal upon degenerating (Property A below); and
    \item a cycle fixed by a rotation must actually \emph{be} a degeneracy (Property B below).
\end{itemize}
A digraph's attractors satisfy both, by \cref{attractor-digraph-injetive-unique-degeneracies,attractor-cycle-autmap-degeneracy}, and \cref{characterization-digraph-attractors} shows that these two conditions are also sufficient.
The theorem gives three further equivalent formulations, of which the most useful is that $K$ decomposes as a coproduct of representables indexed by its non-degenerate orbits --- in effect, that $K$ is freely built from its genuine attractors.
From this description, a formula for counting the non-degenerate cycles of a digraph can be deduced (\cref{non-deg-cycle-from-AB}).

The reader who wants only the statement may skip to \cref{characterization-digraph-attractors}; what precedes it is a sequence of examples establishing that the two properties are genuinely independent and genuinely necessary.

\begin{definition} \label{def:property-AB}
    Let $K$ be a cycle set.
    \begin{enumerate}
        \item We say that $K$ satisfies \emph{Property A} if, for all $m, n \geq 1$ with $n$ dividing $m$, the function $K(\modmap{m}{n}) \from K_n \to K_m$ is injective.
        \item We say $K$ satisfies \emph{Property B} if, for any $n$-cycle $x \in K_n$, if $x\autmap{n}{k} = x$ for some $k$ dividing $n$ then $x = y\modmap{n}{k}$ for some $y \in K_k$.
    \end{enumerate}
\end{definition}
We note that Properties A and B are invariant under isomorphism.

The next three examples are, we think, the most instructive part of this subsection, so we flag explicitly what they are for.
Each exhibits a cycle set that is \emph{not} the attractor cycle set of any digraph, and together they show that Properties A and B are logically independent: \cref{ex:A-not-B} satisfies A but not B, \cref{ex:B-not-A} satisfies B but not A, and \cref{ex:not-AB} satisfies neither.
In particular, neither property may be dropped from \cref{characterization-digraph-attractors}.
It is worth pausing over each to see concretely how a cycle set can fail to come from a digraph; note in particular that, in the first two examples, the digraph realization $|K|$ is the same digraph in both cases, so that $K$ cannot be recovered from $|K|$, which is precisely the failure of the unit map $\eta_K$ to be an isomorphism.

\begin{example} \label{ex:A-not-B}
    Define a cycle set $K$ by
    \[ K_n := \begin{cases}
        \{ \ast_n \} & \text{if $n$ is even} \\
        \varnothing & \text{if $n$ is odd}.
    \end{cases} \]
    The operators $K(\modmap{m}{n})$ and $K(\autmap{n}{i})$ are defined tautologically.
    Then, $K$ satisfies Property A (in fact, it satisfies the conclusion of \cref{unique-degens-from-AB}) but fails Property B, since $\ast_2 \autmap{2}{1} = \ast_2$, but $\ast_2$ is non-degenerate.

    The digraph realization $|K|$ is the digraph with a single looped vertex.
\end{example}
\begin{example} \label{ex:B-not-A}
    Define a cycle set $K$ by
    \[ K_n := \begin{cases}
        \{ \ast_n \} & \text{if } n \geq 2 \\
        \{ 0, 1 \} & \text{if } n = 1.
    \end{cases} \]
    The operators $K(\modmap{m}{n})$ and $K(\autmap{n}{i})$ are defined tautologically.
    Then, $K$ satisfies Property B but not Property A, since $0\modmap{2}{1} = \ast_2 = 1\modmap{2}{1}$ but $0 \neq 1$.

    The digraph realization $|K|$ is again the digraph with a single looped vertex.
\end{example}
\begin{example} \label{ex:not-AB}
    An example of a cycle set which fails both Property A and Property B is the coproduct of the cycle set in \cref{ex:A-not-B} with the cycle set in \cref{ex:B-not-A}.
\end{example}

The proof of item (2) in \cref{attractor-digraph-injetive-unique-degeneracies} may be read verbatim as a proof of the following fact.
\begin{proposition} \label{unique-degens-from-AB}
    Suppose $K$ is a cycle set satisfying Properties A and B.
    Then, every cycle in $K$ is a degeneracy of a unique non-degenerate cycle. \qed
\end{proposition}
Note that if every cycle is a degeneracy of a unique non-degenerate cycle then Property A holds.
We present an example of a cycle set satisfying Property A, but for which the conclusion of \cref{unique-degens-from-AB} does not hold.
\begin{example} \label{ex:A-without-unique-degens}
    Define a cycle set $K$ by
    \[ K_n := \begin{cases}
        \{ \ast_n \} & \text{if } n \geq 2 \\
        \varnothing & \text{if } n = 1.
    \end{cases} \]
    The operators $K(\modmap{m}{n})$ and $K(\autmap{n}{i})$ are defined tautologically.
    Then, $K$ satisfies Property A, but not every cycle is a degeneracy of a unique non-degenerate cycle, since $\ast_6 = \ast_2 \modmap{6}{2} = \ast_3 \modmap{6}{3}$, but both $\ast_2$ and $\ast_3$ are non-degenerate.
    Note that the minimal length of $\ast_6$ is 2, yet $\ast_3$ is a non-degenerate 3-cycle satisfying $\ast_6 = \ast_3 \modmap{6}{3}$.
    
    The digraph realization $|K|$ is the digraph with a single looped vertex.
\end{example}

Before proving that Properties A and B characterize the cycle sets which are cycle sets of attractors of some digraph, we prove a short lemma regarding the representable cycle sets.
For $k \geq 1$, we write $\Cyhat{k}$ for the cycle set represented by $C_k$, i.e.\ the cycle set whose $n$-cycles are the morphisms $(\Cyhat{k})_n := \Cy(C_n, C_k)$; for an $n$-cycle $\varphi \in (\Cyhat{k})_n$, the operators $\varphi \modmap{m}{n}$ and $\varphi \autmap{n}{i}$ are defined by composition of morphisms.

From the realization-attractor adjunction, we extract the unit map
\[ \eta_K \from K \to A|K| \]
which is a map of cycle sets natural in the variable $K$.
Using the explicit description of $|K|$ (before \cref{attractor-functors-right-adj-preserves-lims}), a formula for this map is given by
\[ (\eta_K)_n(x) := i_{(n, x)}, \]
where $i_{(n, x)}$ is the coproduct inclusion $C_n \to |K|$ indexed by the $n$-cycle $x \in K_n$.

\begin{lemma} \label{unit-at-representable-iso}
    For $k \geq 1$, the unit map instantiated at the representable $\Cyhat{k}$ is an isomorphism.
\end{lemma}
\begin{proof}
    The colimit formula for left Kan extensions gives a description of $|\Cyhat{k}|$ as
    \[ |\Cyhat{k}| \cong \colim \left( \Cy \slice C_k \to \Cy \ito \DiGraph \right). \]
    Let $D \from \Cy \to \DiGraph$ denote the diagram on the right.
    The diagram $D$ is indexed by the slice category $\Cy \slice {C_k}$, which has a terminal object $\id[C_k]$.
    Thus, the colimit inclusion $\lambda_{\id[C_k]} \from C_k \to |\Cyhat{k}|$ indexed by $\id[C_k]$ is an isomorphism.
    Using the explicit description of $|\Cyhat{k}|$ as a coproduct, we see that for any $\varphi \from C_n \to C_k$, the colimit inclusion $\lambda_{\varphi} \from C_n \to |\Cyhat{k}|$ indexed by $\varphi$ is exactly the coproduct inclusion $i_{(n, \varphi)}$.
    With this, we compute
    \[ (\eta_{\Cyhat{k}})_n(\varphi) = i_{(n, \varphi)} = \lambda_{\varphi} = \lambda_{\id[C_k]} \circ D(\varphi) = \lambda_{\id[C_k]} \circ \varphi. \]
    Since $\lambda_{\id[C_k]}$ is an isomorphism, it follows that the explicit formula for the unit defines a bijection.
\end{proof}

We will also use the following elementary observation, valid for any adjunction: units are compatible with colimits.
% Stating it separately spares us an explicit diagram chase later on.

\begin{lemma} \label{unit-distributes-over-colimits}
    Let $L \adj R$ be an adjunction between categories $\cat{C}$ and $\cat{D}$, with unit $\eta$, and let $X \cong \colim_{i} X_i$ be a colimit in $\cat{C}$.
    Then the triangle
    \[ \begin{tikzcd}[column sep = large]
        \colim_i X_i \ar[r, "\colim_i \eta_{X_i}"] \ar[rd, "\eta_{\colim_i X_i}"'] & \colim_i RL X_i \ar[d, "\theta"] \\
        {} & RL \colim_i X_i
    \end{tikzcd} \]
    commutes, where $\theta$ is the canonical comparison map.
    Consequently, if each $\eta_{X_i}$ is an isomorphism and $\theta$ is an isomorphism, then $\eta_{\colim_i X_i}$ is an isomorphism.
\end{lemma}
\begin{proof}
    Both composites are maps out of a colimit, so it suffices to check that they agree after precomposition with each colimit inclusion $\lambda_i \from X_i \to \colim_i X_i$.
    The lower composite gives $\eta_{\colim X_i} \circ \lambda_i$, and the upper composite gives $\theta \circ (\colim \eta) \circ \lambda_i = RL(\lambda_i) \circ \eta_{X_i}$ by the definition of $\theta$.
    These agree by naturality of $\eta$ applied to $\lambda_i$.
    The final claim is immediate, since a composite of isomorphisms is an isomorphism.
\end{proof}

\begin{remark} \label{rmk:unit-colimit-application}
    We emphasize that \cref{unit-distributes-over-colimits} holds for a colimit of \emph{any} shape and requires no hypothesis on the morphisms involved.
    In our application, $L = |\mathord{-}|$ and $R = A$, the colimit is a coproduct, and $\theta$ is an isomorphism because $|\mathord{-}|$ preserves coproducts (being a left adjoint) and $A$ preserves coproducts (\cref{attractor-functor-preserves-coprods}).
\end{remark}

For a cycle set $K$ and $k \geq 1$, let $\Orbnd K_k$ denote the set
\[ \Orbnd K_k := \{ x \in K_k \mid x \text{ is non-degenerate} \} \, \big/ \sim, \]
where $\sim$ is the relation given by $x \sim y$ if $x$ and $y$ are in the same orbit under the $(\ZZ / k)$-action.
Note that since $\modmap{m}{n}\autmap{m}{1} = \autmap{n}{1}\modmap{m}{n}$, it follows that $x$ is non-degenerate if and only if any cycle in the orbit of $x$ is non-degenerate.
Hence, we refer to the set $\Orbnd K_k$ as the \emph{non-degenerate orbits} of $K_k$.
\begin{theorem} \label{characterization-digraph-attractors}
    For a cycle set $K$, the following are equivalent:
    \begin{enumerate}
        \item $K$ is isomorphic to the cycle set of attractors of some digraph;
        \item the unit map $\eta_K \from K \to A|K|$ of the realization-attractor adjunction is an isomorphism;
        \item $K$ is isomorphic to a coproduct of representables
        \[ K \cong \coprod\limits_{k \geq 1} \coprod\limits_{[x] \in \Orbnd K_k} \Cyhat{k} \]
        indexed by the non-degenerate orbits of $K$;
        \item $K$ satisfies Properties A and B.
    \end{enumerate}
\end{theorem}
\begin{proof}
    We prove the implications $(2) \implies (1) \implies (4) \implies (3) \implies (2)$ in order.

    \textit{(2) implies (1).}
    This is immediate: if $\eta_K$ is an isomorphism then $K \cong A|K|$, so $K$ is the cycle set of attractors of the digraph $|K|$.

    \textit{(1) implies (4).}
    Property A holds by item (1) of \cref{attractor-digraph-injetive-unique-degeneracies}.
    Property B holds by \cref{attractor-cycle-autmap-degeneracy}.

    \textit{(4) implies (3).}
    The $n$-cycles in the coproduct $\coprod \coprod \Cyhat{k}$ are given by tuples $(k, [x], \varphi)$ where $k \geq 1$, $[x]$ is a non-degenerate orbit of $K_k$, and $\varphi$ is a morphism $C_n \to C_k$ in $\Cy$, with operators acting by post-composition in the last coordinate:
    \[ (k, [x], \varphi)\modmap{m}{n} := (k, [x], \varphi \modmap{m}{n}), \qquad (k, [x], \varphi)\autmap{n}{1} := (k, [x], \varphi \autmap{n}{1}). \]
    For each non-degenerate orbit $[x] \in \Orbnd K_k$ choose a distinguished representative $x \in K_k$, and define
    \[ \Theta \from \coprod\limits_{k \geq 1} \coprod\limits_{[x] \in \Orbnd K_k} \Cyhat{k} \longrightarrow K, \qquad \Theta_n(k, [x], \varphi) := x\varphi . \]
    This commutes with the operators by construction.
    % Rather than verifying that $\Theta$ is bijective in each degree by hand, we exhibit its inverse directly, which is where Properties A and B are used.
    It remains only to use Properties A and B to contruct an inverse to $\Theta$.

    Let $y \in K_n$ be any $n$-cycle.
    Since $K$ satisfies Properties A and B, \cref{unique-degens-from-AB} applies: there is a \emph{unique} non-degenerate cycle $z \in K_k$, for a uniquely determined $k$, with $y = z \modmap{n}{k}$.
    The cycle $z$ lies in a unique non-degenerate orbit $[z] \in \Orbnd K_k$, whose distinguished representative we have called $x$; thus $z = x \autmap{k}{i}$ for some $i \in \ZZ/k$.
    Moreover this $i$ is unique: if $x\autmap{k}{i} = x\autmap{k}{j}$ then $x \autmap{k}{i - j} = x$, and Property B would then exhibit $x$ as a degeneracy unless $i \equiv j$, contradicting non-degeneracy of $x$.
    We may therefore define
    \[ \Xi_n(y) := \big( k, [z], \autmap{k}{i}\modmap{n}{k} \big), \]
    which is well defined by the two uniqueness statements just established.

    It remains to check that $\Xi$ and $\Theta$ are mutually inverse.
    On the one hand, $\Theta_n(\Xi_n(y)) = x \autmap{k}{i} \modmap{n}{k} = z\modmap{n}{k} = y$.
    On the other, given $(k, [x], \varphi)$, factor $\varphi = \autmap{k}{i}\modmap{n}{k}$ uniquely using \cref{graph-maps-between-cycles}; then $\Theta_n(k,[x],\varphi) = x\autmap{k}{i}\modmap{n}{k}$, and since $x\autmap{k}{i}$ is non-degenerate (being in the orbit of the non-degenerate cycle $x$), the uniqueness in \cref{unique-degens-from-AB} forces $\Xi_n$ to return $(k, [x], \autmap{k}{i}\modmap{n}{k}) = (k,[x],\varphi)$.
    Hence $\Theta$ is a levelwise bijection, i.e.\ an isomorphism of cycle sets.

    \textit{(3) implies (2).}
    This is now a formal consequence of \cref{unit-distributes-over-colimits} applied to the adjunction $|\mathord{-}| \adj A$ and the coproduct in (3).
    Indeed, the unit is an isomorphism at each representable $\Cyhat{k}$ by \cref{unit-at-representable-iso}, and the canonical comparison map
    \[ \coprod\limits_{k \geq 1} \coprod\limits_{[x] \in \Orbnd K_k} A|\Cyhat{k}| \longrightarrow A \Big| \coprod\limits_{k \geq 1} \coprod\limits_{[x] \in \Orbnd K_k} \Cyhat{k} \Big| \]
    is an isomorphism because $|\mathord{-}|$ preserves coproducts, being a left adjoint, and $A$ preserves coproducts by \cref{attractor-functor-preserves-coprods}.
    \Cref{unit-distributes-over-colimits} therefore gives that $\eta$ is an isomorphism at $\coprod \coprod \Cyhat{k}$, and hence at $K$, since the two are isomorphic by (3) and $\eta$ is natural.
    % \[ A \Big| \coprod\limits_{k \geq 1} \coprod\limits_{[x] \in \Orbnd K_k} \Cyhat{k} \Big| \overunderset{A\Phi}{\cong}{\to} A \Big( \coprod\limits_{k \geq 1} \coprod\limits_{[x] \in \Orbnd K_k} |\Cyhat{k}| \Big) \overunderset{\Psi}{\cong}{\to} \coprod\limits_{k \geq 1} \coprod\limits_{[x] \in \Orbnd K_k} A|\Cyhat{k}| . \]
    % By definition, 
    % At level $n$, a formula for $(A\Phi)_n$ is given by
    % \[ (A\Phi)_n(c) := \Phi \circ c. \]
    % From this formula, it follows that the triangle
    % \[ \begin{tikzcd}[sep = large]
    %     {} & \displaystyle A \Big| \coprod\limits_{k \geq 1} \coprod\limits_{[x] \in \Orbnd K_k} \Cyhat{k} \Big| \ar[d, "A\Phi"] \\
    %     \displaystyle \coprod\limits_{k \geq 1} \coprod\limits_{[x] \in \Orbnd K_k} \Cyhat{k} \ar[r, "\coprod \coprod \eta_{\Cyhat{k}}"'] \ar[ur, "\eta_{\coprod \coprod \Cyhat{k}}"] & \displaystyle \coprod\limits_{k \geq 1} \coprod\limits_{[x] \in \Orbnd K_k} A|\Cyhat{k}|
    % \end{tikzcd} \]
    % commutes.
    % To show this, we show the levelwise equality $(A\Phi)_n \circ (\eta_{\coprod \coprod \Cyhat{k}})_n = \coprod \coprod (\eta_{\Cyhat{k}})_n$ by instantiating at an $n$-cycle $(k, [x], \varphi)$.
    % Since this is an equality of graph maps, $C_n \to \coprod \coprod |\Cyhat{k}|$, 
    % \begin{align*} 
    %     \big( (A\Phi \circ \eta_{\coprod \coprod \Cyhat{k}})_n(k, [x], \varphi) \big)(j) = \big( \Phi \circ i_{(n, (k, [x], \varphi))} \big)(j) = \Phi(i_{(n, (k, [x], \varphi))}(i)) 
    % \end{align*}
\end{proof}
From the coproduct description, we obtain a formula for counting the non-degenerate cycles of a cycle set $K$ coming from a digraph.
\begin{corollary} \label{non-deg-cycle-from-AB}
    Suppose $K$ is a cycle set satisfying any of the equivalent conditions of \cref{characterization-digraph-attractors}.
    Then, for any $n \geq 1$, the set of non-degenerate orbits is in bijection with the set complement
    \[ \Orbnd K_n \cong \operatorname{Orb} K_n - \coprod\limits_{\substack{d=1 \\ d \text{ divides } n}}^{n-1} \Orbnd K_d . \]
\end{corollary}
\begin{proof}
    We instantiate the isomorphism in item (3) of \cref{characterization-digraph-attractors} at $n$ to get an isomorphism of $(\ZZ / n)$-sets
    \[ K_n \cong \coprod\limits_{k \geq 1} \coprod\limits_{[x] \in \Orbnd K_k} (\Cyhat{k})_n. \]
    Since taking orbits commutes with coproducts, we obtain a bijection of sets
    \[ \operatorname{Orb} K_n \cong \coprod\limits_{k \geq 1} \coprod\limits_{[x] \in \Orbnd K_k} \operatorname{Orb} (\Cyhat{k})_n \cong \coprod\limits_{k \geq 1} \left( \Orbnd K_k \times \operatorname{Orb} (\Cyhat{k})_n \right) . \]
    % The set $(\Cyhat{k})$
    We use \cref{graph-maps-between-cycles} to give an explicit description of the $(\ZZ / n)$-set as
    \[ (\Cyhat{k})_n \cong \begin{cases}
        \{ \modmap{n}{k}, \autmap{k}{1}\modmap{n}{k}, \dots, \autmap{k}{k-1} \modmap{n}{k} \} & \text{if $k$ divides $n$} \\
        \varnothing & \text{otherwise.}
    \end{cases} \]
    In the non-empty case, the action is given by $1 \cdot \autmap{k}{i}\modmap{n}{k} = \autmap{k}{i}\modmap{n}{k}\autmap{n}{1} = \autmap{k}{i+1 \, \operatorname{mod} k} \modmap{n}{k}$.
    In particular, the action is transitive, hence the previous bijection simplifies to
    \[ \operatorname{Orb} K_n \cong \coprod\limits_{\substack{d = 1 \\ d \text{ divides } n}}^{d=n} \Orbnd K_d. \]
    The desired isomorphism then follows by isolating the $d = n$ term of the coproduct.
\end{proof}

% \subsection{Semi-direct products}

\section{The decomposition theorem}\label{Section:Decomposition}

This section is the natural endpoint of \cref{sec:framework}, where the machinery built there bears fruit.
The goal is a single theorem, \cref{modularity-theorem}, which generalizes the decomposition theorem of \cite{kadelka2023modularity} (quoted as \cref{old-decomposition-theorem} in the introduction) from Boolean networks to arbitrary discrete dynamical systems.
We highlight that this theorem is where our focus lies: it says that the attractors of a Boolean network built as a semi-direct product can be assembled from the attractors of the pieces.
What we will show is that this is an instance of a general principle, namely that a semi-direct product is a certain pullback (\cref{semiprod-pb}) and that the attractor functor, being a composite of right adjoints, preserves it.

This relies on the notion of a semi-direct product of dynamical systems (\cref{semidirect-product}), which is a generalization of the semi-direct product of Boolean networks developed in \cite{kadelka2023modularity}.
For this, we work in some generality, considering dynamical systems not necessarily in the category of sets, but in any category with finite products.

The theorem asserts the existence of a decomposition once a semi-direct product structure is given; the separate and equally important question of how to \emph{find} such a structure is the subject of \cref{sec:wiring}.

\begin{definition} \label{semidirect-product}
    Let $\cat{C}$ be a category with finite products.
    Given morphisms $f \from X \to X$, $g \from E \times Y \to Y$ and $p \from X \to E$, the \emph{semi-direct product} of $f$ and $g$ along $p$ is the morphism $f \semiprod{p} g \from X \times Y \to X \times Y$ defined to be the composite
    \[ \begin{tikzcd}
        X \times Y \ar[rr, dotted, "f \semiprod{p} g"] \ar[rd, "{(\proj_X, p \circ \proj_X, \proj_Y)}"'] & {} & X \times Y \\
        {} & X \times E \times Y \ar[ur, "{(f \circ \proj_X, g \circ \proj_{E \times Y})}"'] & {}
    \end{tikzcd} \] 
\end{definition}
When $\cat{C} = \Set$, the formula for the semi-direct product can be written as
\[ (f \semiprod{p} g)(x, y) = \big( f(x), g(p(x), y) \big). \]

\begin{remark}[Terminology] \label{rmk:cascade-product}
    Constructions of this shape appear under several names in the literature.
    In automata theory and semigroup theory, an endomorphism of $X \times Y$ whose first coordinate depends only on $X$ is called a \emph{cascade product} (or \emph{cascade composition}); it is the basic building block of the Krohn--Rhodes decomposition theorem \cite{krohn-rhodes}, which asserts that every finite automaton divides a cascade of simple groups and flip-flops.
    In dynamics, the same construction is called a \emph{skew product}, with $(X,f)$ the base and $(Y, g)$ the fibre.
    We retain the name ``semi-direct product'' for continuity with \cite{kadelka2023modularity}, where the term is adopted in analogy with the semi-direct product of groups.
\end{remark}

\begin{remark}[The object $E$ is not essential] \label{rmk:E-redundant}
    It is worth observing that the object $E$ carries no information that is not already present in the pair $(p, g)$: for any $f, g, p$ as above we have
    \[ f \semiprod{p} g = f \semiprod{\id[X]} \big( g \circ (p \times \id[Y]) \big), \]
    as one checks immediately from the defining composite, so that every semi-direct product is one taken along an identity map.
    We nonetheless retain $E$ in the definition for two reasons.
    First, it is what makes \cref{semiprod-pb} natural to state: there the interest lies precisely in varying $X$ along a map $\alpha$ while keeping $E$ and $g$ fixed, and it is the factorization of $p$ through $E$ that makes the square a pullback.
    Second, $E$ corresponds to the ``external inputs'' of a network with external inputs in the sense of \cite{kadelka2023modularity}, so keeping it makes the comparison with that paper transparent.
    Where $E$ genuinely adds nothing, as in \cref{independent-inputs-implies-semiprod,good-decomposition-implies-semiprod}, we suppress it.
\end{remark}

In $\Set$, the semi-direct product defines a discrete dynamical system $(X \times Y, f \semiprod{p} g)$.
The projection map $\proj_X \from X \times Y \to X$ is equivariant with respect to $f \semiprod{p} g$ and $f$; that is, we have a commutative square
\[ \begin{tikzcd}
    X \times Y \ar[r, "\proj_X"] \ar[d, "f \semiprod{p} g"'] & X \ar[d, "f"] \\
    X \times Y \ar[r, "\proj_X"] & X
\end{tikzcd} \] 
thus $\proj_X \from (X \times Y, f \semiprod{p} g) \to (X, f)$ is a morphism of discrete dynamical systems.

Before proceeding, we note that this definition generalizes the semi-direct product of \cite{kadelka2023modularity}.
In fact, the use of abstraction and category-theoretic language makes our definition simpler, while recovering the notion in \cite{kadelka2023modularity} as an example.

The next proposition warrants an informal explanation.
Our aim is to show that the attractor functor turns a semi-direct product into a decomposition of attractors.
The only leverage we have over the attractor functor is that it preserves limits, so the entire strategy must be to express the semi-direct product as a limit of some kind.
That is exactly what \cref{semiprod-pb} does: it says that a semi-direct product $f \semiprod{p} g$ sitting over a base $(X,f)$ is a \emph{pullback} of a semi-direct product over any larger base $(X', f')$ that $X$ maps into compatibly.
Applied with $(X, f) = (\ZZ/k, \succ)$ a single attractor of the base and $\alpha = c$ the inclusion of that attractor, it says that the part of the system lying over one attractor of the base is cut out by a pullback square, and pullbacks are preserved by the attractor functor.
This is the whole proof of \cref{modularity-theorem} in outline; everything else is bookkeeping over the orbits.

\begin{proposition} \label{semiprod-pb}
    Let $\alpha \from (X, f) \to (X', f')$ be a morphism of dynamical systems and let $p, p'$ be morphisms fitting into a commutative triangle
    \[ \begin{tikzcd}
        X \ar[rr, "\alpha"] \ar[rd, "p"'] & {} & X' \ar[ld, "p'"] \\
        {} & E & {}
    \end{tikzcd} \]
    Then, for any function $g \from E \times Y \to Y$, the map
    \[ \alpha \times Y \from X \times Y \to X' \times Y \]
    ascends to a morphism $(X \times Y, f \semiprod{p} g) \to (X' \times Y, f' \semiprod{p'} g)$ in $\DDS$ such that the square
    \[ \begin{tikzcd}
        (X \times Y, f \semiprod{p} g) \ar[r, "{\alpha \times Y}"] \ar[d, "\proj_{X}"'] & (X' \times Y, f' \semiprod{p'} g) \ar[d, "\proj_{X'}"] \\
        (X, f) \ar[r, "\alpha"] & (X', f')
    \end{tikzcd} \]
    is a pullback diagram. 
\end{proposition}
\begin{proof}
    We verify that $\alpha \times Y$ is equivariant by computing
    \[ (\alpha \times Y) \big( (f \semiprod{p} g)(x, y) \big) = \big( \alpha(f(x)), g(p(x), y) \big) = \big( f'(\alpha(x)), g(p'(\alpha(x)), y) \big) = (f' \semiprod{p'} g) \big( (\alpha \times Y) (x, y) \big). \]
    The square
    \[ \begin{tikzcd}
        (X \times Y, f \semiprod{p} g) \ar[r, "{\alpha \times Y}"] \ar[d, "\proj_{X}"'] & (X' \times Y, f' \semiprod{p'} g) \ar[d, "\proj_{X'}"] \\
        (X, f) \ar[r, "\alpha"] & (X', f')
    \end{tikzcd} \]
    commutes by definition of $\alpha \times Y$.

    An explicit description of the pullback is given by
    \[ P := \{ (x, x', y) \in X \times X' \times Y \mid \alpha(x) = x' \} \]
    with the endofunction $f_P$ defined by 
    \[ f_P(x, x', y) = \big( f(x), f'(x'), g(p'(x'), y) \big). \]
    Define a map $\Phi \from X \times Y \to P$ by
    \[ \Phi(x, y) = (x, \alpha(x), y). \]
    The map $\Phi$ admits an inverse $\Phi^{-1} \from P \to X \times Y$ given by
    \[ \Phi^{-1}(x, x', y) = (x, y). \]
    The diagram
    \[ \begin{tikzcd}
        {} & X \times Y \ar[d, "\Phi"] \ar[ld, "\proj_X"'] \ar[rd, "\alpha \times Y"] & {} \\
        X & P \ar[r, "{\proj_{X' \times Y}}"'] \ar[l, "{\proj_{X}}"] & X' \times Y
    \end{tikzcd} \]
    commutes by construction, so it remains only to verify that $\Phi$ is equivariant.
    This is a straightforward calculation
    \begin{align*}
        \Phi \big( (f \semiprod{p} g)(x, y) \big) &= \big( f(x), \alpha(f(x)), g(p(x), y) \big) \\
        &= \big( f(x), f'(\alpha(x)), g(p'(\alpha(x)), y) \big) \\
        &= f_P(x, \alpha(x), y) \\
        &= f_P(\Phi(x, y)). \qedhere
    \end{align*}
\end{proof}
% \begin{corollary} \label{semiprod-pb-along}
%     Let $(X, f)$ be a discrete dynamical system and $\alpha \from (X', f') \to (X, f)$ be a map of dynamical systems.
%     Then, for any $g \from E \times Y \to Y$ and $p \from X \to E$, we have a pullback diagram
%     \[ \begin{tikzcd}[row sep = large]
%         (X' \times Y, f_A \semiprod{pi} g) \ar[r, "i \times Y"] \ar[d, "\proj_{A}"'] \ar[rd, phantom, "\pbtick" very near start, xshift=-3ex] & (X \times Y, f \semiprod{p} g) \ar[d, "\proj_X"] \\
%         (X', f') \ar[r, "i"] & (X, f)
%     \end{tikzcd} \tag*{\qedsymbol} \]
% \end{corollary}

Having established the requisite generalization of the semi-direct product, we are now ready to state and prove our decomposition theorem.
Recall from \cref{cycles-of-a-dynamial-system-correspondence} that, for a discrete dynamical system $(X, f)$, the $n$-cycles of the cycle set of attractors $AS(X, f)$ are in one-to-one correspondence with maps of dynamical systems $(\ZZ / n) \to (X, f)$.
Going forward, we will not distinguish between these.
For instance, if $c \in AS(X, f)_n$ is an $n$-cycle and $\alpha \from (X, f) \to (Y, g)$ is a map of dynamical systems then we may write $\alpha c$ for the composite map $\alpha \circ c \from (\ZZ / n, \succ) \to (Y, g)$, which is an $n$-cycle of $AS(Y, g)$.
In particular, a graph map between cycles $\varphi \from C_m \to C_n$ is an $m$-cycle in $AS(\ZZ / n, \succ)$, which we identify as a map of dynamical systems $\varphi \from (\ZZ / m, \succ) \to (\ZZ / n, \succ)$.

Before stating the theorem, we fix notation carefully, since several indices are in play at once.
Let $(X, f)$ be a dynamical system and, for $n \geq 1$, let $\omega \in \Orb(AS(X,f)_n)$ be an orbit of $n$-cycles.
Write $k_\omega := \lvert \omega \rvert$ for the size of this orbit.
By item (3) of \cref{attractor-digraph-injetive-unique-degeneracies}, $k_\omega$ is also the minimal length of any cycle in $\omega$, and by item (2) each $\tilde{c} \in \omega$ is a degeneracy $\tilde{c} = c \modmap{n}{k_\omega}$ of a unique non-degenerate cycle $c \in AS(X,f)_{k_\omega}$.
As $\tilde{c}$ ranges over $\omega$, these $c$ range over a single orbit of $AS(X,f)_{k_\omega}$, so we may choose one and call it $c_\omega$.
Under the identification of \cref{cycles-of-a-dynamial-system-correspondence}, $c_\omega$ is a map of dynamical systems $c_\omega \from (\ZZ / k_\omega, \succ) \to (X, f)$, and hence $p \, c_\omega \from \ZZ / k_\omega \to E$ is a legitimate gluing map.

\begin{theorem} \label{modularity-theorem}
    Let $(X \times Y, f \semiprod{p} g)$ be a semi-direct product and let $n \geq 1$.
    With notation as above, we have an isomorphism of $(\ZZ / n)$-sets
    \[ AS(X \times Y, f \semiprod{p} g)_n \cong \coprod\limits_{\omega \in \Orb (AS(X, f)_n)} AS\big(\ZZ / k_\omega \times Y, \ \succ \semiprod{p \, c_\omega} g\big)_n. \]
    Moreover, up to isomorphism, the dynamical system $(\ZZ / k_\omega \times Y, \succ \semiprod{p \, c_\omega} g)$ does not depend on the choice of $c_\omega$ within its orbit.
\end{theorem}

\begin{remark} \label{rmk:subscript-is-pc}
    We draw attention to the subscript on the right-hand side, which is $p \, c_\omega$ and not $p \, \tilde{c}$ for $\tilde{c} \in \omega$.
    This is deliberate and is forced by the domains: $\tilde{c}$ is an $n$-cycle, i.e.\ a map out of $(\ZZ/n, \succ)$, whereas the system appearing on the right has base $(\ZZ / k_\omega, \succ)$, so the gluing map must be precomposed with $c_\omega$, which is defined on $\ZZ / k_\omega$.
    Of course $k_\omega$ divides $n$ and $\tilde{c} = c_\omega \modmap{n}{k_\omega}$, so the two are closely related, but only $p \, c_\omega$ typechecks.
\end{remark}
% Before proving \cref{modularity-theorem}, we prove a short lemma about $G$-sets.
% \begin{lemma} \label{g-set-mono}
%     Let $G$ be a (discrete) group and $f \from A \to B$ be an injective map between $G$-sets.
%     If the action on $B$ is transitive and $A$ is nonempty then $f$ is an isomorphism
% \end{lemma}
% \begin{proof}
%     It suffices to show $f$ is surjective.
%     Fix an element $a \in A$.
%     Since $B$ is transitive, for any $b$, there exists $g \in G$ such that $g \cdot f(a) = b$.
%     Then, $f(g \cdot a) = g \cdot f(a) = b$, as desired.
% \end{proof}
\begin{proof}
    Every $(\ZZ / n)$-set is the coproduct of its orbits; the orbits partition the underlying set, and the action restricts to each block of the partition.
    % (This requires no appeal to the orbit--stabilizer theorem, which concerns the \emph{size} of an orbit rather than this decomposition.)
    The map
    \[ AS(\proj_{X})_n \from AS(X \times Y, f \semiprod{p} g)_n \to AS(X, f)_n \]
    is equivariant and hence sends orbits into orbits.
    Given an orbit $\omega$ of $AS(X, f)_n$, let $\proj_X^{-1}\omega$ denote the set of orbits of $AS(X \times Y, f \semiprod{p} g)_n$ which are sent into $\omega$.
    With this, we obtain an isomorphism
    \[ AS(X \times Y, f \semiprod{p} g)_n \cong \coprod\limits_{\omega \in \Orb(AS(X, f)_n)} \left( \coprod\limits_{[d] \in \proj_X^{-1}\omega} [d] \right) \tag{1} \]
    in $\ZnSet{n}$.
    We consider each orbit $[d]$ as an object in the slice category $\ZnSet{n} \slice AS(X, f)_n$ by restricting the map $AS(\proj_X)_n$.
    Since coproducts in the slice category are computed underlyingly, the isomorphism (1) ascends to an isomorphism in $\ZnSet{n} \slice AS(X, f)_n$.

    % That is,
    % \[ \proj_X^{-1}[c] := \{ d \in AS(X \times Y, f \semiprod{p} g)_n \mid AS(\proj_X)_n(d) \in [c] \}. \]
    % Then, we may re-arrange this coproduct as
    % \[ AS(X \times Y, f \semiprod{p} g)_n \cong \coprod\limits_{[c] \in \Orb(X, f)} \coprod\limits_{} \]

    % For an orbit $[c]$ of $AS(X, f)_n$, let $k$ denote the cardinality of $[c]$ and suppose $c$ is a distinguished element.
    % By definition, $c$ is a graph map $C_n \to S(X, f)$, and by \cref{attractor-digraph-injetive-unique-degeneracies}, there is a unique $\tilde{c} \from C_{k} \to S(X, f)$ such that $c = \tilde{c} \modmap{n}{k}$.
    % The cycle graph $C_k$ is the state space of the discrete dynamical system $(\ZZ / k, \succ)$.
    % By \cref{state-space-functor-full-faithful}, there is a unique map of dynamical systems $\alpha \from (\ZZ / k, \succ) \to (X, f)$ such that $\tilde{c} = S(\alpha)$.

    Fix an orbit $\omega$ and abbreviate $k := k_\omega$ and $c := c_\omega$.
    By \cref{semiprod-pb} applied to $\alpha := c \from (\ZZ/k, \succ) \to (X, f)$, the square
    \[ \begin{tikzcd}
        (\ZZ / k \times Y, \succ \semiprod{pc} g) \ar[r, "c \times Y"] \ar[d, "\proj_{\ZZ / k}"'] \ar[rd, phantom, "\pbtick" very near start] & (X \times Y, f \semiprod{p} g) \ar[d, "\proj_{X}"] \\
        (\ZZ / k, \succ) \ar[r, "c"] & (X, f)
    \end{tikzcd} \]
    is a pullback.
    Each of the functors in the composite
    \[ \DDS \xrightarrow{S} \DiGraph \xrightarrow{A} \cySet \xrightarrow{\operatorname{ev}_n} \ZnSet{n} \]
    preserves limits by \cref{cor:S_preserves_lim_colim,attractor-functors-right-adj-preserves-lims,ev-preserves-limits} respectively. 
    We therefore deduce that the diagram
    \[ \begin{tikzcd}[sep = large]
        AS(\ZZ / k \times Y, \succ \semiprod{pc} g)_n \ar[r, "{AS(c \times Y)_n}"] \ar[d, "{\proj_{\ZZ / k}}"'] \ar[rd, phantom, "\pbtick" very near start] & AS(X \times Y, f \semiprod{p} g)_n \ar[d, "{AS(\proj_{X})_n}"] \\
        AS(\ZZ / k, \succ)_n \ar[r, "AS(c)_n"] & AS(X, f)_n
    \end{tikzcd} \]
    is again a pullback.

    % We now pull the decomposition of Equation (1) back along $AS(c)_n$.
    The category $\ZnSet{n}$ is infinitary extensive --- equivalently, coproducts in it are stable under pullback, as they are computed underlyingly in $\Set$ --- so the isomorphism of Equation (1) yields an isomorphism between pullbacks:
    \[ AS(\ZZ / k \times Y, \succ \semiprod{pc} g)_n \cong \coprod\limits_{\omega' \in \Orb(AS(X, f)_n)} \left( \coprod\limits_{[d] \in \proj_X^{-1}\omega'} AS(c)_n^*[d] \right), \tag{2} \]
    where $AS(c)^*_n[d]$ denotes the pullback:
    \[ \begin{tikzcd}[sep = large]
        AS(c)^*_n[d] \ar[r] \ar[d] \ar[rd, phantom, "\pbtick" very near start] & {[d]} \ar[d, "AS(\proj_X)_n"] \\
        AS(\ZZ / k, \succ)_n \ar[r, "AS(c)_n"] & AS(X, f)_n
    \end{tikzcd} \]
    The pullback $AS(c)^*_n[d]$ is empty whenever $\omega' \neq \omega$, since the bottom and right maps then have disjoint images.
    The bottom map $AS(c)_n$ is injective with image exactly $\omega$, which has $k$ elements; since $AS(\ZZ/k, \succ)_n$ also has $k$ elements, $AS(c)_n$ is an isomorphism onto $\omega$.
    Thus, the top map is an isomorphism, so we may re-write the isomorphism in (2) as
    \[ AS(\ZZ / k \times Y, \succ \semiprod{pc} g)_n \cong \coprod\limits_{[d] \in \proj_X^{-1}\omega} [d] . \]
    Combining with the isomorphism in (1) yields
    \[ AS(X \times Y, f \semiprod{p} g)_n \cong \coprod\limits_{\omega \in \Orb (AS(X, f)_n)} AS\big(\ZZ / k_\omega \times Y, \succ \semiprod{p \, c_\omega} g\big)_n, \]
    as desired.

    Finally, we check independence of the choice of $c_\omega$.
    Any two choices differ by an element of $\ZZ / k$: if $c$ and $c'$ are non-degenerate $k$-cycles in the same orbit then $c' = c \autmap{k}{i}$ for some $i$.
    By \cref{semiprod-pb} applied to the automorphism $\autmap{k}{i} \from (\ZZ/k, \succ) \to (\ZZ/k, \succ)$, which satisfies $p c' = (p c) \autmap{k}{i}$, the map $\autmap{k}{i} \times Y$ is an isomorphism of dynamical systems
    \[ (\ZZ / k \times Y, \succ \semiprod{p c'} g) \cong (\ZZ / k \times Y, \succ \semiprod{p c} g), \]
    as desired.
\end{proof}

We illustrate \cref{modularity-theorem} by working it out completely on the running example.

\begin{example}[Running example, continued] \label{ex:running-decomposition}
    Let $(X, f)$ be the system of \cref{ex:running}, so $X = \{a, b, c\}$ with $f(a) = a$, $f(b) = c$, $f(c) = b$.
    Take $E := \{0, 1\}$ with
    \[ p(a) = p(c) = 0, \qquad p(b) = 1, \]
    take $Y := \{0,1\}$, and let $g \from E \times Y \to Y$ be $g(e, y) := e \oplus y$, where $\oplus$ denotes addition modulo $2$.
    The semi-direct product is the system on $X \times Y$, a set with six elements, given by
    \[ (f \semiprod{p} g)(x, y) = \big( f(x), \, p(x) \oplus y \big). \]
    Explicitly, its state space graph is
    \[ \begin{tikzpicture}[scale=1.0]
        \node[vertex, label={[label distance=-1ex]below:$(a,0)$}] (a0) at (0,-0.35) {};
        \node[vertex, label={[label distance=-1ex]below:$(a,1)$}] (a1) at (1.9,-0.35) {};

        \node[vertex, label={[label distance=-1ex]above:$(b,0)$}] (b0) at (4.4,0.45) {};
        \node[vertex, label={[label distance=-1ex]above:$(c,1)$}] (c1) at (6.2,0.45) {};
        \node[vertex, label={[label distance=-1ex]below:$(b,1)$}] (b1) at (6.2,-1.15) {};
        \node[vertex, label={[label distance=-1ex]below:$(c,0)$}] (c0) at (4.4,-1.15) {};

        \draw[thick, ->] (a0) to [out=45, in=125, looseness=6] (a0);
        \draw[thick, ->] (a1) to [out=45, in=125, looseness=6] (a1);

        \draw[thick, ->] (b0) to (c1);
        \draw[thick, ->] (c1) to (b1);
        \draw[thick, ->] (b1) to (c0);
        \draw[thick, ->] (c0) to (b0);
    \end{tikzpicture} \]
    The two states over the fixed point $a$ are themselves fixed, while the four states over the $2$-cycle $\{b,c\}$ form a single $4$-cycle.
    Let us confirm this against the theorem, taking $n = 2$.
    From \cref{ex:running-attractors}, $AS(X,f)_2 = \{ (aa), (bc), (cb) \}$ has two orbits:
    \begin{itemize}
        \item $\omega_1 = \{ (aa) \}$, of size $k_{\omega_1} = 1$, with non-degenerate representative $c_{\omega_1} = (a) \from (\ZZ/1, \succ) \to (X,f)$.
        Here $p \, c_{\omega_1} \from \ZZ/1 \to E$ is constant at $p(a) = 0$, so $\succ \semiprod{p c_{\omega_1}} g$ is the system on $\ZZ/1 \times Y \cong \{0,1\}$ given by $y \mapsto 0 \oplus y = y$, the identity.
        Its state space graph is two looped vertices, so it has exactly two $2$-cycles.
        \item $\omega_2 = \{ (bc), (cb) \}$, of size $k_{\omega_2} = 2$, with non-degenerate representative $c_{\omega_2} = (bc) \from (\ZZ/2, \succ) \to (X,f)$.
        Here $p \, c_{\omega_2} \from \ZZ / 2 \to E$ sends $0 \mapsto p(b) = 1$ and $1 \mapsto p(c) = 0$, so $\succ \semiprod{p c_{\omega_2}} g$ is the system on $\ZZ/2 \times Y$ given by
        \[ (i, y) \mapsto \big( i + 1, \ (p c_{\omega_2})(i) \oplus y \big). \]
        Starting at $(0,0)$ we get $(0,0) \mapsto (1,1) \mapsto (0,1) \mapsto (1,0) \mapsto (0,0)$: a single $4$-cycle through all four states.
        A $4$-cycle contributes no $2$-cycles at all, so this system has no $2$-cycles.
    \end{itemize}
    The theorem therefore predicts $\lvert AS(X \times Y, f \semiprod{p} g)_2 \rvert = 2 + 0 = 2$.
    Reading off the state space graph above, the closed walks of length $2$ are $\big((a,0)(a,0)\big)$ and $\big((a,1)(a,1)\big)$ --- exactly two, as predicted, since the $4$-cycle over $\omega_2$ contributes nothing in degree $2$.

    % This example repays a second look, because it shows the theorem doing real work.
    Note that the base system $(X,f)$ has an attractor of length $2$, and the fibre system $g$ is an involution, yet the composite has a $4$-cycle: the ``twisting'' by $p$ has doubled the period.
    This is precisely the phenomenon that \cref{old-decomposition-theorem} was designed to capture, and which the naive guess $\mathcal{A}(f \semiprod{p} g) = \mathcal{A}(f) \times \mathcal{A}(g)$ would get wrong.
\end{example}

\section{Wiring diagrams and effective decomposition} \label{sec:wiring}

\Cref{modularity-theorem} tells us what a semi-direct product decomposition buys us, but it presupposes that we have one.
This section addresses the complementary question: given a system, how do we \emph{find} such a decomposition, or determine that none exists?

The answer, and the point of this section, is that the question can be settled by inspecting a small combinatorial invariant of the system --- its wiring diagram --- rather than its state space.
This is an instance of a general phenomenon: a well-chosen diagram can make structure visible that is buried in an equivalent but less perspicuous representation \cite{larkin:diagram}.
Moreover, in the biological setting the wiring diagram is often the representation actually available to the modeller \cite{shilts:wiring_diagram}.
For a system $f \from A^n \to A^n$, the state space has $\lvert A \rvert^n$ elements, whereas the wiring diagram has $n$ vertices.
Our main result, \cref{semiprod-iff-wiring-map}, says that $f$ admits a semi-direct product decomposition if and only if $W(f)$ admits a graph map to a fixed two-vertex digraph $\DEdge$, a condition that can be checked directly.
Combining this with \cref{modularity-theorem} gives \cref{wiring-map-implies-attractors-decompose}, which is the practical payoff of the paper: a check on $n$ vertices certifies that the attractors decompose.

We recall from \cref{rmk:wiring-terminology} that our use of the term ``wiring diagram'' is the one standard in the Boolean network literature, and is unrelated to the operadic notion of the same name.

\subsection{The setting: regular categories and supports} \label{subsec:regular}

The results of this section hold well beyond the category of sets, and we set up the appropriate level of generality now.
A reader interested only in the case $\cat{C} = \Set$ and $A$ a non-empty set loses nothing by keeping that case in mind: there every hypothesis below is automatic, and $A^0$ is the one-point set.

Throughout this section, $\cat{C}$ denotes a \emph{regular} category: a category with finite limits in which every morphism factors as a regular epimorphism followed by a monomorphism, and in which regular epimorphisms are stable under pullback.
We refer to \cite[\S A1.3]{johnstone:elephant} for the theory.
Regularity is a mild assumption: every topos (in particular $\Set$, $\DDS$, and every presheaf category), every abelian category, and every variety of algebras is regular.

% There is one convention which does real work and which we state explicitly.
% We thank an anonymous referee for introducing us to the following convention
% The following definition establishes our convention

\begin{definition} \label{def:support-convention}
    Let $A$ be an object of a regular category $\cat{C}$.
    The \emph{support} of $A$, written $\supp(A)$, is the subobject of the terminal object $1$ obtained by factoring the unique morphism $A \to 1$ as a regular epimorphism followed by a monomorphism:
    \[ A \twoheadrightarrow \supp(A) \rightarrowtail 1. \]
\end{definition}
Throughout, we adopt the notational convention
\[ A^0 := \supp(A), \]
in contrast with the more standard convention $A^0 := 1$.
The projection $A^1 \to A^0$ is understood to be the regular epimorphism $A \twoheadrightarrow \supp(A)$ in this factorization.
An object $A$ is \emph{well-supported} if $A \to 1$ is a regular epimorphism, equivalently if $\supp(A) = 1$, in which case the convention agrees with the usual one.

\begin{remark} \label{rmk:why-support}
    In $\Set$, $\supp(A)$ is the empty set if $A = \varnothing$ and the one-point set otherwise, so if $A$ is a non-empty set then the convention changes nothing.
    % Similarly in $\DDS$, an object is well-supported if and only if its underlying set is nonempty.
    Its purpose is to make \cref{projections-regular-epi} below true without any hypothesis on $A$ (such as non-emptiness in the case of $A \in \Set$).
    % The alternative --- requiring a global element $1 \to A$, so that the projections admit sections --- is a genuinely stronger condition and one that fails in cases of interest.
\end{remark}

% The following is the one fact about this setup that we use, and everything in this section rests on it.
Our motivation for working with regular categories is the following result.

\begin{proposition} \label{projections-regular-epi}
    Let $A$ be an object of a regular category $\cat{C}$, and adopt the convention $A^0 = \supp(A)$.
    Then for every $n \geq 1$ and every $i \in \{1, \dots, n\}$, the projection
    \[ \proj_{\neq i} \from A^n \to A^{n-1} \]
    is a regular epimorphism.
    More generally, $\proj_{\neq S} \from A^n \to A^{n - \lvert S \rvert}$ is a regular epimorphism for every $S \subseteq \{1, \dots, n\}$.
\end{proposition}
\begin{proof}
    By symmetry we may assume $i = n$, so that $\proj_{\neq i}$ is the projection $A^{n-1} \times A \to A^{n-1}$.

    For $n = 1$ this is the morphism $A \to \supp(A)$, which is a regular epimorphism by \cref{def:support-convention}.

    For $n \geq 2$, note first that $\supp(A)$ is a subterminal object, and that both $A^{n-1}$ and $A$ admit (necessarily unique) morphisms to it.
    For any $X, Y$ equipped with morphisms to a subterminal $U$, the canonical map $X \times Y \to X \times_U Y$ is an isomorphism, since the commutativity condition defining the pullback is automatic.
    Hence the square
    \[ \begin{tikzcd}
        A^{n-1} \times A \ar[r] \ar[d, "\proj_{\neq n}"'] \ar[rd, phantom, "\pbtick" very near start] & A \ar[d] \\
        A^{n-1} \ar[r] & \supp(A)
    \end{tikzcd} \]
    is a pullback.
    Its right-hand vertical map $A \to \supp(A)$ is a regular epimorphism, and regular epimorphisms in a regular category are stable under pullback, so $\proj_{\neq n}$ is a regular epimorphism.

    The last statement follows since $\proj_{\neq S}$ is a composite of morphisms of the above form, and regular epimorphisms compose.
\end{proof}

Since regular epimorphisms are in particular epimorphisms, we obtain the following, which is what makes \cref{def:independent-input} well behaved.

\begin{corollary} \label{projection-epi}
    For every $S \subseteq \{1, \dots, n\}$, the projection $\proj_{\neq S}$ is an epimorphism.
    Consequently, a morphism factors through $\proj_{\neq S}$ in at most one way. \qed
\end{corollary}

\subsection{Independent inputs} \label{subsec:independent}

We can now define what it means for a morphism to be independent of an input.

\begin{definition} \label{def:independent-input}
    Let $f \colon A^n \to B$ be a morphism in $\cat{C}$.
    The morphism $f$ is \emph{independent} of input $i$ if it factors through the projection $\proj_{\neq i}$, i.e.\ if there exists a (necessarily unique, by \cref{projection-epi}) morphism $\overline{f}$ making the triangle
    \[
    \begin{tikzcd}[column sep=large, row sep=large]
        A^n
        \ar[r, "f"]
        \ar[d, "\proj_{\neq i}", swap]
    &
        B
    \\
        A^{n-1}
        \ar[ru, dotted, "\overline{f}", swap]
    \end{tikzcd}
    \]
    commute.
    More generally, for $S \subseteq \{1, \dots, n\}$, we say $f$ is \emph{independent of the inputs in $S$} if it factors through $\proj_{\neq S}$.
\end{definition}

As one might expect, we say a morphism depends on a particular input if it is not independent of this input.
\begin{definition} \label{def:depends-on-input}
    A morphism $f \colon A^n \to B$ in $\cat{C}$ \emph{depends} on input $i$ if it is not independent of input $i$.
\end{definition}

\begin{example} \label{ex:independent-in-Set}
    In $\Set$ with $A$ non-empty, a function $f \from A^n \to B$ is independent of input $i$ precisely when
    \[ f(a_1, \dots, a_{i-1}, a, a_{i+1}, \dots, a_n) = f(a_1, \dots, a_{i-1}, a', a_{i+1}, \dots, a_n) \]
    for all $a, a' \in A$ and all choices of the remaining coordinates; that is, the value of $f$ does not vary as the $i$-th argument varies.
    This is the familiar elementwise notion, and \cref{def:independent-input} is its faithful generalization.
\end{example}

The key technical result of this subsection is that independence of several inputs, taken one at a time, implies independence of all of them at once.
This is what allows us to pass from the edgewise data recorded by the wiring diagram to an actual factorization.
For this, we use the following pushout property of squares of projections in a regular category.

\begin{lemma} \label{projection-square-pushout}
    Let $n \geq 1$ and let $i \neq j$ be two indices in $\{1, \dots, n+1\}$.
    Then the square of projections
    \[ \begin{tikzcd}[sep = large]
        A^{n+1} \ar[r, "\proj_{\neq j}"] \ar[d, "\proj_{\neq i}"'] & A^{n} \ar[d, "\proj_{\neq i}"] \\
        A^{n} \ar[r, "\proj_{\neq j}"] & A^{n-1}
    \end{tikzcd} \]
    is both a pullback and a pushout in $\cat{C}$.
\end{lemma}
\begin{proof}
    Reordering coordinates, we may identify $A^{n+1} \cong R \times A \times A$, where $R := A^{n-1}$ carries the coordinates other than $i$ and $j$, and identify both copies of $A^n$ with $R \times A$ in such a way that the two upper maps become the projections
    \[ \pi_{12}, \pi_{13} \from R \times A \times A \longrightarrow R \times A \]
    discarding the third and second factor respectively, and both lower maps become the projection $\pi \from R \times A \to R$.

    That the square is a pullback is standard: $R \times A \times A$ is exactly the pullback $(R \times A) \times_R (R \times A)$, i.e.\ the kernel pair of $\pi$.

    For the pushout property, fix $u, v \from R \times A \to Z$ with $u \pi_{12} = v \pi_{13}$.
    Let $\delta$ denote the morphism
    \[ \id[R] \times \Delta_A \from R \times A \to R \times A \times A\] 
    induced by the diagonal of $A$, so that $\pi_{12}\delta = \pi_{13} \delta = \id$.
    Whiskering the equality $u \pi_{12} = v \pi_{13}$ with $\delta$ gives
    \[ u = u \pi_{12} \delta = v \pi_{13} \delta = v, \]
    so the two maps coincide; call the common map $u$.
    Then $u \pi_{12} = u \pi_{13}$, i.e.\ $u$ coequalizes the kernel pair of $\pi$.
    By \cref{projections-regular-epi}, $\pi$ is a regular epimorphism, hence is the coequalizer of its own kernel pair, and therefore $u$ factors uniquely as $u = w \pi$ for a unique $w \from R \to Z$.
    This $w$ is the required map out of the pushout corner, and it is unique because $\pi$ is an epimorphism.
\end{proof}

\begin{lemma} \label{multiple-independent-lift}
    Let $f \from A^n \to B$ be a morphism and let $S \subseteq \{1, \dots, n\}$.
    If $f$ is independent of input $i$ for every $i \in S$, then $f$ is independent of the inputs in $S$; that is, $f$ factors uniquely through $\proj_{\neq S} \from A^n \to A^{n - \lvert S \rvert}$.
\end{lemma}
\begin{proof}
    We induct on $\lvert S \rvert$.
    The cases $\lvert S \rvert = 0$ and $\lvert S \rvert = 1$ are trivial and hold by definition respectively, uniqueness in each case coming from \cref{projection-epi}.

    Suppose $\lvert S \rvert = k \geq 2$ and write $S = S' \sqcup \{\ell\}$ with $\lvert S' \rvert = k - 1$.
    By the inductive hypothesis, $f$ factors as $f = g \circ \proj_{\neq S'}$ for a unique $g \from A^{n-k+1} \to B$.
    Let $\bar{\ell}$ denote the coordinate of $A^{n-k+1}$ corresponding to $\ell$.
    It suffices to show that $g$ is independent of input $\bar{\ell}$, for then the base case gives $g = w \circ \proj_{\neq \bar{\ell}}$, whence
    \[ f = w \circ \proj_{\neq \bar{\ell}} \circ \proj_{\neq S'} = w \circ \proj_{\neq S} \]
    as required, with uniqueness again from \cref{projection-epi}.

    To that end, reorder coordinates to identify $A^n \cong R \times A \times A^{S'}$, where $R := A^{n-k}$ carries the coordinates outside $S$, the middle factor is the $\ell$-th, and $A^{S'} := A^{k-1}$ carries the coordinates in $S'$.
    Under this identification, $\proj_{\neq S'}$ becomes the projection onto $R \times A$ and $\proj_{\neq \ell}$ becomes the projection onto $R \times A^{S'}$.
    Since $f$ is independent of $\ell$, we may write $f = h \circ \proj_{\neq \ell}$ for some $h \from R \times A^{S'} \to B$.

    Consider $W := R \times A \times A \times A^{S'}$ together with the two projections
    \[ \alpha, \beta \from W \longrightarrow R \times A \times A^{S'} \cong A^n \]
    which discard the third and the second factor respectively.
    Writing $\pi \from W \to R \times A \times A$ for the projection discarding $A^{S'}$, we have
    \[ \proj_{\neq S'} \circ \alpha = \pi_{12} \circ \pi, \qquad \proj_{\neq S'} \circ \beta = \pi_{13} \circ \pi, \qquad \proj_{\neq \ell} \circ \alpha = \proj_{\neq \ell} \circ \beta, \]
    the last because both sides discard the $\ell$-th coordinate and retain the same remaining ones.
    Therefore
    \[ g \circ \pi_{12} \circ \pi = f \circ \alpha = h \circ \proj_{\neq \ell} \circ \alpha = h \circ \proj_{\neq \ell} \circ \beta = f \circ \beta = g \circ \pi_{13} \circ \pi . \]
    Since $S'$ is non-empty, $\pi$ is a regular epimorphism by \cref{projections-regular-epi}, in particular an epimorphism, so it may be cancelled:
    \[ g \circ \pi_{12} = g \circ \pi_{13} \from R \times A \times A \longrightarrow B. \]
    But $(\pi_{12}, \pi_{13})$ is the kernel pair of $\proj_{\neq \bar{\ell}} \from R \times A \to R$, which is a regular epimorphism by \cref{projections-regular-epi} and hence the coequalizer of that kernel pair.
    Therefore $g$ factors through $\proj_{\neq \bar{\ell}}$, i.e.\ $g$ is independent of input $\bar{\ell}$, as required.
\end{proof}

% \begin{remark} \label{rmk:admissibility-superseded}
%     Readers familiar with the approach via sections may find it helpful to know how the two compare.
%     One can instead require that each projection $\proj_{\neq i}$ admit a \emph{section}, which for $n = 1$ amounts to the existence of a global element $1 \to A$; independence of $f$ from input $i$ is then equivalent to the equation $f = f \circ s \circ \proj_{\neq i}$ for any such section $s$.
%     That route is available, but the hypothesis is strictly stronger than what is needed --- as noted in \cref{rmk:why-support}, it fails for systems without fixed points --- and the combinatorics of composing sections makes the proof of \cref{multiple-independent-lift} considerably longer.
%     The formulation above applies to \emph{every} object of every regular category.
% \end{remark}

Given a morphism $f \colon A \to B^n$, we write $f_i$ for the composite $\proj_i \circ f$.
Note that $f = (f_1, f_2, \ldots, f_n)$.

\subsection{Semi-direct products and independent inputs} \label{subsec:semiprod-independent}

We are interested in characterizing semi-direct product morphisms $g \semiprod{p} h \from A^{m+n} \to A^{m+n}$ in terms of the independent inputs of certain coordinate functions.
The following result is the forward implication for this equivalence: if a morphism $A^{m+n} \to A^{m+n}$ is a semi-direct product then the first $m$ outputs are independent of the last $n$ inputs.
\begin{proposition} \label{semiprod-independent-inputs}
    Let $m, n \geq 1$.
    Suppose $f \from A^{m+n} \to A^{m+n}$ is a semi-direct product $f = g \semiprod{p} h$ of maps $g \from A^m \to A^m$ and $h \from E \times A^n \to A^n$ along $p \from A^m \to E$.
    Then for every $i \in \{ 1, \dots, m \}$, the coordinate morphism $f_i \from A^{m+n} \to A$ is independent of the inputs $\{ m+1, \dots, m+n \}$.
\end{proposition}
\begin{proof}
    Fix $i \in \{ 1, \dots, m \}$.
    By definition of the semi-direct product, we have a commutative diagram
    \[ \begin{tikzcd}[sep = large]
        A^{m} \times A^n \ar[r, "g \semiprod{p} h"] \ar[d, "{(\proj_{A^m}, \, p \circ \proj_{A^m}, \, \proj_{A^n})}"'] & A^m \times A^n \ar[d, "\proj_i"] \\
        A^m \times E \times A^n \ar[r, "\proj_{i} \circ \, g \circ \proj_{A^m}"] & A
    \end{tikzcd} \]
    Following the left-bottom composite, we see that $f_i$ equals $g_i \circ \proj_{A^m}$, which by definition factors through the projection $\proj_{A^m} = \proj_{\neq m+1, \dots, m+n}$.
    Hence $f_i$ is independent of the inputs $\{m+1, \dots, m+n\}$.
\end{proof}

The reverse implication says that this property characterizes semi-direct products.
That is, if the first $m$ outputs are independent of the last $n$ inputs then the morphism can be written as a semi-direct product.
\begin{theorem} \label{independent-inputs-implies-semiprod}
    Let $f \from A^{m+n} \to A^{m+n}$ be a morphism, and suppose that for every $i \in \{ 1, \dots, m \}$ and every $j \in \{ m+1, \dots, m+n \}$, the coordinate morphism $f_i \from A^{m+n} \to A$ is independent of input $j$.
    Then there exist morphisms
    \[ g \from A^m \to A^m, \qquad h \from A^m \times A^n \to A^n \]
    such that $f = g \semiprod{\id[A^m]} h$.
\end{theorem}
\begin{proof}
    Fix $i \in \{1, \dots, m\}$.
    By hypothesis $f_i$ is independent of each of the inputs $m+1, \dots, m+n$ individually, so by \cref{multiple-independent-lift} it is independent of all of them at once; that is, it factors as $f_i = g_i \circ \proj_{\leq m}$ for a unique $g_i \from A^m \to A$.
    \[ \begin{tikzcd}
        A^{m+n} \ar[r, "f_i"] \ar[d, "\proj_{\leq m}"'] & A \\
        A^{m} \ar[ur, dotted, "g_i"']
    \end{tikzcd} \]
    Assemble these into $g := (g_1, \dots, g_m) \from A^m \to A^m$.

    Now simply set $h := (f_{m+1}, \dots, f_{m+n}) \from A^{m+n} \to A^n$, regarded as a morphism $A^m \times A^n \to A^n$ under the identification $A^{m+n} = A^m \times A^n$; thus $h_j = f_{m+j}$ for $j \in \{1, \dots, n\}$.
    We claim $f = g \semiprod{\id[A^m]} h$.
    Indeed, unfolding \cref{semidirect-product} with $E = A^m$ and gluing map $\id[A^m]$, the right-hand side sends $(x, y) \in A^m \times A^n$ to $\big( g(x), h(x, y) \big)$, so its coordinate morphisms are
    \[ (g \semiprod{\id[A^m]} h)_k = \begin{cases}
        g_k \circ \proj_{\leq m} & \text{if } k \leq m, \\
        h_{k-m} & \text{if } k \geq m+1.
    \end{cases} \]
    For $k \leq m$ we have $g_k \circ \proj_{\leq m} = f_k$ by construction, and for $k \geq m+1$ we have $h_{k-m} = f_{m + (k-m)} = f_k$ by definition of $h$.
    Since the coordinate morphisms agree, $f = g \semiprod{\id[A^m]} h$.
\end{proof}

\begin{remark} \label{rmk:no-projection-needed}
    An earlier formulation of \cref{independent-inputs-implies-semiprod} produced a semi-direct product along a projection $A^\iota \from A^m \to A^{\ell}$,
    where $\iota \from \{ 1, \dots, \ell \} \ito \{ 1, \dots, m \}$ is an injection which omits any and all coordinates that the last $n$ outputs are independent of.
    By contrast, the current formulation instead takes the gluing map to be the identity, absorbing the projection into $h$.
    % By \cref{rmk:E-redundant}, one may instead take the gluing map to be the identity, absorbing the projection into $h$.
    If one wants the smallest possible $E$, it is better to restrict $h$ along $A^i$; this is a matter of economy of presentation rather than of expressive power, since the class of morphisms expressible as a semi-direct product is unchanged.
\end{remark}

\subsection{The wiring diagram} \label{subsec:wiring-diagram}

To an endomorphism $f \from A^n \to A^n$, one can associate a digraph known as the \emph{wiring diagram} of $f$, which records the dependent inputs of each coordinate morphism of $f$.
\begin{definition} \label{wiring-diagram}
    The \emph{wiring diagram} $W(f)$ of a morphism $f \colon A^n \to A^n$ is the directed graph whose vertices are the numbers $1$, $2$, \ldots, $n$, with an edge $i \to j$ if and only if $f_j$ depends on input $i$.
\end{definition}
\begin{example} \label{ex:wiring-diagram-basic}
    Consider the case where $A = \Ftwo$ is the field with 2 elements and $f \from \Ftwo^2 \to \Ftwo^2$ is defined by
    \[ f(x_1, x_2) := (x_1 + x_2, x_1). \]
    The coordinate function $f_1$ depends on both inputs 1 and 2, whereas the coordinate function $f_2$ depends only on input 1.
    Thus, the wiring diagram of $f$ is:
    
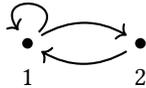
\[ \begin{tikzpicture}
        \node[vertex, label={south:1}] (A) {};
        \path (A) -- +(1.5, 0) node[vertex, label={south:2}] (B) {};

        \draw[thick, ->] (A) to [out=45, in=125, looseness=6] (A);
        \draw[thick, ->] (B) to [bend left] (A);
        \draw[thick, ->] (A) to [bend left] (B);
    \end{tikzpicture} \]
\end{example}

\begin{example}[A running example for this section] \label{ex:wiring-running}
    We fix a Boolean network to which we will return at each stage below.
    Let $A = \Ftwo$ and let $f \from \Ftwo^3 \to \Ftwo^3$ be given by
    \[ f(x_1, x_2, x_3) := (x_2, \ x_1, \ x_1 \wedge x_3). \]
    Reading off the dependencies: $f_1 = x_2$ depends on input $2$ only; $f_2 = x_1$ depends on input $1$ only; and $f_3 = x_1 \wedge x_3$ depends on inputs $1$ and $3$.
    Hence the wiring diagram is, drawing the vertices in the order $2$, $1$, $3$ so that no two edges overlap,
    \[ \begin{tikzpicture}
        \node[vertex, label={south:2}] (B) at (0,0) {};
        \node[vertex, label={south:1}] (A) at (1.8,0) {};
        \node[vertex, label={south:3}] (C) at (3.6,0) {};

        \draw[thick, ->] (B) to [bend left] (A);
        \draw[thick, ->] (A) to [bend left] (B);
        \draw[thick, ->] (A) to (C);
        \draw[thick, ->] (C) to [out=45, in=125, looseness=6] (C);
    \end{tikzpicture} \]
    Note that there is no edge from $3$ to $1$ or from $3$ to $2$: the third coordinate influences nothing but itself.
\end{example}

Let $\DEdge$ denote the digraph consisting of a single directed edge $0 \to 1$ and two self-loops.
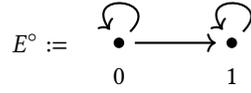
\[ \begin{tikzpicture}
    \node[vertex, label={south:0}] (A) at (0, 0) {};
    \path (0, 0) -- +(1.5, 0)  node[vertex, label={south:1}] (B) {};

    \path (0, 0) -- (-1.05, 0) node {$\DEdge := $};

    \draw[thick, ->] (A) to (B);
    \draw[thick, ->] (A) to [out=45, in=125, looseness=6] (A);
    \draw[thick, ->] (B) to [out=45, in=125, looseness=6] (B);
\end{tikzpicture} \]

We now come to the main result of this section.
Our final characterization of semi-direct product morphisms is in terms of the wiring diagram: a morphism $f \from A^n \to A^n$ can be written as a semi-direct product if and only if its wiring diagram $W(f)$ admits a graph map to $\DEdge$.
The forward implication is \cref{semiprod-implies-good-decomposition}, and the reverse implication is \cref{good-decomposition-implies-semiprod}.
Note that while the forward implication holds as stated, the reverse implication holds only up to a permutation of the set $\{ 1, \dots, n \}$, since a semi-direct product decomposition as we have defined it requires the two blocks of coordinates to be consecutive.

\begin{proposition} \label{semiprod-implies-good-decomposition}
    Let $m, n \geq 1$.
    Suppose $f \from A^{m+n} \to A^{m+n}$ is a semi-direct product $f = g \semiprod{p} h$ of maps $g \from A^m \to A^m$ and $h \from E \times A^n \to A^n$ along $p \from A^m \to E$.
    Then, the wiring diagram $W(f)$ admits a graph map $W(f) \to \DEdge$.
\end{proposition}
\begin{proof}
    Define a graph map $\varphi \from W(f) \to \DEdge$ by
    \[ \varphi(i) := \begin{cases}
        0 & i \leq m \\
        1 & i \geq m+1.
    \end{cases} \]
    To prove that $\varphi$ is a graph map, we must show that there are no edges in $W(f)$ from a vertex in $\{ m+1, \dots, m+n \}$ to a vertex in $\{ 1, \dots, m \}$.
    Unfolding the definition of $W(f)$, this is exactly \cref{semiprod-independent-inputs}.
\end{proof}

If $\sigma \in \Sigma_k$ is a permutation of the set $\{ 1, \dots, k \}$ and $f \from A^k \to A^k$ is a morphism in a category with products then we write $\sigma \cdot f \from A^k \to A^k$ for the morphism obtained by conjugation with $A^{\sigma}$, i.e.\ the morphism $A^{\sigma} \circ f \circ A^{\sigma^{-1}}$.
\[ \begin{tikzcd}
    A^k \ar[r, dotted, "\sigma \cdot f"] \ar[d, "A^{\sigma^{-1}}"'] & A^k \\
    A^k \ar[r, "f"] & A^k \ar[u, "A^{\sigma}"']
\end{tikzcd} \]

The following lemma says that conjugating by a permutation simply relabels the wiring diagram.
It is intuitively evident: permuting the names of the coordinates cannot change which coordinates influence which; we record it only to make the bookkeeping in \cref{good-decomposition-implies-semiprod} precise.

\begin{lemma} \label{permutation-relabels-wiring}
    Let $f \from A^k \to A^k$ be a morphism and $\sigma \in \Sigma_k$ a permutation.
    Then $(\sigma \cdot f)_i$ is independent of input $j$ if and only if $f_{\sigma(i)}$ is independent of input $\sigma(j)$.
    Consequently $\sigma$ defines an isomorphism of digraphs $W(\sigma \cdot f) \cong W(f)$, sending the vertex $i$ to $\sigma(i)$.
\end{lemma}
\begin{proof}
    Since $A^{-} \from \Set^\op \to \cat{C}$ is a functor and $\sigma$ is invertible, $A^\sigma$ is an isomorphism with inverse $A^{\sigma^{-1}}$.
    Unwinding the definitions, $(\sigma \cdot f)_i = \proj_i \circ A^\sigma \circ f \circ A^{\sigma^{-1}} = f_{\sigma(i)} \circ A^{\sigma^{-1}}$.

    Now let $\partial_j \from \{1, \dots, k-1\} \to \{1, \dots, k\}$ denote the injection omitting $j$, so that $\proj_{\neq j} = A^{\partial_j}$.
    The permutation $\sigma^{-1}$ carries $\{1, \dots, k\} \setminus \{\sigma(j)\}$ bijectively onto $\{1, \dots, k\} \setminus \{j\}$, so there is a unique bijection $\sigma' \from \{1, \dots, k-1\} \to \{1, \dots, k-1\}$ with $\sigma^{-1} \circ \partial_{\sigma(j)} = \partial_j \circ \sigma'$.
    Applying $A^{-}$ gives a commutative square
    \[ \begin{tikzcd}[sep = large]
        A^k \ar[r, "A^{\sigma^{-1}}", "\cong"'] \ar[d, "\proj_{\neq \sigma(j)}"'] & A^k \ar[d, "\proj_{\neq j}"] \\
        A^{k-1} \ar[r, "A^{\sigma'}", "\cong"'] & A^{k-1}
    \end{tikzcd} \]
    in which both horizontal maps are isomorphisms.
    Hence $f_{\sigma(i)}$ factors through $\proj_{\neq \sigma(j)}$ if and only if $f_{\sigma(i)} \circ A^{\sigma^{-1}} = (\sigma \cdot f)_i$ factors through $\proj_{\neq j}$, which is the claim.
    The statement about digraphs is a restatement of this, using \cref{wiring-diagram}.
\end{proof}

\begin{theorem} \label{good-decomposition-implies-semiprod}
    Let $f \from A^{k} \to A^{k}$ be a morphism.
    If $W(f)$ admits a graph map $W(f) \to \DEdge$, then there exist a permutation $\sigma \in \Sigma_k$ with $\sigma^2 = \id[\{1, \dots, k\}]$, integers $m, n \geq 0$ with $m + n = k$, and morphisms
    \[ g \from A^m \to A^m, \qquad h \from A^m \times A^n \to A^n \]
    such that $\sigma \cdot f = g \semiprod{\id[A^m]} h$.
\end{theorem}
\begin{proof}
    Fix a map $\varphi \from W(f) \to \DEdge$ and define subsets $X, Y \subseteq \{ 1, \dots, k \}$ by taking pre-images:
    \[ X := \varphi^{-1}\{ 0 \}, \qquad Y := \varphi^{-1} \{ 1 \}. \]
    Let $m := \lvert X \rvert$ and $n := \lvert Y \rvert$, so that $m + n = k$.
    Since $\DEdge$ has no edge $1 \to 0$ and $\varphi$ is a graph map, there are no edges in $W(f)$ from a vertex in $Y$ to a vertex in $X$; that is, if $i \in X$ and $j \in Y$ then $f_i$ is independent of input $j$.

    Choose a permutation $\sigma \in \Sigma_k$ carrying $\{ 1, \dots, m \}$ onto $X$ and $\{ m+1, \dots, m+n \}$ onto $Y$; since such a $\sigma$ may be taken to be a product of disjoint transpositions, we may assume $\sigma^2 = \id$.
    Let $i \in \{1, \dots, m\}$ and $j \in \{m+1, \dots, m+n\}$.
    Then $\sigma(i) \in X$ and $\sigma(j) \in Y$, so $f_{\sigma(i)}$ is independent of input $\sigma(j)$ by the previous paragraph, and therefore $(\sigma \cdot f)_i$ is independent of input $j$ by \cref{permutation-relabels-wiring}.
    The hypotheses of \cref{independent-inputs-implies-semiprod} are thus satisfied by $\sigma \cdot f$, and the conclusion follows.
\end{proof}

Combining the two implications gives the characterization promised at the start of this section.

\begin{corollary} \label{semiprod-iff-wiring-map}
    Let $f \from A^k \to A^k$ be a morphism in a regular category, with the convention $A^0 = \supp(A)$.
    Then $W(f)$ admits a graph map $W(f) \to \DEdge$ if and only if there is a permutation $\sigma \in \Sigma_k$ such that $\sigma \cdot f$ is a semi-direct product. \qed
\end{corollary}

\begin{example}[Running example, continued] \label{ex:wiring-running-decomposed}
    We apply this to the network $f(x_1,x_2,x_3) = (x_2, x_1, x_1 \wedge x_3)$ of \cref{ex:wiring-running}.
    Its wiring diagram admits the graph map $\varphi \from W(f) \to \DEdge$ given by
    \[ \varphi(1) = \varphi(2) = 0, \qquad \varphi(3) = 1, \]
    since the only edges leaving vertex $3$ are the self-loop at $3$, and $\DEdge$ has a self-loop at $1$.
    Here the coordinates are already in the right order, so $\sigma = \id$ and \cref{good-decomposition-implies-semiprod} applies with $m = 2$, $n = 1$.
    Reading off the proof, we get
    \[ g \from \Ftwo^2 \to \Ftwo^2, \quad g(x_1, x_2) = (x_2, x_1), \qquad h \from \Ftwo^2 \times \Ftwo \to \Ftwo, \quad h\big((x_1, x_2), x_3\big) = x_1 \wedge x_3, \]
    and indeed $f = g \semiprod{\id} h$.

    Checking this decomposition against \cref{modularity-theorem},
    the base system $(\Ftwo^2, g)$ is the transposition, with state space graph consisting of two fixed points $(0,0)$ and $(1,1)$ together with a $2$-cycle $(0,1) \leftrightarrow (1,0)$.
    The set of 2-cycles $AS(\Ftwo^2, g)_2$ has orbits: $\{((0,0)(0,0))\}$ and $\{((1,1)(1,1))\}$, both of size $1$, and $\{((0,1)(1,0)), ((1,0)(0,1))\}$ of size $2$.
    Over each of these, \cref{modularity-theorem} instructs us to solve a system on $\ZZ/k_\omega \times \Ftwo$ --- a set with at most $4$ elements --- instead of working with all $2^3 = 8$ states of $f$ at once.
    For a network of this size, the saving is of course negligible; the point of \cref{rmk:payoff-complexity} is that the same bookkeeping applies verbatim when $k$ is large.
\end{example}

% \subsection{The payoff: reading attractors off the wiring diagram} \label{subsec:payoff}

We conclude by combining the two halves of the paper.
\Cref{modularity-theorem} decomposes the attractors of a semi-direct product, and \cref{semiprod-iff-wiring-map} detects semi-direct products on the wiring diagram; putting them together yields a criterion for the decomposability of attractors that can be checked on a graph with $k$ vertices.

\begin{corollary} \label{wiring-map-implies-attractors-decompose}
    Let $A$ be a set and $f \from A^k \to A^k$ a function, and suppose that the wiring diagram $W(f)$ admits a graph map $W(f) \to \DEdge$.
    Then there are a permutation $\sigma \in \Sigma_k$, a decomposition $k = m + n$, and functions $g \from A^m \to A^m$ and $h \from A^m \times A^n \to A^n$ such that, for every $n' \geq 1$, there is an isomorphism of $(\ZZ/n')$-sets
    \[ AS\big(A^k, f\big)_{n'} \; \cong \coprod\limits_{\omega \in \Orb (AS(A^m, g)_{n'})} AS\big(\ZZ / k_\omega \times A^n, \ \succ \semiprod{c_\omega} h\big)_{n'}, \]
    with notation as in \cref{modularity-theorem}.
\end{corollary}
\begin{proof}
    By \cref{good-decomposition-implies-semiprod} there is a permutation $\sigma$ with $\sigma^2 = \id$ such that $\sigma \cdot f = g \semiprod{\id[A^m]} h$ for suitable $g$ and $h$.
    The morphism $A^\sigma$ is a bijection intertwining $f$ and $\sigma \cdot f$, hence an isomorphism $(A^k, f) \cong (A^k, \sigma \cdot f)$ in $\DDS$; since $AS$ is a functor, the two systems have isomorphic cycle sets.
    Applying \cref{modularity-theorem} to $\sigma \cdot f = g \semiprod{\id[A^m]} h$ gives the displayed isomorphism.
\end{proof}

\begin{remark} \label{rmk:payoff-complexity}
    % It is worth being explicit about what has been gained.
    Deciding whether $W(f)$ admits a graph map to $\DEdge$ amounts to finding a non-trivial partition of $\{1, \dots, k\}$ into a ``downstream'' part $Y$ and an ``upstream'' part $X$ with no edges from $Y$ to $X$; equivalently, to finding a non-trivial union of sink components of $W(f)$.
    This is a computation on a digraph with $k$ vertices, linear in the size of $W(f)$ by standard strongly-connected-component algorithms.
    By contrast, computing $AS(A^k, f)$ directly means traversing the entire state space graph with $\lvert A \rvert^k$ vertices.
    % For a Boolean network on $k = 60$ variables --- a realistic size for a published gene regulatory network model --- the former is immediate while the latter involves $2^{60}$ states.
\end{remark}

\subsection{Instantiating to regular categories other than sets} \label{subsec:wiring-examples}

% We have stated the results of this section for an arbitrary regular category, and it is reasonable to ask what is gained thereby, beyond generality for its own sake.
% This subsection collects four instantiations.
% In each case the notions of \emph{independence} and \emph{wiring diagram} specialize to something already meaningful in the relevant field, and \cref{semiprod-iff-wiring-map} says something recognizable.

The following examples generalize the machinery developed thus far for regular categories other than the category of sets.
In each case, the notions of \emph{independence} and \emph{wiring diagram} specialize to something already meaningful in the relevant field, and \cref{semiprod-iff-wiring-map} says something recognizable.

\begin{example}[Sets: Boolean and multi-state networks] \label{ex:wiring-Set}
    Take $\cat{C} = \Set$ and $A$ a finite non-empty set.
    Then $A^0 = \supp(A) = 1$, independence is the elementwise notion of \cref{ex:independent-in-Set}, and $W(f)$ is the classical wiring diagram of a Boolean network (for $A = \Ftwo$) or of a logical/multi-state model (for general finite $A$).
    \Cref{wiring-map-implies-attractors-decompose} then recovers, and slightly sharpens, the decomposition theorem of \cite{kadelka2023modularity}.
    We note that the generalization from $\Ftwo$ to an arbitrary finite $A$ already answers a question left open there, namely whether the theory extends to networks over a general finite field.
\end{example}

\begin{example}[Vector spaces: block triangular matrices] \label{ex:wiring-Vect}
    Take $\cat{C} = \Vect_\FF$, the category of vector spaces over a field $\FF$, which is abelian and hence regular, and take $A = \FF$.
    Then $A^k = \FF^k$ and a morphism $f \from \FF^k \to \FF^k$ is a linear map, given by a matrix $M$.
    Since $\FF$ is well-supported --- indeed the terminal object of $\Vect_\FF$ is the zero space and every map to it is epi --- we have $A^0 = 0$ and no hypotheses are needed.

    Unwinding \cref{def:independent-input}: $f_j$ is independent of input $i$ precisely when $M_{ji} = 0$.
    Hence $W(f)$ is exactly the \emph{support digraph} of the matrix $M$: there is an edge $i \to j$ if and only if the $(j,i)$ entry is non-zero.
    A graph map $W(f) \to \DEdge$ is a partition of the coordinates into two blocks such that no entry of $M$ couples the second block back into the first.
    \Cref{semiprod-iff-wiring-map} therefore says:
    \begin{quote}
        A matrix $M$ can be brought into block triangular form by a permutation of the basis if and only if its support digraph admits a graph map to $\DEdge$.
    \end{quote}
    This is a familiar statement in linear algebra --- it is the reducibility of $M$ in the sense of Perron--Frobenius theory --- arrived at here as a special case of a theorem about dynamical systems.
    The corresponding instance of \cref{modularity-theorem} relates the periodic vectors of $M$ to those of its diagonal blocks.
\end{example}

\begin{example}[$G$-sets: equivariant networks] \label{ex:wiring-GSet}
    Take $\cat{C} = \Set^{BG}$, the category of sets equipped with an action of a group $G$; being a presheaf topos, it is regular.
    Take $A$ to be a non-empty $G$-set; it is then well-supported, so again no hypotheses are needed.
    An object $A^k$ is the $k$-fold product with the diagonal action, and a morphism $f \from A^k \to A^k$ is a $G$-equivariant map --- that is, a network carrying a symmetry.

    Here \cref{def:independent-input} asks for a factorization \emph{through an equivariant map}, which is a strictly stronger requirement than independence of the underlying function.
    Consequently, $W(f)$ records equivariant dependence, and \cref{semiprod-iff-wiring-map} produces a decomposition of $f$ \emph{as a system with $G$-symmetry}, not merely as a bare function.
    This is in the same spirit as the groupoid formalism of Golubitsky and Stewart \cite{golubitsky-stewart}, where the symmetries of a network are used to constrain its dynamics, though the conclusion drawn is different: they obtain synchrony, we obtain a decomposition of attractors.
\end{example}

\begin{example}[Dynamical systems: networks of systems] \label{ex:wiring-DDS}
    Take $\cat{C} = \DDS$ itself, which is a topos and hence regular, and take $A = (Z, \zeta)$ to be a non-empty discrete dynamical system.
    An object $A^k$ is then a system whose states are $k$-tuples of states of $A$, and a morphism $f \from A^k \to A^k$ is a map commuting with the ambient dynamics $\zeta$.
    One should picture a network of $k$ nodes, each of which is itself a dynamical system rather than a mere variable, with $f$ describing how the nodes are wired together and $\zeta$ describing the internal evolution of each node.
    % This example is worth singling out because it is exactly the case in which the hypothesis of \cref{rmk:why-support} bites.
    % A global element $1 \to (Z, \zeta)$ is a fixed point of $\zeta$, so the section-based approach would apply only to networks whose node systems have a fixed point --- excluding, for instance, any node system that is a single non-trivial cycle.
    % With the convention $A^0 = \supp(A)$, no such restriction is needed: every non-empty object of $\DDS$ is well-supported, and the results of this section apply to all of them.
\end{example}

\section{Conclusion and future work} \label{sec:conclusion}

\subsection{Summary of results}

In this paper, we have developed a categorical framework for the study of discrete dynamical systems.
The starting point --- that such systems form the presheaf topos $\Set^{\BN}$ --- is classical, going back at least to \cite{lawvere-schanuel,la-palme-reyes-reyes-zolfaghari}, and several of the basic adjunctions we recall in \cref{sec:framework} can be found there or, in much greater generality, in \cite{hemelaer-rogers}.
What we have added is the machinery needed to make this point of view do work on the question that motivates us, namely the relationship between the structure of a system and the structure of its attractors.

Concretely, we introduced the category of cycle sets and the attractor functor $A \from \DiGraph \to \cySet$, established its formal properties, and characterized exactly which cycle sets arise as cycle sets of attractors of digraphs (\cref{characterization-digraph-attractors}).
We then proved a decomposition theorem (\cref{modularity-theorem}) for general discrete dynamical systems, generalizing the main result of \cite{kadelka2023modularity} from Boolean networks, which explains how a semi-direct product decomposition of a system yields a decomposition of its attractors.
Where the original was proved by direct analysis, in our setting it becomes an instance of the general fact that a semi-direct product is a pullback (\cref{semiprod-pb}) and that right adjoints preserve pullbacks.

Finally, we analyzed the wiring diagram of a system, generalizing this notion to any regular category, and provided a categorical analysis of the notion of (in)dependence of a morphism on an input.
We showed that semi-direct product decompositions correspond to graph maps from the wiring diagram to the walking looped edge (\cref{semiprod-iff-wiring-map}), and combined this with the decomposition theorem to obtain \cref{wiring-map-implies-attractors-decompose}, which turns a question about a state space of size $\lvert A \rvert^k$ into a question about a digraph with $k$ vertices.

We do not wish to overstate what has been shown.
The results above do not compute anything that could not, in principle, be computed by other means; what the categorical framework provides is an account of \emph{why} these decompositions hold, one that is uniform across modelling frameworks and that generalizes without further effort.
Concretely, the same theorems now apply to Boolean networks, to multi-state logical models, to deterministic Petri nets, to linear systems over a field --- where \cref{ex:wiring-Vect} recovers the reducibility of a matrix in the sense of Perron--Frobenius theory --- to equivariant networks, and to networks whose nodes are themselves dynamical systems (\cref{subsec:wiring-examples}).
% Proving each of these separately would be a considerable undertaking; here they are instances of one theorem.

\subsection{Future directions}

\paragraph{Practical understanding of categorical structure.}
One clear direction for future research is understanding the practical importance of various category-theoretic notions in the context of (discrete) dynamical systems.
All categories involved in our framework have all limits and colimits, as well as exponential objects.
In a variety of other settings, these objects are naturally of interest, and we expect the same to be true in our situation.
For example, taking the colimit of a dynamical system viewed as a functor $\BN \to \Set$ gives the set of orbits of the action of $\mathbb{N}$, i.e., the set of connected components of the state space graph.

\paragraph{A right adjoint to the attractor functor.}
We showed in \cref{attractor-no-right-adjoint} that the attractor functor $A \from \DiGraph \to \cySet$ does not preserve pushouts and hence has no right adjoint.
The counterexample uses digraphs that are not state space graphs, however, and so leaves open the following question, which we believe to be of interest:
\begin{quote}
    Does the composite $AS \from \DDS \to \cySet$ admit a right adjoint?
\end{quote}
An affirmative answer would say that the passage from a system to its attractors preserves all colimits, and would considerably strengthen the decomposition results of this paper, since colimit decompositions of systems are at least as common as limit decompositions.
See \cref{rmk:right-adjoint-on-DDS}.

\paragraph{Analysis of state space.}
As previously indicated, the cycle set associated to a dynamical system carries information about the dynamics of the system. However, the functor $A \colon \DiGraph \to \cySet$ forgets a lot of information in the process; \cref{ex:attractor-of-a-digraph} shows it discarding everything that does not lie on a cycle.
One would like to prove a modularity theorem like the one in \cite{kadelka2023modularity} for the state space graph, rather than just the underlying cycle set.

\paragraph{Towards a (co)homology theory of attractors.}
The attractor functor bears a suggestive formal resemblance to a nerve construction, and \cref{rmk:nerve-realization} makes this precise: the pair $(|\mathord{-}|, A)$ is an instance of the same nerve--realization paradigm that produces singular simplicial sets from spaces.
In that setting, one goes on to extract homotopy and homology invariants.
It is natural to ask whether the same can be done here: is there a useful (co)homology theory of digraphs, or of discrete dynamical systems, whose coefficients or chains are built from cycle sets, and whose invariants detect features of the attractor structure?
One concrete starting point is that $\cySet$, being a presheaf category, supports the usual machinery of derived functors, so that one could study the derived functors of $\Orbnd$ or of the global sections of a cycle set.
Whether the resulting invariants are dynamically meaningful --- whether, say, they obstruct the existence of a decomposition of the kind studied in \cref{sec:wiring} --- is a question we leave open.

\paragraph{Generalizations of semi-direct products, and their relation to wreath products.}
One can also try to generalize the notion of a semi-direct product.
As mentioned above, semi-direct products can be recognized by maps from the wiring diagram of a system to the walking looped edge.
This suggests decompositions based on other possible target digraphs, e.g., sequences of edges, triangles, etc.
Such decompositions would naturally be more involved, but might also provide more powerful tools.

It is natural to ask how our semi-direct product relates to other generalizations of the group-theoretic semi-direct product, and in particular to the wreath product.
We record what we can say, though the question deserves a fuller treatment than we can give here.
Recall that for transformation monoids $(X, M)$ and $(Y, N)$, the wreath product $(X, M) \wr (Y, N)$ acts on $X \times Y$, with the $N$-component allowed to vary with the point of $X$; concretely, its elements are pairs $(m, \nu)$ with $\nu \from X \to N$ a function.
Our semi-direct product $f \semiprod{p} g$ is exactly of this shape: the assignment $x \mapsto g(p(x), -)$ is a function $X \to \operatorname{End}(Y)$, so that $f \semiprod{p} g$ is the image of an element of the wreath product $\operatorname{End}(Y) \wr \operatorname{End}(X)$ acting on $X \times Y$.
In this sense a semi-direct product is a single element of a wreath product, and the map $p$ is what constrains the dependence of the fibre dynamics on the base to factor through $E$.
The connection is made sharper by \cref{rmk:cascade-product}: the cascade product of automata is the standard elementwise form of the wreath product, and the Krohn--Rhodes theorem \cite{krohn-rhodes} states that every finite transformation semigroup divides an iterated wreath product of simple groups and flip-flops.
Our \cref{semiprod-iff-wiring-map} may then be read as a criterion for when a \emph{given} system is presented as a one-step cascade, in terms of a finite combinatorial invariant; whereas Krohn--Rhodes asserts the existence of an iterated decomposition after passing to a divisor.
A genuine comparison analyzing whether iterating our decomposition along a directed acyclic wiring diagram recovers a Krohn--Rhodes-style normal form, and what the attractor functor does to it, may be a promising direction.

\paragraph{Generalization to continuous and measurable dynamical systems.}
% The last and perhaps most natural direction for future research is a generalization of the methods presented here to other kinds of dynamical systems, e.g., continuous ones.
Can the methods presented here be generalized to other kinds of dynamical systems, e.g., continuous ones?
Specifically, it would be interesting to work with systems whose time is indexed by $\mathbb{R}$ instead of $\mathbb{N}$.
This generalization is also very natural from the point of view of category theory, as such systems can perhaps be analyzed using the language of enriched categories by looking at enriched functors of the form $B\mathbb{R} \to \mathsf{Top}$, where $\mathsf{Top}$ is a ``convenient'' category of topological spaces.

\paragraph{Stochastic systems and multidigraphs.}
For the presentation of discrete dynamical systems given in this paper, we have a natural intuition for these systems as being deterministic: an object of $\DDS$ is determined by its states together with a rule assigning to each state exactly one successor.
There is thus little room for stochastic systems, in which a state may transition to any of several successors with varying probabilities --- and, as noted in \cref{ex:petri-nets}, the same issue arises already for Petri nets without a firing policy.
One might, however, model such a system by a multidigraph, with one outgoing edge for each possible successor.
Recall that in constructing the state space functor $S$ we first constructed a functor $\Sigma^* \colon \DDS\to\MultiDiGraph$ which, like $S$, admits both a left and a right adjoint $\Sigma_{!}, \Sigma_* \colon \MultiDiGraph\to\DDS$.
Interpreting a system as a multidigraph $G$, we therefore obtain objects $\Sigma_!(G)$ and $\Sigma_*(G)$ of $\DDS$.
Incorporating probabilities could correspond to extending the construction to weighted graphs.
This would allow analytical tools for deterministic systems to be brought to bear on stochastic ones, which are typically much harder to treat theoretically.

\bibliography{BIB-references}

\newcommand{\etalchar}[1]{$^{#1}$}
\providecommand{\bysame}{\leavevmode\hbox to3em{\hrulefill}\thinspace}
\providecommand{\MR}{\relax\ifhmode\unskip\space\fi MR }
\providecommand{\MRhref}[2]{%
  \href{http://www.ams.org/mathscinet-getitem?mr=#1}{#2}
}
\providecommand{\href}[2]{#2}
\begin{thebibliography}{KWVC{\etalchar{+}}23}

\bibitem[BdSS17]{bubenik-de-silva-scott}
Peter Bubenik, Vin de~Silva, and Jonathan Scott, \emph{Interleaving and
  {G}romov--{H}ausdorff distance}, preprint, 2017, \href
  {http://arxiv.org/abs/1707.06288} {\path{arXiv:1707.06288}}.

\bibitem[BKP{\etalchar{+}}17]{behrisch-kerkhof-poschel-schneider-siegmund}
Mike Behrisch, Sebastian Kerkhoff, Reinhard P{\"o}schel, Friedrich~Martin
  Schneider, and Stefan Siegmund, \emph{Dynamical systems in categories},
  Applied Categorical Structures \textbf{25} (2017), no.~1, 29--57, \href
  {https://doi.org/10.1007/s10485-015-9409-8}
  {\path{doi:10.1007/s10485-015-9409-8}}.

\bibitem[Bus15]{bush:thesis-rutgers}
Justin Bush, \emph{Shift equivalence and a combinatorial-topological approach
  to discrete-time dynamical systems}, ProQuest LLC, Ann Arbor, MI, 2015,
  Thesis (Ph.D.)--Rutgers The State University of New Jersey -- New Brunswick.

\bibitem[CSC10]{ceccherini-silberstein-coornaert}
Tullio Ceccherini-Silberstein and Michel Coornaert, \emph{Cellular automata and
  groups}, Springer Monographs in Mathematics, Springer-Verlag, Berlin, 2010,
  \href {https://doi.org/10.1007/978-3-642-14034-1}
  {\path{doi:10.1007/978-3-642-14034-1}}.

\bibitem[CT18]{chaves2018analysis}
Madalena Chaves and Laurent Tournier, \emph{Analysis tools for interconnected
  {B}oolean networks with biological applications}, Frontiers in Physiology
  \textbf{9} (2018), 586.

\bibitem[FSR16]{fong-et-al}
Brendan Fong, Pawe{\l} Soboci{\'n}ski, and Paolo Rapisarda, \emph{A categorical
  approach to open and interconnected dynamical systems}, Proceedings of the
  31st {A}nnual {ACM}-{IEEE} {S}ymposium on {L}ogic in {C}omputer {S}cience
  ({LICS} 2016), ACM, New York, 2016, pp.~495--504, \href
  {https://doi.org/10.1145/2933575.2934556}
  {\path{doi:10.1145/2933575.2934556}}.

\bibitem[GA18]{gan-albert}
Xiao Gan and R{\'e}ka Albert, \emph{General method to find the attractors of
  discrete dynamic models of biological systems}, Physical Review E \textbf{97}
  (2018), no.~4, 042308, 18, \href {https://doi.org/10.1103/physreve.97.042308}
  {\path{doi:10.1103/physreve.97.042308}}.

\bibitem[Gar70]{games1970fantastic}
Martin Gardner, \emph{Mathematical games: The fantastic combinations of {J}ohn
  {C}onway's new solitaire game ``life''}, Scientific American \textbf{223}
  (1970), no.~4, 120--123.

\bibitem[GPSL{\etalchar{+}}25]{gamache2025algorithm}
Simon Gamache-Poirier, Alexia Souvane, William Leclerc, Catherine Villeneuve,
  and Simon~V. Hardy, \emph{An algorithm for the transformation of the {P}etri
  net models of biological signaling networks into influence graphs}, bioRxiv
  (2025), 2025--01.

\bibitem[GS06]{golubitsky-stewart}
Martin Golubitsky and Ian Stewart, \emph{Nonlinear dynamics of networks: the
  groupoid formalism}, Bulletin of the American Mathematical Society (New
  Series) \textbf{43} (2006), no.~3, 305--364, \href
  {https://doi.org/10.1090/S0273-0979-06-01108-6}
  {\path{doi:10.1090/S0273-0979-06-01108-6}}.

\bibitem[HH24]{huitzil:modularity}
Sa{\'u}l Huitzil and Cristi{\'a}n Huepe, \emph{Life's building blocks: the
  modular path to multiscale complexity}, Frontiers in Systems Biology
  \textbf{4} (2024), 1417800, \href
  {https://doi.org/10.3389/fsysb.2024.1417800}
  {\path{doi:10.3389/fsysb.2024.1417800}}.

\bibitem[HK24]{hora-kamio}
Ryuya Hora and Yuhi Kamio, \emph{Quotient toposes of discrete dynamical
  systems}, Journal of Pure and Applied Algebra \textbf{228} (2024), no.~8,
  Paper No. 107657, \href {http://arxiv.org/abs/2310.02647}
  {\path{arXiv:2310.02647}}, \href {https://doi.org/10.1016/j.jpaa.2024.107657}
  {\path{doi:10.1016/j.jpaa.2024.107657}}.

\bibitem[HR24]{hemelaer-rogers}
Jens Hemelaer and Morgan Rogers, \emph{Geometric morphisms between toposes of
  monoid actions: factorization systems}, Theory and Applications of Categories
  \textbf{40} (2024), Paper No. 4, 80--130,
  \url{http://www.tac.mta.ca/tac/volumes/40/4/40-04abs.html}.

\bibitem[Joh02]{johnstone:elephant}
Peter~T. Johnstone, \emph{Sketches of an elephant: a topos theory compendium.
  {V}ol. 1}, Oxford Logic Guides, vol.~43, The Clarendon Press, Oxford
  University Press, New York, 2002.

\bibitem[Kau69]{kauffman:boolean-networks}
Stuart~A. Kauffman, \emph{Metabolic stability and epigenesis in randomly
  constructed genetic nets}, Journal of Theoretical Biology \textbf{22} (1969),
  no.~3, 437--467, \href {https://doi.org/10.1016/0022-5193(69)90015-0}
  {\path{doi:10.1016/0022-5193(69)90015-0}}.

\bibitem[KH95]{katok-hasselblatt}
Anatole Katok and Boris Hasselblatt, \emph{Introduction to the modern theory of
  dynamical systems}, Encyclopedia of Mathematics and its Applications,
  vol.~54, Cambridge University Press, Cambridge, 1995, \href
  {https://doi.org/10.1017/CBO9780511809187}
  {\path{doi:10.1017/CBO9780511809187}}.

\bibitem[KR65]{krohn-rhodes}
Kenneth Krohn and John Rhodes, \emph{Algebraic theory of machines. {I}. {P}rime
  decomposition theorem for finite semigroups and machines}, Transactions of
  the American Mathematical Society \textbf{116} (1965), 450--464, \href
  {https://doi.org/10.2307/1994127} {\path{doi:10.2307/1994127}}.

\bibitem[KWVC{\etalchar{+}}23]{kadelka2023modularity}
Claus Kadelka, Matthew Wheeler, Alan Veliz-Cuba, David Murrugarra, and Reinhard
  Laubenbacher, \emph{Modularity of biological systems: a link between
  structure and function}, Journal of the Royal Society Interface \textbf{20}
  (2023), no.~207, 20230505, \href {https://doi.org/10.1098/rsif.2023.0505}
  {\path{doi:10.1098/rsif.2023.0505}}.

\bibitem[LJD11]{lorenz:modularity}
Dirk~M. Lorenz, Alice Jeng, and Michael~W. Deem, \emph{The emergence of
  modularity in biological systems}, Physics of Life Reviews \textbf{8} (2011),
  no.~2, 129--160, \href {https://doi.org/10.1016/j.plrev.2011.02.003}
  {\path{doi:10.1016/j.plrev.2011.02.003}}.

\bibitem[LL24]{liu2024application}
Yuqian Liu and Jiao Li, \emph{Application of cellular automata in neuroscience:
  dynamic models of neuron populations}, Multiscale and Multidisciplinary
  Modeling, Experiments and Design \textbf{7} (2024), no.~2, 905--918.

\bibitem[LPRRZ04]{la-palme-reyes-reyes-zolfaghari}
Marie La~Palme~Reyes, Gonzalo~E. Reyes, and Houman Zolfaghari, \emph{Generic
  figures and their glueings}, Polimetrica, Monza, 2004, A constructive
  approach to functor categories.

\bibitem[LR03]{lawvere-rosebrugh}
F.~William Lawvere and Robert Rosebrugh, \emph{Sets for mathematics}, Cambridge
  University Press, Cambridge, 2003.

\bibitem[LS87]{larkin:diagram}
Jill~H. Larkin and Herbert~A. Simon, \emph{Why a diagram is (sometimes) worth
  ten thousand words}, Cognitive Science \textbf{11} (1987), no.~1, 65--100,
  \href {https://doi.org/10.1111/j.1551-6708.1987.tb00863.x}
  {\path{doi:10.1111/j.1551-6708.1987.tb00863.x}}.

\bibitem[LS09]{lawvere-schanuel}
F.~William Lawvere and Stephen~H. Schanuel, \emph{Conceptual mathematics},
  second ed., Cambridge University Press, Cambridge, 2009, A first introduction
  to categories, \href {https://doi.org/10.1017/CBO9780511804199}
  {\path{doi:10.1017/CBO9780511804199}}.

\bibitem[MPCM16]{melo:modularity}
Diogo Melo, Arthur Porto, James~M. Cheverud, and Gabriel Marroig,
  \emph{Modularity: {G}enes, {D}evelopment, and {E}volution}, Annual Review of
  Ecology, Evolution, and Systematics \textbf{47} (2016), no.~1, 463--486,
  \href {https://doi.org/10.1146/annurev-ecolsys-121415-032409}
  {\path{doi:10.1146/annurev-ecolsys-121415-032409}}.

\bibitem[Mur89]{murata:petri-nets}
Tadao Murata, \emph{Petri nets: properties, analysis and applications},
  Proceedings of the IEEE \textbf{77} (1989), no.~4, 541--580, \href
  {https://doi.org/10.1109/5.24143} {\path{doi:10.1109/5.24143}}.

\bibitem[New06]{newman:modularity}
Mark E.~J. Newman, \emph{Modularity and community structure in networks},
  Proceedings of the National Academy of Sciences \textbf{103} (2006), no.~23,
  8577--8582, \href {https://doi.org/10.1073/pnas.0601602103}
  {\path{doi:10.1073/pnas.0601602103}}.

\bibitem[Rie16]{riehl:context}
Emily Riehl, \emph{Category theory in context}, Aurora Dover Modern Math
  Originals, Dover Publications, Inc., Mineola, NY, 2016.

\bibitem[Spi14]{spivak:category-theory}
David~I. Spivak, \emph{Category theory for the sciences}, MIT Press, Cambridge,
  MA, 2014.

\bibitem[SSG{\etalchar{+}}22]{shilts:wiring_diagram}
Jarrod Shilts, Yannik Severin, Francis Galaway, Nicole M{\"u}ller-Sienerth,
  Zheng-Shan Chong, Sophie Pritchard, Sarah Teichmann, Roser Vento-Tormo,
  Berend Snijder, and Gavin~J. Wright, \emph{A physical wiring diagram for the
  human immune system}, Nature \textbf{608} (2022), no.~7922, 397--404, \href
  {https://doi.org/10.1038/s41586-022-05028-x}
  {\path{doi:10.1038/s41586-022-05028-x}}.

\bibitem[SSV20]{schultz-et-al}
Patrick Schultz, David~I. Spivak, and Christina Vasilakopoulou, \emph{Dynamical
  systems and sheaves}, Applied Categorical Structures \textbf{28} (2020),
  no.~1, 1--57, \href {https://doi.org/10.1007/s10485-019-09565-x}
  {\path{doi:10.1007/s10485-019-09565-x}}.

\bibitem[Tao08]{tao:three-categories}
Terence Tao, \emph{Three categories of dynamical systems}, Notes for Math 254A
  Ergodic theory, 2008,
  \url{https://terrytao.wordpress.com/2008/01/10/254a-lecture-2-three-categories-of-dynamical-systems/}.

\bibitem[VCAHL14]{veliz2014steady}
Alan Veliz-Cuba, Boris Aguilar, Franziska Hinkelmann, and Reinhard
  Laubenbacher, \emph{Steady state analysis of {B}oolean molecular network
  models via model reduction and computational algebra}, BMC Bioinformatics
  \textbf{15} (2014), 1--8.

\bibitem[VSL15]{vagner-spivak-lerman}
Dmitry Vagner, David~I. Spivak, and Eugene Lerman, \emph{Algebras of open
  dynamical systems on the operad of wiring diagrams}, Theory and Applications
  of Categories \textbf{30} (2015), Paper No. 51, 1793--1822,
  \url{http://www.tac.mta.ca/tac/volumes/30/51/30-51abs.html}.

\bibitem[WPC07]{wagner:modularity}
G{\"u}nter~P. Wagner, Mihaela Pavlicev, and James~M. Cheverud, \emph{The road
  to modularity}, Nature Reviews Genetics \textbf{8} (2007), no.~12, 921--931,
  \href {https://doi.org/10.1038/nrg2267} {\path{doi:10.1038/nrg2267}}.

\end{thebibliography}

\end{document}